\documentclass[moor]{informs1}              

\usepackage[english]{babel}
\usepackage{amsmath,bbm}
\usepackage{graphicx}
\usepackage{amsfonts}
\usepackage{caption}
\usepackage{subcaption}
\usepackage{todonotes}
\usepackage{mathtools}
\usepackage{float}
\usepackage{natbib}
\usepackage[utf8]{inputenc}
\usepackage[title]{appendix}
\usepackage{etoolbox}
\usepackage{etex}
\usepackage{tikz}

\usepackage{chngpage}

\RequirePackage[normalem]{ulem} 
\RequirePackage{color}\definecolor{RED}{rgb}{1,0,0}\definecolor{BLUE}{rgb}{0,0,1} 

\usetikzlibrary{trees}
\usetikzlibrary{chains,shapes.multipart}
\usetikzlibrary{shapes}
\usetikzlibrary{automata,positioning}

\theoremstyle{definition}

\theoremstyle{remark}

\numberwithin{equation}{section}

 \renewcommand \newline {\\ \indent}
 
\let\ForAll\forall
\renewcommand{\forall}[0]{\ForAll\;}

\let\Exists\exists
\renewcommand{\exists}[0]{\Exists\;}

\makeatletter
\newcommand{\vast}{\bBigg@{4}}
\newcommand{\Vast}{\bBigg@{5}}
\makeatother

\newcommand{\ArrivalProbability}{p^{(N)}}
\newcommand{\SuccessProbability}{q^{(N)}}
 
\newcommand{\Arrival}[1]{\ifstrequal{#1}{}{A^{(N)}}{A^{(N)}\left(#1\right)}}

\newcommand{\TildeArrival}[1]{\tilde{A}^{(N)}\left(#1\right)}
\newcommand{\ArrivalBernoulli}[1]{X^{(N)}(#1)}

\newcommand{\ArrivalBernoulliUpper}[1]{X^{(u,N)}(#1)}

\newcommand{\ArrivalUpper}[1]{A^{(u,N)}(#1)}
\newcommand{\Service}[2]{\ifstrequal{#2}{}{S_{#1}^{(N)}}{S_{#1}^{(N)}\left(#2\right)}}
\newcommand{\TildeService}[2]{\tilde{S}_{#1}^{(N)}\left(#2\right)}
\newcommand{\ServiceBernoulli}[2]{Y_{#1}^{(N)}\left(#2\right)}

\newcommand{\ServiceBernoulliUpper}[2]{Y_{#1}^{(u,N)}\left(#2\right)}
\newcommand{\ServiceLower}[2]{\ifstrequal{#2}{}{S_{#1}^{(l,N)}}{S_{#1}^{(l,N)}\left(#2\right)}}

\newcommand{\ServiceUpper}[2]{S_{#1}^{(u,N)}\left(#2\right)}

\newcommand{\TildeQueue}[2]{\ifstrequal{#2}{}{\tilde{R}_{#1}^{(N)}}{\tilde{R}_{#1}^{(N)}\left(#2\right)}}

\newcommand{\ScaledProcess}{\frac{1}{N\sqrt{\log N}}\left(\Arrival{{tN^{3}\log N}}-\Service{i}{{tN^{3}\log N}}\right)}
\newcommand{\ScaledProcessLine}{\frac{\left(\Arrival{{tN^{3}\log N}}-\Service{i}{{tN^{3}\log N}}\right)}{(N\sqrt{\log N})}}

\newcommand{\Calculation}[1]{\the\numexpr 2*#1}

\newcommand{\MaxQueueLength}[1]{Q_{(\alpha,\beta)}^{(N)}\left(#1\right)}

\newcommand{\MaxQueueLengthLine}[1]{Q_{(\alpha,\beta)}^{(N)}\big(#1\big)}

\newcommand{\QueueZero}[1]{Q_{#1}^{(N)}(0)}

\newcommand{\StartPositionGeneral}[2]
{\frac{\Arrival{#2N^3\log N}-\Service{#1}{#2N^3\log N}+\QueueZero{#1}}{N\log N}}

\newcommand{\FluidProcess}[1]{%
	\mathchoice
		{\ifstrequal{#1}{0}
{\frac{\MaxQueueLength{0}}{N\sqrt{\log N}}}
{\ifstrequal{#1}{1}
{\frac{\MaxQueueLength{N^3}}{N\sqrt{\log N}}}{\frac{\MaxQueueLength{#1N^3}}{N\sqrt{\log N}}}}
}
		{\ifstrequal{#1}{0}
{\frac{\MaxQueueLengthLine{0}}{\big(N\sqrt{\log N}\big)}}
{\ifstrequal{#1}{1}
{\frac{\MaxQueueLengthLine{N^3}}{\big(N\sqrt{\log N}\big)}}{\frac{\MaxQueueLengthLine{#1N^3}}{\big(N\sqrt{\log N}\big)}}}
}
		{\ifstrequal{#1}{0}
{\frac{\MaxQueueLengthLine{0}}{\big(N\sqrt{\log N}\big)}}
{\ifstrequal{#1}{1}
{\frac{\MaxQueueLengthLine{N^3}}{\big(N\sqrt{\log N}\big)}}{\frac{\MaxQueueLengthLine{#1N^3}}{\big(N\sqrt{\log N}\big)}}}
}
		{\ifstrequal{#1}{0}
{\frac{\MaxQueueLengthLine{0}}{\big(N\sqrt{\log N}\big)}}
{\ifstrequal{#1}{1}
{\frac{\MaxQueueLengthLine{N^3}}{\big(N\sqrt{\log N}\big)}}{\frac{\MaxQueueLengthLine{#1N^3}}{\big(N\sqrt{\log N}\big)}}}
}
}

\newcommand{\FluidProcessTwo}[1]{%
	\mathchoice
		{\ifstrequal{#1}{0}
{\frac{\MaxQueueLength{0}}{N\log N}}
{\ifstrequal{#1}{1}
{\frac{\MaxQueueLength{N^3\log N}}{N\log N}}
{\ifstrequal{#1}{\infty}
{\frac{\MaxQueueLength{\infty}}{N\log N}}
{\frac{\MaxQueueLength{#1N^3\log N}}{N\log N}}}}
}
		{\ifstrequal{#1}{0}
{\frac{\MaxQueueLengthLine{0}}{\big(N\log N\big)}}
{\ifstrequal{#1}{1}
{\frac{\MaxQueueLengthLine{N^3\log N}}{\big(N\log N\big)}}
{\ifstrequal{#1}{\infty}
{\frac{\MaxQueueLength{\infty}}{\big(N\log N\big)}}
{\MaxQueueLengthLine{#1N^3\log N}\big/\big(N\log N\big)}}}
}
		{\ifstrequal{#1}{0}
{\frac{\MaxQueueLengthLine{0}}{\big(N\log N\big)}}
{\ifstrequal{#1}{1}
{\frac{\MaxQueueLengthLine{N^3\log N}}{\big(N\log N\big)}}
{\ifstrequal{#1}{\infty}
{\frac{\MaxQueueLength{\infty}}{\big(N\log N\big)}}
{\frac{\MaxQueueLengthLine{#1N^3\log N}}{\big(N\log N\big)}}}}
}
		{\ifstrequal{#1}{0}
{\frac{\MaxQueueLengthLine{0}}{\big(N\log N\big)}}
{\ifstrequal{#1}{1}
{\frac{\MaxQueueLengthLine{N^3\log N}}{\big(N\log N\big)}}{\frac{\MaxQueueLengthLine{#1N^3\log N}}{\big(N\log N\big)}}}
}
}

\newcommand{\FluidProcessThree}[1]{%
	\mathchoice
		{\ifstrequal{#1}{0}
{\frac{\MaxQueueLength{0}}{N\log N}}
{\ifstrequal{#1}{1}
{\frac{\MaxQueueLength{N^3\log N}}{N\log N}}
{\ifstrequal{#1}{\infty}
{\frac{\MaxQueueLength{\infty}}{N\log N}}
{\frac{\MaxQueueLength{#1^{(N)}}}{N\log N}}}}
}
		{\ifstrequal{#1}{0}
{\frac{\MaxQueueLengthLine{0}}{\big(N\log N\big)}}
{\ifstrequal{#1}{1}
{\frac{\MaxQueueLengthLine{N^3\log N}}{\big(N\log N\big)}}
{\ifstrequal{#1}{\infty}
{\frac{\MaxQueueLength{\infty}}{\big(N\log N\big)}}
{\frac{\MaxQueueLengthLine{#1^{(N)}}}{\big(N\log N\big)}}}}
}
		{\ifstrequal{#1}{0}
{\frac{\MaxQueueLengthLine{0}}{\big(N\log N\big)}}
{\ifstrequal{#1}{1}
{\frac{\MaxQueueLengthLine{N^3\log N}}{\big(N\log N\big)}}
{\ifstrequal{#1}{\infty}
{\frac{\MaxQueueLength{\infty}}{\big(N\log N\big)}}
{\frac{\MaxQueueLengthLine{#1^{(N)}}}{\big(N\log N\big)}}}}
}
		{\ifstrequal{#1}{0}
{\frac{\MaxQueueLengthLine{0}}{\big(N\log N\big)}}
{\ifstrequal{#1}{1}
{\frac{\MaxQueueLengthLine{N^3\log N}}{\big(N\log N\big)}}{\frac{\MaxQueueLengthLine{#1^{(N)}}}{\big(N\log N\big)}}}
}
}

\newcommand{\ThetaArrivalUpper}{\theta_A^{(u,N)}}

\newcommand{\ThetaServiceUpper}[1]{\theta_{#1}^{(u,N)}}

\newcommand{\MaximumProcess}[2]{M^{(N)}_{#1}(#2)}
\newcommand{\DifferenceNormal}{\frac{c_{t}}{N\sqrt{\log N}}\frac{1}{1+|y|^3}}

\newcommand{\StartRVTriangIndep}[1]{U_{#1}^{(N)}}
\newcommand{\StartRVIndep}[1]{U_{#1}}
\newcommand{\StartRVTriangDep}[1]{V_{#1}^{(N)}}

\newcommand{\LimitN}{\overset{N\to\infty}{\longrightarrow}}
\newcommand{\LimitD}{\overset{d}{\longrightarrow}}
\newcommand{\LimitP}{\overset{\mathbb{P}}{\longrightarrow}}

\DeclarePairedDelimiterX{\abs}[1]\lvert\rvert{\ifblank{#1}{\,\cdot\,}{#1}}
\newcommand{\expect}{\operatorname{\mathbb{E}}\expectarg}

\DeclarePairedDelimiterX{\expectarg}[1]{[}{]}{%
  \ifnum\currentgrouptype=16 \else\begingroup\fi
  \activatebar#1
  \ifnum\currentgrouptype=16 \else\endgroup\fi
}
\newcommand{\probability}{\operatorname{\mathbb{P}}\probarg}
\DeclarePairedDelimiterX{\probarg}[1]{(}{)}{%
  \ifnum\currentgrouptype=16 \else\begingroup\fi
  \activatebar#1
  \ifnum\currentgrouptype=16 \else\endgroup\fi
}
\newcommand{\innermid}{\nonscript\;\delimsize\vert\nonscript\;}
\newcommand{\activatebar}{%
  \begingroup\lccode`\~=`\|
  \lowercase{\endgroup\let~}\innermid
  \mathcode`|=\string"8000
}

\let\oldfrac\frac
\renewcommand{\frac}[2]{%
  \mathchoice
    {\oldfrac{#1}{#2}}
    {#1/#2}
    {#1/#2}
    {#1/#2}
}

\newif\ifhideproofs
\ifhideproofs
  \usepackage{environ}
  \NewEnviron{hide}{}
  \let\proof\hide
  \let\endproof\endhide
\fi

\usepackage{natbib}
 \NatBibNumeric
 \def\bibfont{\small}%
 \def\BIBand{and}%
 \bibpunct[, ]{[}{]}{,}{n}{}{,}%

\usepackage[hypertexnames=false,unicode,colorlinks=true,breaklinks=true,bookmarks=true,urlcolor=blue,
     citecolor=blue,linkcolor=blue,bookmarksopen=false,draft=false,psdextra=true]{hyperref}

\def\EMAIL#1{\href{mailto:#1}{#1}}
\def\URL#1{\href{#1}{#1}}         

\TheoremsNumberedBySection  

\EquationsNumberedBySection 

\begin{document}


 \RUNAUTHOR{Schol, Vlasiou, and Zwart}

 \RUNTITLE{Large fork-join queues with nearly deterministic arrival and service times}

\TITLE{Large fork-join queues with nearly deterministic arrival and service times}

\ARTICLEAUTHORS{%
\AUTHOR{Dennis Schol}
\AFF{Eindhoven University of Technology,
P.O. Box 513, 5600 MB
Eindhoven, The Netherlands, \EMAIL{c.schol@tue.nl} \URL{}}
\AUTHOR{Maria Vlasiou}
\AFF{Eindhoven University of Technology, University of Twente,
P.O. Box 513, 5600 MB
Eindhoven, The Netherlands, \EMAIL{m.vlasiou@tue.nl} \URL{}}
\AUTHOR{Bert Zwart}
\AFF{Eindhoven University of Technology, CWI,
P.O. Box 513, 5600 MB
Eindhoven, The Netherlands, \EMAIL{bert.zwart@cwi.nl} \URL{}}
} 

\ABSTRACT{%
In this paper, we study an $N$ server fork-join queue with nearly deterministic arrival and service times. Specifically, we present a fluid limit for the maximum queue length as $N\to\infty$. This fluid limit depends on the initial number of tasks. In order to prove these results, we develop extreme value theory and diffusion approximations for the queue lengths.
}%


\KEYWORDS{queueing network; heavy traffic; fluid limit; extreme value theory}
\MSCCLASS{60K25}
\ORMSCLASS{Primary: Queues: networks; limit theorems; secondary: probability: Markov processes; random walk; stochastic model applications}
\HISTORY{Received December 24, 2019; revised December 15, 2020; accepted March 23, 2021.}

\maketitle

%
\section[Introduction]{Introduction.}
Fork-join queues are widely studied in many applications, such as communication systems and production processes. However, due to the fact that all service stations see exactly the same arrival process, which is the main characteristic of fork-join queues, these fork-join queues are very challenging to analyze.  Hence, there are only a few exact results, which are mainly for systems in stationarity and are restricted to fork-join queues with two service stations.

In this paper, we focus on a fork-join queue where the number of service stations is large. Our objective is to analyze the queue length of the longest queue. We explore a discrete-time fork-join queue where the arrival and service times are nearly deterministic. In addition, we consider a heavily loaded system. That is, we assume that the arrival rate to a queue times the expected service time of that queue, i.e.\ the traffic intensity per queue $\rho_N$, depends on the number of service stations $N$ and satisfies $(1-\rho_N)N^2\LimitN \beta$, with $\beta>0$. Our main result is a fluid limit of the maximum queue length of the system as $N$ goes to infinity, which holds under very mild conditions on the distribution of the number of jobs at time 0.

Both the model and the scaling studied in this paper are inspired by assembly systems. In particular, we are inspired by problems faced by original equipment manufacturers (OEMs) that assemble thousands of components, each produced using specialized equipment, into complex systems. Examples of such OEMs are Airbus and ASML. If one component is missing, the final product cannot be assembled, giving rise to costly delays. In reality, for some components, OEMs may hedge the shortage risk by investing in capacity or by keeping an inventory of finished components. However, we study the maximum queue length, which is only relevant for components where there is no inventory. As such, our model is a somewhat stylized model of reality.

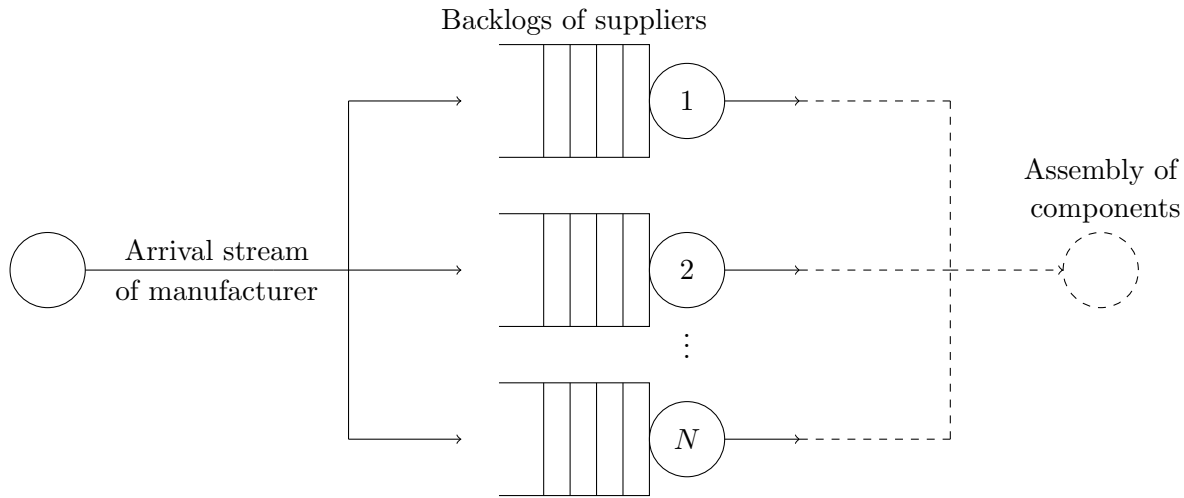
\begin{figure}[H]
\centering
\begin{tikzpicture}[level distance=3cm,
level 1/.style={sibling distance=2cm},
edge from parent fork right]

\node[] (Root) [red]  (A) {}
    child[line width=0,grow=right,draw opacity=0] {
    node[]  {}    
    child[draw=black,draw opacity=0,->] { node{} }     
    child[draw=black,draw opacity=0,->] { node{} }
    child[draw=black,draw opacity=0,->] { node {} }
}
;
\node [draw=none, shift={(8.5cm,-0.9cm)}] (A) {\vdots};
\draw (8.5cm,-2.25cm) circle [radius=0.5cm] node {$N$};
\draw (8.5cm,0cm) circle [radius=0.5cm] node {2};
\draw (8.5cm,2.25cm) circle [radius=0.5cm] node{1};
\draw (6,3) -- ++(2cm,0) -- ++(0,-1.5cm) -- ++(-2cm,0)node[above,yshift=1.5cm,xshift=1cm]{\text{Backlogs of suppliers}};
\foreach \i in {1,...,4}
  \draw (8cm-\i*10pt,3) -- +(0,-1.5cm) ;
  \draw (6,0.75) -- ++(2cm,0) -- ++(0,-1.5cm) -- ++(-2cm,0);
\foreach \i in {1,...,4}
  \draw (8cm-\i*10pt,0.75) -- +(0,-1.5cm);
  
  \draw (6,-1.5) -- ++(2cm,0) -- ++(0,-1.5cm) -- ++(-2cm,0);
\foreach \i in {1,...,4}
  \draw (8cm-\i*10pt,-1.5) -- +(0,-1.5cm);
  \draw[->] (9,-2.25) -- (10,-2.25) node[above]{};
    \draw[->] (9,0) -- (10,0) node[above]{};
      \draw[->] (9,2.25) -- (10,2.25) node[above]{};
    
      \draw (0,0) circle [radius=0.5cm] node{};
      \draw (0.5,0) -- (3,0) node[above,xshift=-0.75cm]{\text{Arrival stream}} node[below,xshift=-0.75cm]{\text{of manufacturer}};
    \draw[->] (3,0) -- (5.5,0);
    \draw (4,2.25)--(4,-2.25);
    \draw[->] (4,2.25) -- (5.5,2.25);
    \draw[->] (4,-2.25) -- (5.5,-2.25);
    
    \draw[dashed](10,0)--(12,0);
    \draw[dashed](10,-2.25)--(12,-2.25);    
    \draw[dashed](10,2.25)--(12,2.25);
    
    \draw[dashed](12,2.25)--(12,0);
    \draw[dashed](12,-2.25)--(12,0);
    
    \draw[dashed,->](12,0)--(13.5,0);
    
     \draw[dashed] (14,0) circle [radius=0.5cm] node[above,yshift=1cm,xshift=0cm]{\text{Assembly of}}
node[above,yshift=0.5cm,xshift=0cm]{\text{ components}};
\end{tikzpicture}
\caption{Fork-join queue with $N$ servers}
\label{fig: tikzpicture fork-join queue}
\end{figure}

An interesting question is whether the manufacturer can produce on schedule. To answer this question, we consider a make-to-order system, i.e.\ suppliers only produce when they have an order, and we assume that the manufacturer sends orders to all the suppliers at the same time. Now, we can model this process by a fork-join queueing system, where the various servers represent suppliers, jobs in the system represent orders requested by the manufacturer and queue lengths in front of each server represent the number of unfinished components each supplier has. As the slowest supplier determines the delay that the manufacturer observes, we wish to study the longest queue. Additionally, we consider a supply chain network operating under full capacity, which is indeed the situation in this industry. Last, we capture the property that in high-tech manufacturing arrival and service times have a low variance by considering nearly deterministic arrival and service times. A visualization of the fork-join queue as a simple representation of a high-tech supply chain is given in Figure \ref{fig: tikzpicture fork-join queue}. Note that in this paper we focus on the backlogs of the suppliers and not on the assembly phase.
We now turn to a survey of related literature. As mentioned, the earliest literature on fork-join queues focuses on systems with two service stations. Analytic results, such as asymptotics on limiting distributions, can be found in \cite{baccelli1985two, flatto1984two, de1988fredholm, wright1992two}. However, due to the complexity of fork-join queues, these results cannot be expanded to fork-join queues with more than two service stations. Thus, most of the work on fork-join queues with more than two service stations is focused on finding approximations of performance measures. For example, an approximation of the distribution of the response time in M/M/s fork-join queues is given in Ko and Serfozo \cite{ko2004response}. Upper and lower bounds for the mean response time of servers, and other performance measures, are given by Nelson, Tantawi \cite{nelson1988approximate} and Baccelli, Makowski \cite{baccelli1989queueing}.

A common property of the aforementioned classic literature is that it mainly focuses on steady-state distributions or other one-dimensional performance measures. Some work on the heavy-traffic process limit has been done, for example, Varma \cite{varma1990heavy} derives a heavy-traffic analysis for fork-join queues, and shows weak convergence of several processes, such as the joint queue lengths in front of each server. Furthermore, Nguyen \cite{nguyen1993processing} proves that various appearing limiting processes are in fact multi-dimensional reflected Brownian motions. In \cite{nguyen1994trouble}, Nguyen extends this result to a fork-join queue with multiple job types. Lu and Pang study fork-join networks in \cite{lu2015gaussian,lu2017heavy,lu2017heavy2}. In \cite{lu2015gaussian}, they investigate a fork-join network where each service station has multiple servers under non-exchangeable synchronization and operates in the quality-driven regime. They derive functional central limit theorems for the number of tasks waiting in the waiting buffers for synchronization and for the number of synchronized jobs. In \cite{lu2017heavy}, they extend this analysis to a fork-join network with a fixed number of service stations, each having many servers, where the system operates in the Halfin-Whitt regime. In \cite{lu2017heavy2}, the authors investigate these heavy-traffic limits for a fixed number of infinite-server stations, where services are dependent and could be disrupted. Finally, we mention Atar, Mandelbaum and Zviran \cite{atar2012control}, who investigate the control of a fork-join queue in heavy traffic by using feedback procedures. Our work contributes to this literature on the process-level analysis of fork-join networks. To be precise, we derive a fluid limit of the stochastic process that keeps track of the largest queue length. This study seems to be the first explicit process-level approximation of a large fork-join queue.

Moreover, our work also adds to the literature on queueing systems with nearly deterministic arrivals and services.
The only research line on queueing systems with nearly deterministic service times that we are aware of is Sigman and Whitt \cite{ sigman2011heavy,sigman2011heavystat},
who investigate the G/G/1 and G/D/N queues and establish heavy-traffic results on waiting times, queue lengths and other performance measures in stationarity, as well as functional central limit theorems on the waiting time and on other performance measures. In these papers, they
distinguish two cases, one in which $(1-\rho_N)\sqrt{N}\LimitN\beta$ and one in which $(1-\rho_N)N\LimitN\beta$, with $\rho_N$
the traffic intensity and $\beta$ some constant.

We now turn to an overview of the techniques that we use in this paper. Because of the fact that we aim to obtain a fluid limit of a maximum of $N$ queue lengths, we mainly use techniques from extreme value theory in our proofs. This is, however, quite a challenge, since, on the one hand, the queue lengths of the servers are mutually dependent. On the other hand, most results on extreme values hinge heavily on the assumption of mutual independence. Furthermore, we consider a fork-join queue where the arrival and service probabilities depend on $N$, which makes the queue lengths triangular arrays with respect to $N$. This makes our paper also rather unusual, as studies on triangular arrays are rare. One paper on this subject, relevant for us, is Anderson, Coles and H\"usler \cite{anderson1997maxima}, where they study the maximum of a sum of a large number of triangular arrays.

In order to get fluid limits for the maximum queue lengths, we need to study diffusion limits for the individual queue lengths. We thus combine ideas from the literature on extreme value theory with literature on diffusion approximations, which we show in Section \ref{subsec: scaling explanation}. In order to be able to analyze the queue lengths through diffusion approximations, we impose a heavy-traffic assumption, namely
$(1-\rho_N)N^2\to\beta$. Then, for each separate queue length, we have a reflected Brownian motion as diffusion approximation. By using the well-known formula for the cumulative distribution function of a reflected Brownian motion (cf.\ Harrison \cite[p.~49]{harrison1985brownian}), we investigate the maximum of $N$ independent reflected Brownian motions to get an idea of the scaling of the maximum queue length.

Now, we give a brief sketch of how we apply these ideas to prove the fluid limit, we start by considering the slightly simpler scenario that each queue is empty at time 0. Because we want to prove a fluid limit that holds uniformly on compact intervals, we need to prove pointwise convergence of the process and tightness of the collection of processes. Our first step in proving this is by showing that each queue length is in distribution the same as a supremum of an arrival process minus a service process. We then show in Section \ref{subsec: scaling explanation} that under the temporal scaling of $tN^3\log N$ and the spatial scaling of $N\log N$, the arrival process minus a drift term converges to $-\beta t$, as $N\to\infty$. Furthermore, we derive under that same temporal scaling but under a spatial scaling of $N\sqrt{\log N}$, that the centralized service process satisfies the central limit theorem. This scaled centralized service process is given in Equation \eqref{eq: normal like random variable}. We use the non-uniform Berry-Ess\' een inequality, which is described by Michel in \cite{michel1976constant}, to deduce the convergence rate of the cumulative distribution function of this scaled centralized service process to the cumulative distribution function of a normally distributed random variable, which is given in Equation \eqref{eq: inequality michel}. It turns out that this convergence rate is fast enough so that we can replace the scaled centralized service process with a normally distributed random variable in the expression of the maximum queue length in order to get the same limit. By Pickands' result \cite{pickands1968moment} on the convergence of moments of the maximum of $N$ scaled random variables, we know that the expectation of the maximum of standard normally distributed random variables divided by $\sqrt{\log N}$ converges to $\sqrt{2}$, as $N\to\infty$. This gives us the convergence of the maximum of $N$ scaled centralized service processes. After we have obtained these limiting results for the scaled arrival and service process, we use these, together with Doob's maximal submartingale inequality to prove convergence in probability of the maximum queue length, we show this in Section \ref{subsubsec: pointwise convergence}. Finally, in Section \ref{subsubsec: tightness} we use Doob's maximal submartingale inequality to bound the probability that the process makes large jumps and prove that this probability is small so that the maximum queue length is a tight process.

After we have considered the maximum queue length for the process with empty queues at time 0, we then turn to the scenario that the length of each queue at time 0 is identically distributed. In this case, we can use Lindley's recursion to express the maximum queue length as the pairwise maximum of the maximum queue length with empty queues at time 0 and a part depending on the number of jobs at time 0, this formula is given in Equation \eqref{eq: max queue length representation}. How to prove the fluid limit for the first part is already sketched above. In order to derive a fluid limit for the latter part, we first observe that this part equals the maximum of $N$ times the sum of the number of jobs at time 0 at each server plus the number of arrivals minus the number of services at each server. Following a similar path as earlier, we can prove that the scaled centralized service process at server $i$ behaves like a normally distributed random variable. Thus, we have to analyze a maximum of $N$ pairwise sums of normally distributed random variables and random variables describing the number of jobs at time 0, which is stated in more detail in Lemma \ref{lem: pointwise convergence startposition}.

In Lemma \ref{lem: approximation starting point convergence} we prove a convergence result of this maximum, this is quite a challenge because we need to apply extreme value theory on pairwise sums. In order to do this, we use the results from Davis, Mulrow \& Resnick \cite{davis1988almost} and Fisher \cite{fisher1969limiting} on the convergence of samples of random variables to limiting sets. The authors prove convergence results of the convex hull of $\{(\frac{Z^{(1)}_i}{b_N},\ldots,\frac{Z^{(k)}_i}{b_N})_{i\leq N}\}$ to a limiting set, as $N\to\infty$, with $(Z^{(j)}_i,i\leq N)$ i.i.d., $Z^{(j)}_i$ and $Z^{(l)}_m$ are independent and $b_N$ is a proper scaling sequence. We show in the proof of Lemma \ref{lem: extreme value convergence} that these results can be extended in establishing convergence of extreme values of $\max_{i\leq N}\sum_{j=1}^k \frac{Y^{(j)}_i}{a_N^{(j)}}$, where $a_N^{(l)}$ and $a_N^{(m)}$ are not necessarily the same, which is a stand-alone result of independent interest. We did not find this extension in other literature. The result in Lemma \ref{lem: approximation starting point convergence} follows from Lemma \ref{lem: extreme value convergence}.

The rest of the paper is organized as follows. In Section \ref{sec: Model}, we describe the fork-join system in more detail; we give a definition of the arrival and service processes and we present a scaled version of the queueing model. In Section \ref{subsec: fluid limit model description}, we introduce the fluid limit and explain it heuristically. We elaborate a bit more on the scaling and the shape of the fluid limit in Sections \ref{subsec: scaling explanation} and \ref{subsec: shape fluid limit}. Furthermore, we give some examples and numerical results in Section \ref{subsec: examples and numerics}. We finish with some concluding remarks in Section \ref{sec: future work}. The proof of the fluid limit is given in Section \ref{sec: Proofs}. In Appendix \ref{app: taylor approx theta}, we elaborate on the convergence of the upper bound that was given in Lemma \ref{lem: upper bound pointwise convergence t large}. We prove in Appendix \ref{app: extreme values} a convergence result of $\max_{i\leq N}\sum_{j=1}^k \frac{Y^{(j)}_i}{a_N^{(j)}}$. In Appendix \ref{app: proof lemma convergence fourth moment}, we prove the lemmas stated in Section \ref{subsubsec: useful lemmas}. An overview of all notation is given in Appendix \ref{app: Notation}.


\section[Model description and main results]{Model description and main results.}\label{sec: Model}

We now turn to a formal definition of the fork-join queue that we study. We consider a fork-join queue with integer-valued arrivals and services. In this queueing system, there is one arrival process. The arriving tasks are divided into $N$ subtasks which are completed by $N$ servers. We assume that both the number of arrivals and services per time step are Bernoulli distributed. The parameters of the Bernoulli random variables depend on the number of servers. This is formalized in Definitions \ref{def: arrival process} and \ref{def: service process i server}.
\begin{definition}[Arrival process]\label{def: arrival process}
The random variable $\Arrival{n}$ indicates the number of arrivals up to time $n$ and equals
\begin{align*}
\Arrival{n}=\sum_{j=1}^{\lfloor n\rfloor} \ArrivalBernoulli{j}
\end{align*}
with $\ArrivalBernoulli{j}$ indicating whether or not there is an arrival at time $j$. $\ArrivalBernoulli{j}$ is a Bernoulli random variable with parameter $\ArrivalProbability$. So,

\begin{align*}
\ArrivalBernoulli{j}=\left\{\begin{matrix}
1&\text{ w.p. }&  &\ArrivalProbability,\\ 
0&\text{ w.p. }& 1-&\ArrivalProbability.
\end{matrix}\right.
\end{align*}

\end{definition}

\begin{definition}[Service process $i$-th server]\label{def: service process i server}
The random variable $\Service{i}{n}$ describes the number of potentially completed tasks of the $i$-th server in the fork-join queue at time $n$ with

\begin{align*}
\Service{i}{n}=\sum_{j=1}^{\lfloor n\rfloor} \ServiceBernoulli{i}{j},
\end{align*}
where $\ServiceBernoulli{i}{j}$ is a Bernoulli random variable with parameter $\SuccessProbability$ indicating whether the $i$-th server completed a service at time $j$. 
\begin{align*}
\ServiceBernoulli{i}{j}=\left\{\begin{matrix}
1&\text{ w.p. }&  &\SuccessProbability,\\ 
0&\text{ w.p. }& 1-&\SuccessProbability.
\end{matrix}\right.
\end{align*}
Both $\ArrivalProbability$ and $\SuccessProbability$ are taken as functions of $N$, which we specify in Definition \ref{def: maximum queue length} below.
\end{definition}

We assume that for all $N\geq 1$ the random variables $(\ArrivalBernoulli{j},j\geq 1)$ are mutually independent for all $j$ and $(\ServiceBernoulli{i}{j},j\geq 1,i\leq N)$ are mutually independent for all $j$ and $i$. We also assume that an incoming task can be completed in the same time slot as in which the task arrived. Finally, we assume that $\ArrivalBernoulli{j}$ and $\ServiceBernoulli{i}{j}$ are independent, in other words, $\ServiceBernoulli{i}{j}$ could still be 1 while there are no tasks to be served at server $i$ at time $j$. Due to this assumption, we have on the one hand the beneficial situation that $(\Arrival{n},n\geq 0)$ and $(\Service{i}{n},n\geq 0)$ are independent processes, but on the other hand, we should be careful with defining the queue length. However, it is a well-known result that we can use Lindley's recursion, and write the queue length of the $i$-th server at time $n$ as
\begin{align*}
\sup_{0\leq k\leq n}\left[\Big(\Arrival{n}-\Arrival{k}\Big)-\left(\Service{i}{n}-\Service{i}{k}\right)\right],
\end{align*}
provided that the queue length is 0 at time 0. This is in distribution equal to
\begin{align*}
\sup_{0\leq k\leq n}\left(\Arrival{k}-\Service{i}{k}\right).
\end{align*}
As can be seen in this expression, the queue lengths of different servers are mutually dependent, since the arrival process is the same. When at time 0 there are already jobs in the queue, then we can, after again applying Lindley's recursion, write the queue length of the $i$-th server at time $n$ as 
\begin{align*}
\max\left(\sup_{0\leq k\leq n}\left[\Big(\Arrival{n}-\Arrival{k}\Big)-\left(\Service{i}{n}-\Service{i}{k}\right)\right],\QueueZero{i}+\Arrival{n}-\Service{i}{n}\right),
\end{align*}
with $\QueueZero{i}$ the number of jobs in front of the $i$-th server at time 0. Observe that the queue length of the $i$-th server equals the maximum of the queue length when the number of jobs at time 0 would be 0, and a random variable that depends on the initial number of jobs.

The aim of this work is to investigate the behavior of the fork-join queue when the number of servers $N$ is very large. The main objective is deriving the distribution of the largest queue, as this represents the slowest supplier, which is the bottleneck for the manufacturer. Therefore, we define in Definition \ref{def: maximum queue length} a random variable indicating the maximum queue length at time $n$. Furthermore, we explore this model in the heavy-traffic regime. To this end, we let $\ArrivalProbability$ and $\SuccessProbability$ go to 1 at similar rates, so that the arrivals and services are nearly deterministic processes. 
\begin{definition}[Maximum queue length at time $n$]\label{def: maximum queue length}
Let $\ArrivalProbability=1-\frac{\alpha}{N}-\frac{\beta}{N^2}$ and $\SuccessProbability=1-\frac{\alpha}{N}$, with $\alpha,\beta>0$. Let $\MaxQueueLength{n}$ be the maximum queue length of $N$ parallel servers at time $n$, with $\MaxQueueLength{0}=0$. Then
\begin{align}\label{eq: maxqueuelength definition}
\MaxQueueLength{n}=\max_{i\leq N}\sup_{0\leq k\leq n}\left[\Big(\Arrival{n}-\Arrival{k}\Big)-\left(\Service{i}{n}-\Service{i}{k}\right)\right].
\end{align}
\end{definition}
So,
\begin{align}\label{eq: maxqueuelength simplification distribution}
\MaxQueueLength{n}\overset{d}=\max_{i\leq N}\sup_{0\leq k\leq n}\left(\Arrival{k}-\Service{i}{k}\right),
\end{align}
under the assumption that $\MaxQueueLength{0}=0$. From these choices of $\ArrivalProbability$ and $\SuccessProbability$, it follows that the traffic intensity $\rho_N$ of a single queue satisfies $(1-\rho_N)N^2\to \beta$, as $N\to\infty$. Furthermore, if $\QueueZero{i}>0$, the maximum queue length at time $n$ can be written as
\begin{align}\label{eq: max queue length representation}
\MaxQueueLength{n}=\max_{i\leq N}\max\bigg(&\sup_{0\leq k\leq n}\left[\Big(\Arrival{n}-\Arrival{k}\Big)-\left(\Service{i}{n}-\Service{i}{k}\right)\right]\nonumber\\
,&\QueueZero{i}+\Arrival{n}-\Service{i}{n}\bigg).
\end{align}
Observe that we can interchange the order of the $\max_{i\leq N}$ term and the $\max$ term, and rewrite the expression in \eqref{eq: max queue length representation} as the pairwise maximum of two random variables, one random variable is the maximum of $N$ queue lengths with initial condition 0, as given in Equation \eqref{eq: maxqueuelength definition}, and the other is the maximum of $N$ sums of the queue length at time 0 plus the number of arrivals minus the number of services.
\subsection[Fluid limit]{Fluid limit.}\label{subsec: fluid limit model description}
As we just have formally defined the fork-join queue that we study, with the particular nearly deterministic setting, we now state and explain the main result of this paper. 
Our central result is a fluid approximation for the rescaled maximum queue length process, which is given in Theorem \ref{thm: fluid limit}. We prove that under a certain spatial and temporal scaling the maximum queue length converges to a continuous function, which depends on time $t$.

There is, however, not a straightforward procedure in choosing the temporal and spatial scaling, there are namely more possibilities that lead to a non-trivial limit. For instance, when we choose a temporal scaling of $N^3$ and a spatial scaling of $N\sqrt{\log N}$, we get the fluid limit that is given in Proposition \ref{cor: other fluid limit}. Here, we assume that the initial condition is 0.
\begin{proposition}[Temporal scaling of $N^3$ and spatial scaling of $N\sqrt{\log N}$]\label{cor: other fluid limit}
For\\ $\MaxQueueLength{0}=0$, $\alpha>0$ and $\beta>0$, with
\begin{enumerate}
\item $\ArrivalProbability=1-\frac{\alpha}{N}-\frac{\beta}{N^{2}}$,
\item $\SuccessProbability=1-\frac{\alpha}{N}$, we have
\end{enumerate}
\begin{align*}
\probability*{\sup_{0\leq s\leq T}\left|\frac{\MaxQueueLength{sN^3}}{N\sqrt{\log N}}-\sqrt{2\alpha s}\right|>\epsilon}
\LimitN 0 \text{ }\forall \epsilon>0.
\end{align*}
\end{proposition}
However, we can also derive a steady-state limit, which is given in Proposition \ref{cor: steady state convergence}.
\begin{proposition}[Steady-state convergence]\label{cor: steady state convergence}
For $\alpha>0$ and $\beta>0$, with
\begin{enumerate}
\item $\ArrivalProbability=1-\frac{\alpha}{N}-\frac{\beta}{N^{2}}$,
\item $\SuccessProbability=1-\frac{\alpha}{N}$, we have
\end{enumerate}
\begin{align*}
\FluidProcessTwo{\infty}\LimitP \frac{\alpha}{2\beta}\text{ as $N\to\infty$}.
\end{align*}
\end{proposition}
As we can see in Proposition \ref{cor: steady state convergence}, to obtain a non-trivial steady-state limit, we need a spatial scaling of $N\log N$. Since this is the only choice that leads to a non-trivial limit, it is a natural choice to look for a fluid limit which also has this spatial scaling. Our main result, stated in Theorem \ref{thm: fluid limit}, is such a fluid limit, and it turns out that for establishing this limit, we need a temporal scaling of $N^3\log N$. In Section \ref{subsec: scaling explanation} we explain why these scalings are natural. We omit the proof of Proposition \ref{cor: other fluid limit}, but we do explain how Proposition \ref{cor: other fluid limit} is connected to Theorem \ref{thm: fluid limit} at the end of this section. Furthermore, we give a proof of Proposition \ref{cor: steady state convergence} in Section \ref{sec: Proofs}.

We now mention and discuss some assumptions under which our main result holds. First of all, we assume that we have nearly deterministic arrivals and services.
\begin{assumption}\label{ass: 1}
$\ArrivalProbability=1-\frac{\alpha}{N}-\frac{\beta}{N^{2}}$ and $\SuccessProbability=1-\frac{\alpha}{N}$, with $\alpha,\beta>0$.
\end{assumption}
Secondly, we have a basic assumption on the initial condition.
\begin{assumption}\label{ass: 2}
$(\QueueZero{i},i\leq N)$ are i.i.d.\ and non-negative for all $N$.
\end{assumption}
Furthermore, we want to prove a fluid limit with a spatial scaling of $N\log N$. Therefore, we need to assume that the maximum number of jobs at time 0 also scales with $N\log N$. In order to do so, we allow $(\QueueZero{i},i\leq N,N\geq 1)$ to be a triangular array, i.e.\ a doubly indexed sequence with $i\leq N$. This is a necessity because otherwise we would be limited to distributions where the maximum scales like $N\log N$, which would lead us to the family of the heavy-tailed distributions for which we do not have convergence in probability of its maximum. Thus in our setting, $\QueueZero{i}$ and $Q_i^{(N+1)}(0)$ do not need to be the same. Consequently, we need to have some regularity on $\QueueZero{i}$ as $N$ increases to be able to prove a limit theorem.
\begin{assumption}\label{ass: 3}
$\frac{\MaxQueueLength{0}}{(N\log N)}\LimitP q(0),$ with $q(0)\geq 0$, as $N\to\infty$, with $\QueueZero{i}=\lfloor r_N \StartRVIndep{i}\rfloor$, where $r_N$ is a scaling sequence.
\end{assumption}
Finally, we can distinguish two cases in which Theorem \ref{thm: fluid limit} holds. 
\begin{assumption}\label{ass: 4}
$\StartRVIndep{i}$ has a finite right endpoint.
\end{assumption}
\begin{assumption}\label{ass: 5}
$\StartRVIndep{i}$ is a continuous random variable and for all $v\in[0,1]$, 
\begin{align*}
\lim_{t\to\infty}\frac{-\log \left(\probability*{\StartRVIndep{i}>vt}\right)}{-\log \left(\probability*{\StartRVIndep{i}>t}\right)}=h(v).
\end{align*}
\end{assumption}
Before stating the theorem, we would like to give two remarks on Assumption \ref{ass: 5}. First of all, the function $h$ has the property that for all $u,v\in[0,1]$, $h(uv)=h(u)h(v)$. Thus, if $h$ is continuous, $h(v)=v^a$, with $a>0$. When $h$ is discontinuous, there are two possibilities: $h(v)=\mathbbm{1}(v>0)$, or $h(v)=\mathbbm{1}(v=1)$, this corresponds to $h(v)=v^a$ with $a=0$ and $a=\infty$, respectively. Secondly, the assumption of continuity of $\StartRVIndep{i}$ can be removed, which would lead to more cumbersome proofs.

\begin{theorem}[Fluid limit with a non-zero initial condition]\label{thm: fluid limit}
If Assumptions \ref{ass: 1}, \ref{ass: 2} and \ref{ass: 3} hold, and either Assumption \ref{ass: 4} or Assumption \ref{ass: 5} holds, then we have $\forall T>0$, that
\begin{align}
\mathbb{P}\bigg(\sup_{0\leq t\leq T}&\left|\FluidProcessTwo{t}-q(t)\right|
>\epsilon\bigg)\LimitN 0\text{ }\forall \epsilon>0,
\end{align}
with 
\begin{align}\label{eq: fluid limit function}
q(t)=\max&\bigg(\bigg(\sqrt{2\alpha t}-\beta t\bigg)\mathbbm{1}\bigg(t<\frac{\alpha}{2\beta^2}\bigg)+\frac{\alpha}{2\beta}\mathbbm{1}\bigg(t\geq\frac{\alpha}{2\beta^2}\bigg),g(t,q(0))-\beta t\bigg).
\end{align}
The function $g(t,q(0))$ has the following properties:
\begin{enumerate}
\item If Assumption \ref{ass: 4} holds, then 
\begin{align}\label{eq: limit scenario 1}
g(t,q(0))=q(0)+\sqrt{2\alpha t}.
\end{align}
\item If Assumption \ref{ass: 5} holds, then 
\begin{align}\label{eq: limit scenario 2}
g(t,q(0))=\sup_{(u,v)}\{\sqrt{2\alpha t}u+q(0)v|u^2+h(v)\leq 1,0\leq u\leq 1, 0\leq v\leq 1\}.
\end{align}
\end{enumerate}
\end{theorem}
There is a connection between Assumptions \ref{ass: 4} and \ref{ass: 5} on $\StartRVIndep{i}$ and extreme value theory. If Assumption \ref{ass: 4} holds, then this means that $\StartRVIndep{i}$ is either a degenerate random variable or is in the domain of attraction of the Weibull distribution. On the other hand, if Assumption \ref{ass: 5} holds, then $\StartRVIndep{i}$ is in the domain of attraction of the Gumbel distribution.

In order to allow dependence between the initial number of jobs at different servers, we can also replace Assumptions \ref{ass: 2} and \ref{ass: 3} with the following assumption.
\begin{assumption}\label{ass: 6}
Let $\QueueZero{i}=\StartRVTriangIndep{i}+\StartRVTriangDep{i},$
with $\StartRVTriangIndep{i}=\lfloor r_N\StartRVIndep{i}\rfloor$, where $(\StartRVIndep{i},i\leq N)$ are i.i.d.\ and non-negative, and satisfy either Assumption \ref{ass: 4} or \ref{ass: 5}. Furthermore, $\StartRVTriangDep{i}$ is non-negative, and $\max_{i\leq N}\frac{\StartRVTriangDep{i}}{(N\log N)}\LimitP 0 $, as $N\to\infty$.
\end{assumption}
When Assumption \ref{ass: 6} is satisfied, there may be mutual dependence between $\QueueZero{i}$ and $\QueueZero{j}$, because  $\StartRVTriangDep{i}$ and $\StartRVTriangDep{j}$ may be mutually dependent.

As can be seen in Theorem \ref{thm: fluid limit}, the fluid limit has an unusual form, $q(t)$ is namely a maximum of two functions. The first part of this maximum is the fluid limit when the initial number of jobs equals 0 and the second part is caused by the initial number of jobs. We elaborate on this more in Section \ref{subsec: shape fluid limit}. The $\log N$ term in the spatial and temporal scaling of the process is also unusual. We show in Section \ref{subsec: scaling explanation} that this is due to the fact that we take a maximum of $N$ random variables, with $N$ large. Scaling terms like $(\log N)^c$ are in this context very natural.

We mentioned earlier that different choices for temporal and spatial scalings lead to a fluid limit. We gave Proposition \ref{cor: other fluid limit} as an example. Since we analyze one and only one system, the two fluid limits that we presented should be connected to each other. An easy way to see this is by observing that from Theorem \ref{thm: fluid limit} it follows that when $\MaxQueueLength{0}=0$,
\begin{align*}
\frac{\MaxQueueLength{tN^3\log N}}{N\log N}\LimitP \sqrt{2\alpha t}-\beta t\text{ as $N\to\infty$},
\end{align*}
for $t<\frac{\alpha}{(2\beta^2)}$. Thus, for all $t>0$ and for $N$ large, we expect that
$\frac{\MaxQueueLength{tN^3}}{(N\sqrt{\log N})}\approx \sqrt{2\alpha t}-\beta\frac{t}{\sqrt{\log N}}\LimitN \sqrt{2\alpha t}.$ This shows heuristically how Proposition \ref{cor: other fluid limit} is connected with Theorem \ref{thm: fluid limit}. The formal proof of Proposition \ref{cor: other fluid limit} is analogous to the proof of Theorem \ref{thm: fluid limit} and is omitted in this paper.
\subsection[Scaling]{Scaling.}\label{subsec: scaling explanation}
In Section \ref{subsec: fluid limit model description}, we presented the fluid limit under the rather unusual temporal scaling of $N^3\log N$ and spatial scaling of $N\log N$. A heuristic justification for these scalings can be given by using extreme value theory and ideas from literature on diffusion approximations. In particular, for the spatial scaling, we argue as follows: as we are interested in the convergence of the maximum queue length, we can use a central limit result to replace each separate queue length with a reflected Brownian motion and use extreme value theory to get a heuristic idea of the convergence of the scaled maximum queue length. To argue this, first observe that the arrival and service processes are binomially distributed random variables, and we can compute the expectation and variance of $\ScaledProcessLine$ as
\begin{align}\label{eq: expectation queue length}
\expect*{\ScaledProcess}=-\beta t\sqrt{\log N}+o_N(1),
\end{align}
and
\begin{align}\label{eq: variance queue length}
&\text{Var}\left(\ScaledProcess\right)\nonumber\\
=&\frac{1}{N^{2}\log N}\lfloor{tN^{3}\log N\rfloor}\left(\left(\frac{\alpha}{N}+\frac{\beta}{N^{2}}\right)\left(1-\frac{\alpha}{N}-\frac{\beta}{N^{2}}\right)+\frac{\alpha}{N}\bigg(1-\frac{\alpha}{N}\bigg)\right)\nonumber\\
=&2\alpha t+o_N(1).
\end{align}
From this, a non-trivial scaling limit can be easily deduced: observe that $\Arrival{{tN^{3}\log N}}-\Service{i}{{tN^{3}\log N}}$ is a sum of independent and identically distributed random variables, so this implies that 
\begin{align*}
\ScaledProcess\overset{d}{\approx} Z_i,
\end{align*}
as $N$ is large, with $Z_i\sim \mathcal{N}\left(-\beta t\sqrt{\log N},2\alpha t\right)$. Furthermore, because $\Arrival{{tN^{3}\log N}}-\Service{i}{{tN^{3}\log N}}$ is, in fact, the difference of two random walks, we also have 
\begin{align*}
\sup_{0\leq n\leq tN^{3}\log N}\frac{1}{N\sqrt{\log N}}\left(\Arrival{{n}}-\Service{i}{{n}}\right)\overset{d}{\approx} R_i(t),
\end{align*}
as $N$ is large, with $R_i(t)$ a reflected Brownian motion for $t$ fixed. We can apply extreme value theory to show that $\max_{i\leq N}R_i(t)$ scales with $\sqrt{\log N}$. This can be deduced from the cumulative distribution function of the reflected Brownian motion which is given in \cite[p.~49]{harrison1985brownian}. Concluding, the proper spatial scaling of the fluid limit in Theorem \ref{thm: fluid limit} is $\frac{1}{(N\log N)}$.

As Equations \eqref{eq: expectation queue length} and \eqref{eq: variance queue length} show, the right temporal and spatial scalings are determined by the choice of the arrival and service probability. When we change the arrival probability to $\ArrivalProbability=1-\frac{\alpha}{N}-\frac{\beta}{N^{1+c}}$, with $c\geq 1$, and keep the service probability the same, we can derive in the same manner, that under a different temporal and spatial scaling of the queueing process, the fluid limit result still holds; we state this in Proposition \ref{cor: other arrivals services}. 
\begin{proposition}[Other arrival and service probabilities]\label{cor: other arrivals services}
For $c\geq 1$, $\alpha>0$ and $\beta>0$, with
\begin{enumerate}
\item $\ArrivalProbability=1-\frac{\alpha}{N}-\frac{\beta}{N^{1+c}}$,
\item $\SuccessProbability=1-\frac{\alpha}{N}$,
\end{enumerate}
and $\MaxQueueLength{0}=O(N^c\log N)$ and satisfies the same assumptions as in Theorem \ref{thm: fluid limit}, then
\begin{align*}
\mathbb{P}\bigg(\sup_{0\leq t\leq T}\bigg|&\frac{\MaxQueueLength{tN^{1+2c}\log N}}{N^c\log N}
-q(t)\bigg|>\epsilon\bigg)
\LimitN 0 \text{ }\forall \epsilon>0.
\end{align*}
\end{proposition}
The proof of this proposition is very similar to the proof of Theorem \ref{thm: fluid limit}. Thus we omit it here.
\subsection[Shape of the fluid limit]{Shape of the fluid limit.}\label{subsec: shape fluid limit}
In Section \ref{subsec: scaling explanation}, we gave a heuristic explanation of the temporal and spatial scaling of the process. Here we do the same for the shape of the fluid limit. First of all, we rewrite the expression in \eqref{eq: max queue length representation} and get that the scaled maximum queue length satisfies
\begin{align}\label{eq: maxqueue length expression}
&\FluidProcessTwo{t}=\nonumber\\
\max\vast(&\max_{i\leq N}\sup_{0\leq s\leq t}\frac{\Big(\Arrival{tN^3\log N}-\Arrival{sN^3\log N}\Big)-\left(\Service{i}{tN^3\log N}-\Service{i}{sN^3\log N}\right)}{N\log N},\nonumber\\
&\max_{i\leq N}\frac{\Arrival{tN^3\log N}+\Service{i}{tN^3\log N}+\QueueZero{i}}{N\log N}\vast).
\end{align}
Now, observe that when $\QueueZero{i}=0$ for all $i$, the pairwise maximum in \eqref{eq: maxqueue length expression} simplifies to the first part of the maximum. Furthermore, it turns out that the first and the second part of this maximum converge to the first and second part of the maximum in \eqref{eq: fluid limit function}, respectively. To see the first limit heuristically, observe that, due to the central limit theorem, 
\begin{align*}
\ScaledProcess \overset{d}{\approx}\vartheta_i+\zeta,
\end{align*}
with $\vartheta_i\sim \mathcal{N}(0,\alpha t)$, independently for all $i$, and $\zeta\sim\mathcal{N}(-\beta t\sqrt{\log N},\alpha t)$.
We can write $\max_{i\leq N} (\vartheta_i+\zeta)=\max_{i\leq N}(\vartheta_i)+\zeta$. Then, by the basic convergence result that the maximum of $N$ i.i.d.\ standard normal random variables scales like $\sqrt{2\log N}$, it is easy to see that
$\frac{\max_{i\leq N} (\vartheta_i+\zeta)}{\sqrt{\log N}}\LimitP \sqrt{2\alpha t}-\beta t\text{ as $N\to\infty$}$. Because of the fact that a queue length which is 0 at time 0, can be written as the supremum of the arrival process minus the service process up to time $t$, the fluid limit yields
$\sup_{0\leq s\leq t}(\sqrt{2\alpha s}-\beta s)$, which equals the first part of the maximum in \eqref{eq: fluid limit function}.

Similarly, for the second part in \eqref{eq: maxqueue length expression} we observe that
\begin{align}\label{eq: split start position}
&\max_{i\leq N}\StartPositionGeneral{i}{t}\nonumber\\
=&\frac{\Arrival{tN^3\log N}-\left(1-\frac{\alpha}{N}\right)tN^3\log N}{N\log N}
+\max_{i\leq N}\frac{\left(1-\frac{\alpha}{N}\right)tN^3\log N-\Service{i}{tN^3\log N}+\QueueZero{i}}{N\log N}.
\end{align}
It is easy to see that the first term converges to $-\beta t$ as $N\to\infty$, and we prove later on that the second term converges to $g(t,q(0))$. This explains the second part of the fluid limit in \eqref{eq: fluid limit function}.

Specific properties of the function $g$ can be deduced. First of all, Assumption \ref{ass: 4} considers the case that $\StartRVIndep{i}$ has a finite right endpoint. In this scenario, we have that $
\frac{\QueueZero{i}}{(N\log N)}=\frac{\lfloor{r_N\StartRVIndep{i}\rfloor}}{(N\log N)}=\frac{\lfloor{N\log N\StartRVIndep{i}\rfloor}}{(N\log N)}\approx \StartRVIndep{i}$. Now, the theorem says that $g(t,q(0))=q(0)+\sqrt{2\alpha t}$. This actually means that for large $N$,
\begin{align*}
&\max_{i\leq N}\left(U_i+\frac{\left(1-\frac{\alpha}{N}\right)tN^3\log N-\Service{i}{tN^3\log N}}{N\log N}\right)\nonumber\\
\approx& \max_{i\leq N}U_i+ \max_{i\leq N}\frac{\left(1-\frac{\alpha}{N}\right)tN^3\log N-\Service{i}{tN^3\log N}}{N\log N}.
\end{align*}
This behavior can be very well explained, because due to the assumption that $\StartRVIndep{i}$ has a finite right endpoint, there will be many observations of $\StartRVIndep{i}$ that are close to the right endpoint, as $N$ becomes large, and thus it will be more and more likely that there is a large observation $\Big(\left(1-\frac{\alpha}{N}\right)(tN^3\log N)-\Service{i^{\star}}{tN^3\log N}\Big)\Big/\Big(N\log N\Big)$, for which the observation $\StartRVIndep{i^{\star}}$ will also be large.

Furthermore, when Assumption \ref{ass: 5} holds, $g(t,q(0))$ can be written as a supremum over a set. To give an idea of why this is the case, we first observe that we can write the last term in \eqref{eq: split start position} as
\begin{align}\label{eq: maximum start position}
\max_{i\leq N}\left(\frac{\left(1-\frac{\alpha}{N}\right)(tN^3\log N)-\Service{i}{tN^3\log N}}{N\log N}+\frac{\QueueZero{i}}{N\log N}\right).
\end{align}
Thus, this maximum can be viewed as a maximum of $N$ pairwise sums of random variables. For any $N>0$, we can write down all the $N$ pairs of random variables as
\begin{align}\label{eq: set random variable}
\left\{\left(\frac{1}{\sqrt{2\alpha t}}\frac{\left(1-\frac{\alpha}{N}\right)(tN^3\log N)-\Service{i}{tN^3\log N}}{N\log N},\frac{1}{q(0)}\frac{\QueueZero{i}}{N\log N}\right)_{i\leq N}\right\}.
\end{align} 
Now, the expression in Equation \eqref{eq: maximum start position} can be written as $\sqrt{2\alpha t}u+q(0)v$ with $(u,v)$ in the set in \eqref{eq: set random variable}, such that $\sqrt{2\alpha t}u+q(0)v$ is maximized. Due to the central limit theorem, the first term in \eqref{eq: set random variable} can be approximated by $\frac{\vartheta_i}{\sqrt{2\alpha t}}$ with $\vartheta_i\sim \mathcal{N}(0,\alpha t)$ when $N$ is large. Therefore, the convex hull of the set in \eqref{eq: set random variable} looks like the convex hull of the set
\begin{align*}
\left\{\left(\frac{1}{\sqrt{2\alpha t}}\frac{\vartheta_i}{\sqrt{\log N}},\frac{1}{q(0)}\frac{\QueueZero{i}}{N\log N}\right)_{i\leq N}\right\}.
\end{align*}
The convex hull of this set can be seen as a random variable, and converges, under an appropriate metric, in probability to the limiting set
\begin{align}\label{eq: limit set}
\{(u,v)|u^2+h(v)\leq 1,-1\leq u\leq 1, 0\leq v\leq 1\},
\end{align}
in $\mathbb{R}^2$, as $N\to\infty$; cf.\ \cite{davis1988almost} and \cite{fisher1969limiting} for details on this. Our intuition says that the limit of the expression in \eqref{eq: maximum start position} is attained at the coordinate $(u,v)$ in the closure of the limiting set given in \eqref{eq: limit set}, such that $\sqrt{2\alpha t}u+q(0)v$ is maximized. We show that this is indeed correct. In fact, we prove this in Lemma \ref{lem: approximation starting point convergence} in a more general context than in \cite{davis1988almost} and \cite{fisher1969limiting}. In \cite{davis1988almost} and \cite{fisher1969limiting}, the authors make the assumption that the scaling sequences are the same, so the analysis is restricted to samples of the type $\{(X_i/a_N,Y_i/a_N)_{i\leq N}\}$. However, we show that for proving convergence of the maximum of the pairwise sum, the scaling sequences do not need to be the same. 
\subsection[Examples and numerics]{Examples and numerics.}\label{subsec: examples and numerics}
In Section \ref{subsec: shape fluid limit}, we showed that the shape of the fluid limit depends on the distribution of the number of jobs at time 0. Here, we give some basic examples of how the fluid limit is influenced by the distribution of the number of jobs at time 0. We also present and discuss some numerical results.

As a first example, for $\StartRVIndep{i}=X_i^+$, with $X_i\sim\mathcal{N}(0,1)$, we can write for $v>0$, $\probability*{\StartRVIndep{i}>v}=\exp(-v^2L(v))$, such that $L$ is slowly varying. Thus for $v\in[0,1]$,
\begin{align*}
h(v)=\lim_{t\to\infty}\frac{-\log \left(\probability*{\StartRVIndep{i}>vt}\right)}{-\log \left(\probability*{\StartRVIndep{i}>t}\right)}=\lim_{t\to\infty}\frac{(vt)^2L(vt)}{t^2L(t)}=v^2.
\end{align*}
Thus,
\begin{align*}
g(t,q(0))=\sup_{(u,v)}\{\sqrt{2\alpha t}u+q(0)v|u^2+v^2\leq 1,-1\leq u\leq 1, 0\leq v\leq 1\}=\sqrt{q(0)^2+2\alpha t}.
\end{align*}  
Concluding,
\begin{align*}
&\max_{i\leq N}\frac{\left(1-\frac{\alpha}{N}\right)(tN^3\log N)-\Service{i}{tN^3\log N}+\QueueZero{i}}{N\log N}\nonumber\\
=&\max_{i\leq N}\frac{\left(1-\frac{\alpha}{N}\right)(tN^3\log N)-\Service{i}{tN^3\log N}+\lfloor{\frac{q(0)N\log N\StartRVIndep{i}}{\sqrt{2\log N}}\rfloor}}{N\log N}\nonumber\\
\LimitP& \sqrt{q(0)^2+2\alpha t}-\beta t\text{ as $N\to\infty$},
\end{align*}
where $r_N=\frac{q(0)N\log N}{\sqrt{2\log N}}$, such that $\frac{\MaxQueueLength{0}}{(N\log N)}\LimitP q(0)$, as $N\to\infty$.\newline
Another example is, when we assume that $\StartRVIndep{i}$ is lognormally distributed, we know that $\probability*{\StartRVIndep{i}>v}=\probability*{X_i>\log v}$, with $X_i\sim\mathcal{N}(0,1)$. Thus, $\probability*{\StartRVIndep{i}>v}=\exp(-\mathbbm{1}(v>0)\log(v)^2L(\log v))$. Then, for $v\in[0,1]$,
\begin{align*}
h(v)=\lim_{t\to\infty}\frac{\mathbbm{1}(v>0)\log(vt)^2L(\log(vt))}{\log(t)^2L(\log(t))}=\mathbbm{1}(v>0).
\end{align*}
In this case, we have that
\begin{align*}
g(t,q(0))=\sup_{(u,v)}\{\sqrt{2\alpha t}u+q(0)v|u^2+\mathbbm{1}(v>0)\leq 1,-1\leq u\leq 1, 0\leq v\leq 1\}=\max(q(0),\sqrt{2\alpha t}).
\end{align*}  
\indent We also consider the case $\probability*{\StartRVIndep{i}>v}= \exp(1-\exp(v))$, then for $v\in[0,1]$,
\begin{align*}
\lim_{t\to\infty}\frac{-\log \left(\probability*{\StartRVIndep{i}>vt}\right)}{-\log \left(\probability*{\StartRVIndep{i}>t}\right)}=\lim_{t\to\infty}\frac{\exp(vt)-1}{\exp(t)-1}=\mathbbm{1}(v=1).
\end{align*}
Then,
\begin{align*}
g(t,q(0))=\sup_{(u,v)}\{\sqrt{2\alpha t}u+q(0)v|u^2+\mathbbm{1}(v=1)\leq 1,-1\leq u\leq 1, 0\leq v\leq 1\}=q(0)+\sqrt{2\alpha t}.
\end{align*}  
\indent As a last example, we observe the scenario that $\probability*{\StartRVIndep{i}>v}= \exp(-vL(v))$, thus $h(v)=v$. Then, 
\begin{align*}
g(t,q(0))=&\sup_{(u,v)}\{\sqrt{2\alpha t}u+q(0)v|u^2+v\leq 1,0\leq u\leq 1, 0\leq v\leq 1\} \nonumber\\
=&\left(q(0)+\frac{\alpha t}{2q(0)}\right)\mathbbm{1}\left(t<\frac{2q(0)^2}{\alpha}\right)+\sqrt{2\alpha t}\mathbbm{1}\left(t\geq\frac{2q(0)^2}{\alpha}\right).
\end{align*}
We would like to give some extra attention to the case where $q(0)=\frac{\alpha}{(2\beta)}$. Then, it is not difficult to see that $q(t)\equiv \frac{\alpha}{(2\beta)}$. Thus, for these choices of $h(v)$ and $q(0)$, the system starts and stays in steady state. One can show that this limit is only obtained for $h(v)=v$, so this gives us some information on the joint steady-state distribution of \emph{all} the queue lengths in the fork-join system.

\begin{figure}[H]
\minipage{0.32\textwidth} \includegraphics[width=\linewidth]{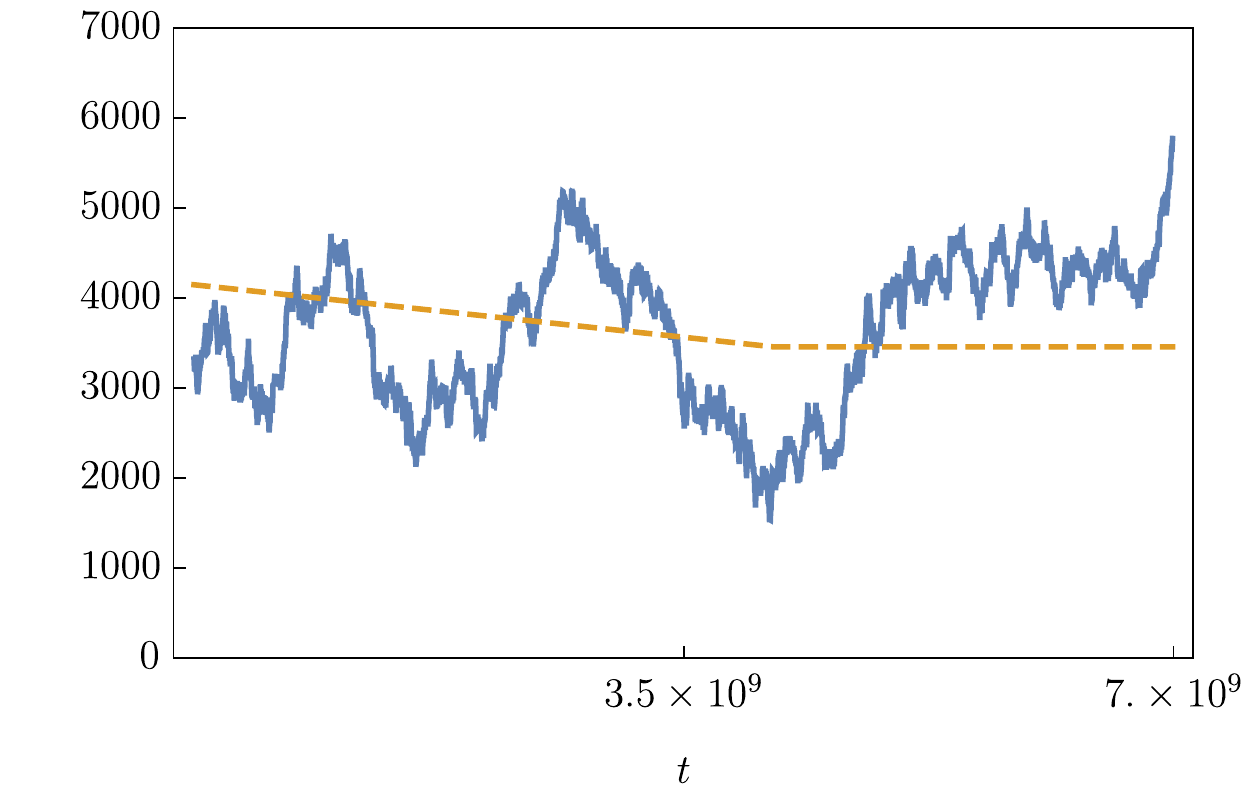}
  \subcaption{$\alpha=1,\beta=1, q(0)=0.6$}\label{fig:a1b1x06}
  \endminipage\hfill
\minipage{0.32\textwidth}
  \includegraphics[width=\linewidth]{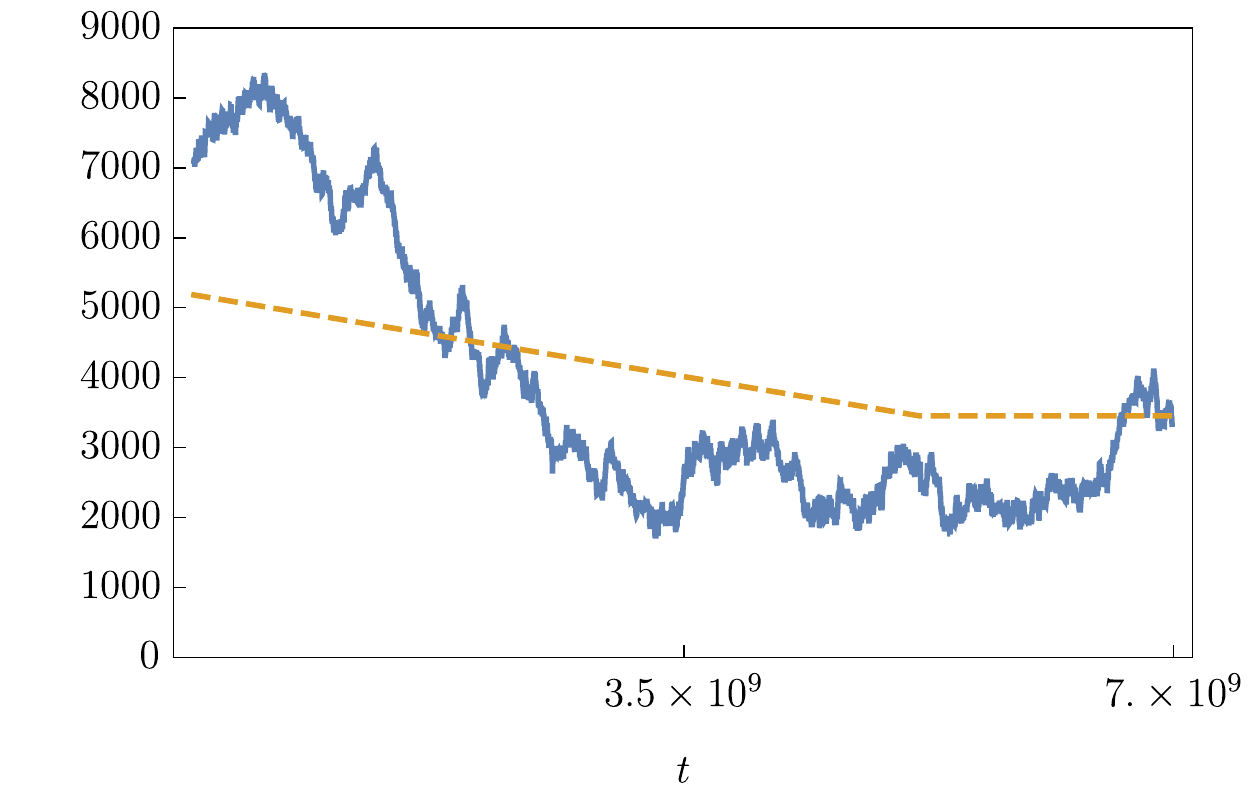}
  \subcaption{$\alpha=1,\beta=1, q(0)=0.75$}\label{fig:a1b1x075}
\endminipage\hfill
\minipage{0.32\textwidth}%
 \includegraphics[width=\linewidth]{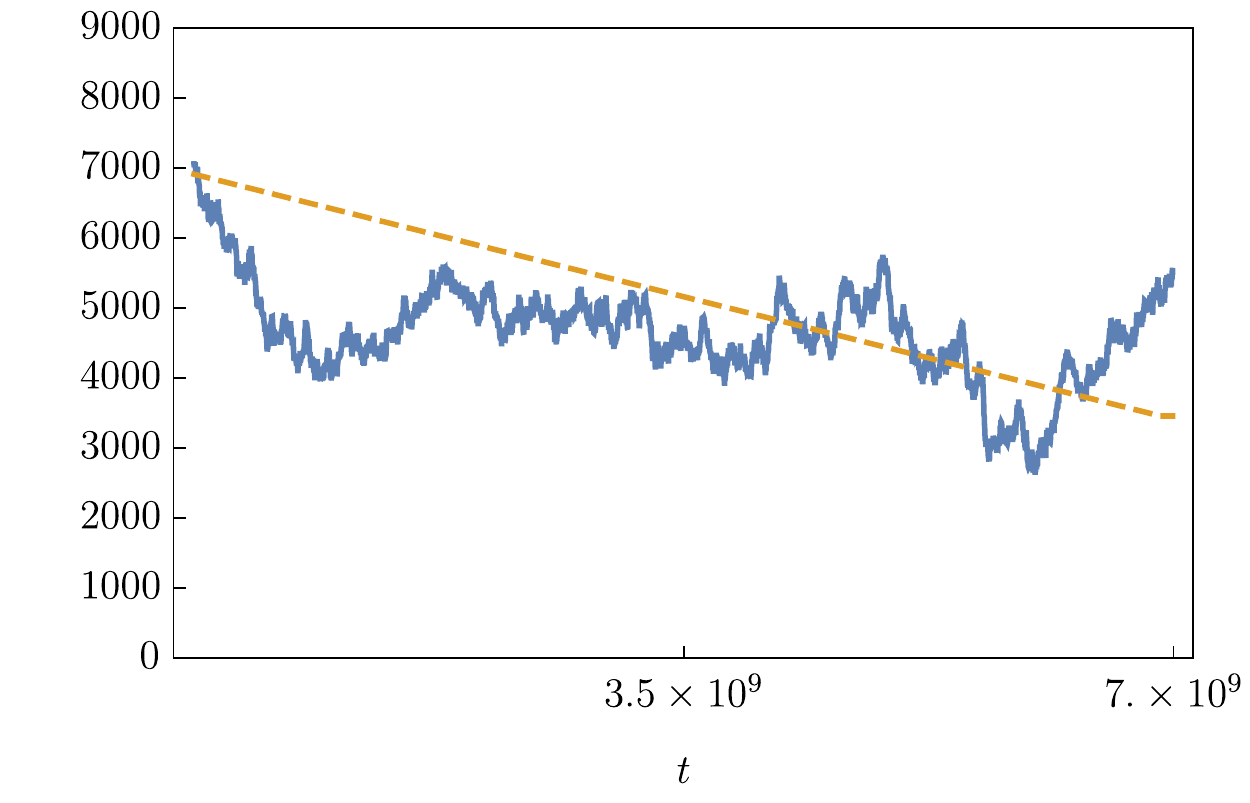}
  \subcaption{$\alpha=1,\beta=1, q(0)=1$}\label{fig:a1b1x1}
\endminipage\hfill
\minipage{0.32\textwidth}%
 \includegraphics[width=\linewidth]{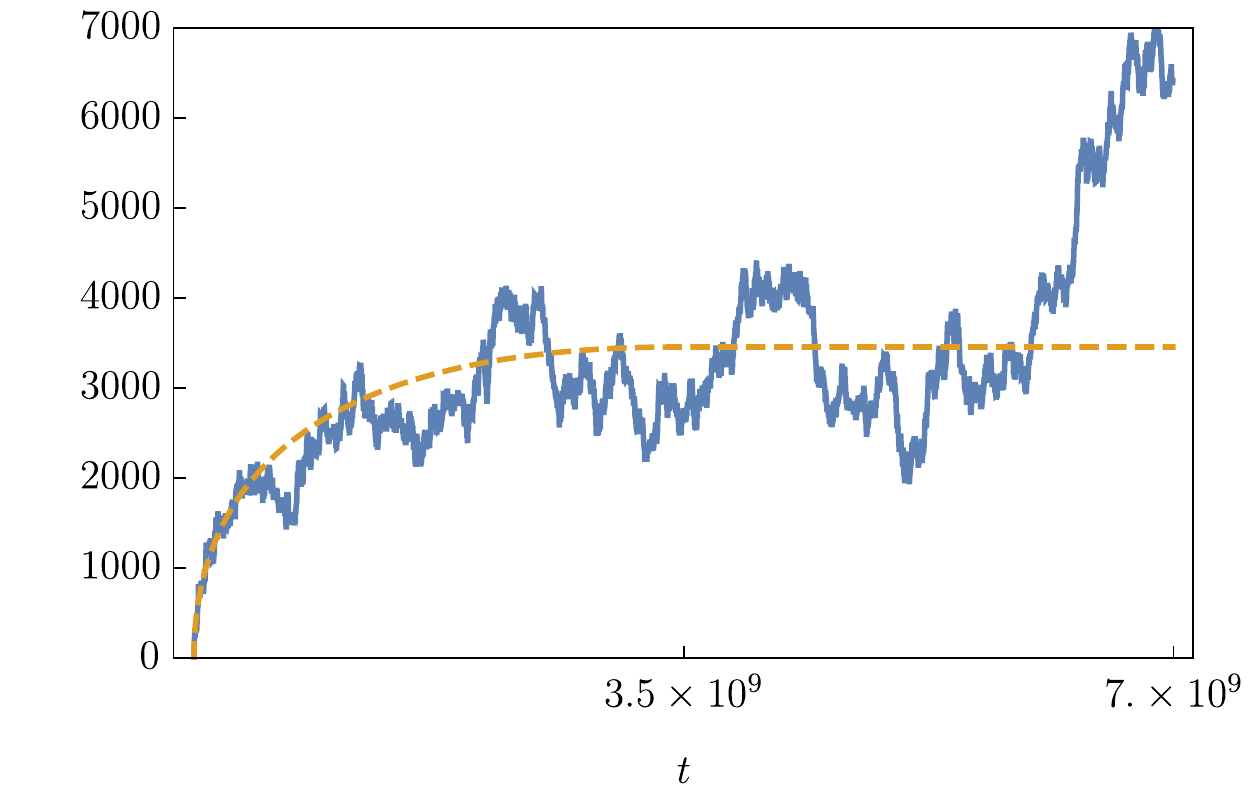}
  \subcaption{$\alpha=1,\beta=1, q(0)=0$}\label{fig:a1b1x0}
\endminipage\hfill
\minipage{0.32\textwidth}%
 \includegraphics[width=\linewidth]{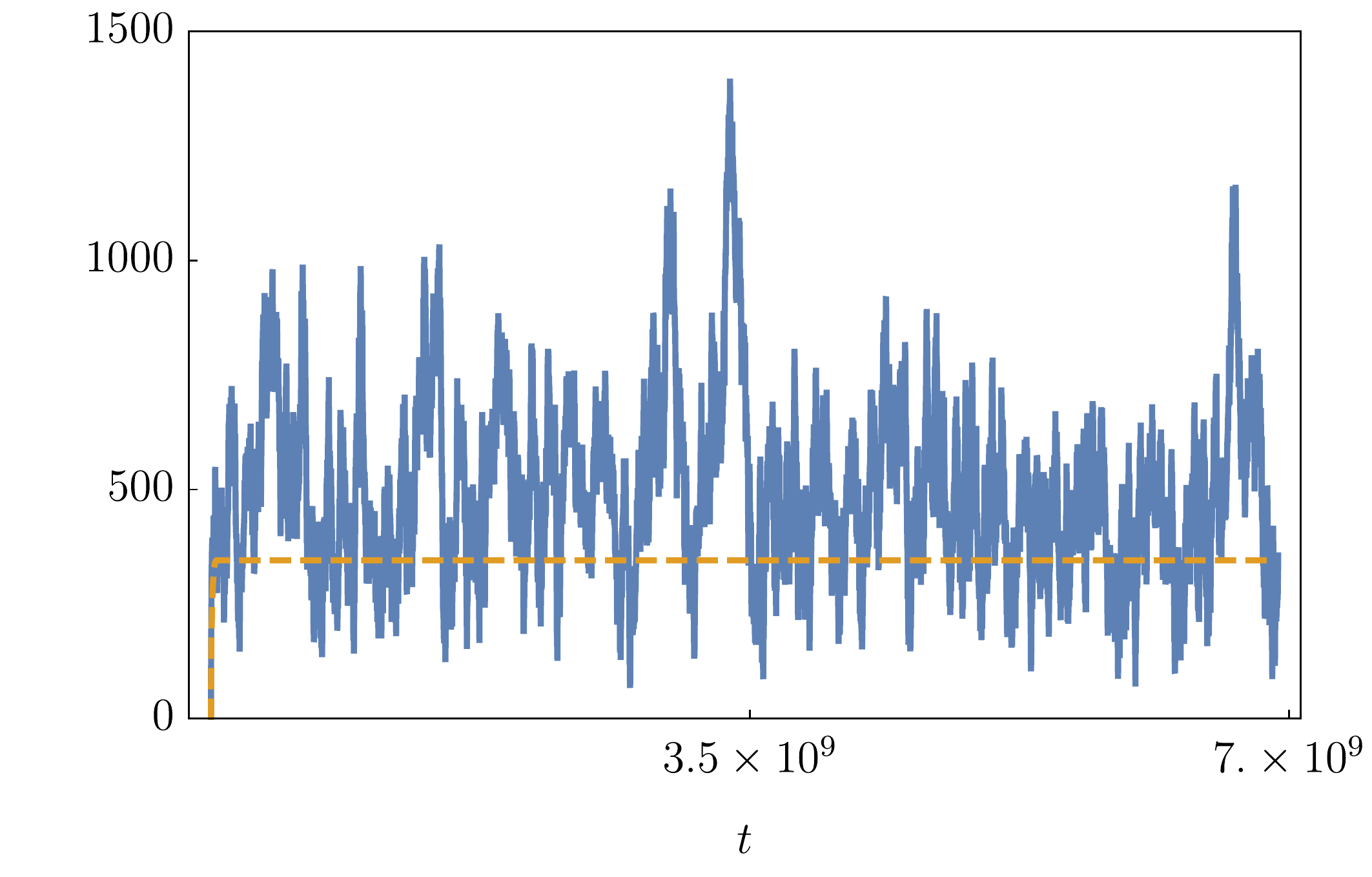}
  \subcaption{$\alpha=1,\beta=10, q(0)=0$}\label{fig:a1b10x0}
\endminipage\hfill
\minipage{0.32\textwidth}%
 \includegraphics[width=\linewidth]{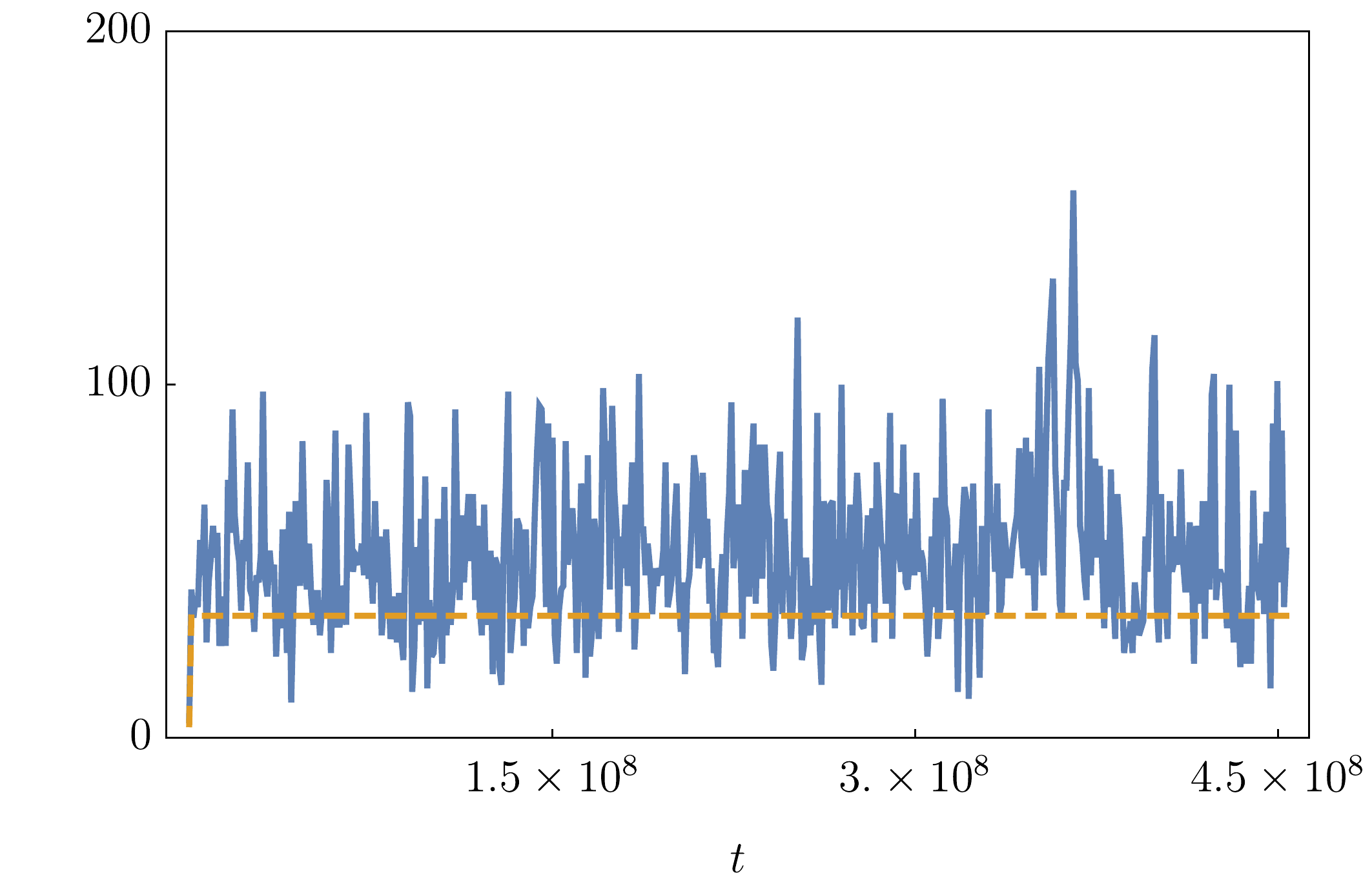}
  \subcaption{$\alpha=1,\beta=100, q(0)=0$}\label{fig:a1b100x0}
\endminipage\hfill
  \caption{Maximum queue length and fluid limit approximation (Thm. \ref{thm: fluid limit}) for $N=1000$}
\label{fig:6fluidplots}
\end{figure}

Now, we turn to some numerical examples. In Figure \ref{fig:6fluidplots}, the simulated maximum queue length is plotted together with the scaled fluid limit $N\log Nq(\frac{t}{(N^3\log N)})$, with $q$ given in Theorem \ref{thm: fluid limit}, and $N=1000$. The queue lengths at time zero in Figures \ref{fig:a1b1x06}, \ref{fig:a1b1x075} and \ref{fig:a1b1x1} are exponentially distributed. These figures show that for $N=1000$, the maximum queue length is not close to its fluid limit.
As these figures show, for $N=1000$, the variance of the maximum queue length is still high. We could however give some heuristic arguments why these results are not very accurate. As mentioned before, we have that
\begin{align*}
\frac{\Arrival{tN^3\log N}-\left(1-\frac{\alpha}{N}\right)(tN^3\log N)}{N\log N}\LimitP -\beta t\text{ as $N\to\infty$},
\end{align*}
which is one building block of the fluid limit.\\ For $\frac{(\Arrival{tN^3\log N}-\left(1-\frac{\alpha}{N}\right)tN^3\log N)}{(N\log N)}$, we can compute the standard deviation. We have for $\alpha=\beta=t=1$ and $N=1000$ that
\begin{align*}
\sqrt{\text{Var}\left(\Arrival{tN^3\log N}-\left(1-\frac{\alpha}{N}\right)(tN^3\log N)\right)}=&\sqrt{\left(1-\frac{\alpha}{N}-\frac{\beta}{N^2}\right)\left(\frac{\alpha}{N}+\frac{\beta}{N^2}\right)\lfloor{tN^3\log N\rfloor}}\nonumber\\
=&2628.26.
\end{align*} 
This is of the order of magnitude of the errors that we see in the figures. 

Another way of seeing that there is a significant deviation is by looking at $\max_{i\leq N}\left(\left(1-\frac{\alpha}{N}\right)tN^3\log N-\Service{i}{tN^3\log N}\right)$. As mentioned in Section \ref{subsec: shape fluid limit}, we have that 
\begin{align*}
\frac{\left(1-\frac{\alpha}{N}\right)tN^3\log N-\Service{i}{tN^3\log N}}{N\sqrt{\log N}}\overset{d}{\approx} \vartheta_i,
\end{align*} 
with $\vartheta_i\sim\mathcal{N}(0,\alpha t)$. Thus, this means that 
\begin{align*}
\max_{i\leq N}\left(\left(1-\frac{\alpha}{N}\right)tN^3\log N-\Service{i}{tN^3\log N}\right)\overset{d}{\approx} \max_{i\leq N}\vartheta_iN\sqrt{\log N}.
\end{align*}
When we choose $N=1000$, $\alpha=t=1$, and simulate enough samples of $ \max_{i\leq N}\vartheta_iN\sqrt{\log N}$, we observe a standard deviation which is higher than 900.

In Figures \ref{fig:a1b1x06}, \ref{fig:a1b1x075} and \ref{fig:a1b1x1}, the high standard deviation is also caused by the distribution of the number of jobs at time 0. For example, for $E_i\sim \text{Exp}(\frac{1}{N})$, i.i.d.\ for all $i$, and $N=1000$, we have that $\sqrt{\text{Var}\left(\max_{i\leq N}E_i\right)}=1282.16$, so this is also of the order of magnitude of the errors that we see.

\begin{figure}[H]
\minipage{0.32\textwidth} \includegraphics[width=\linewidth]{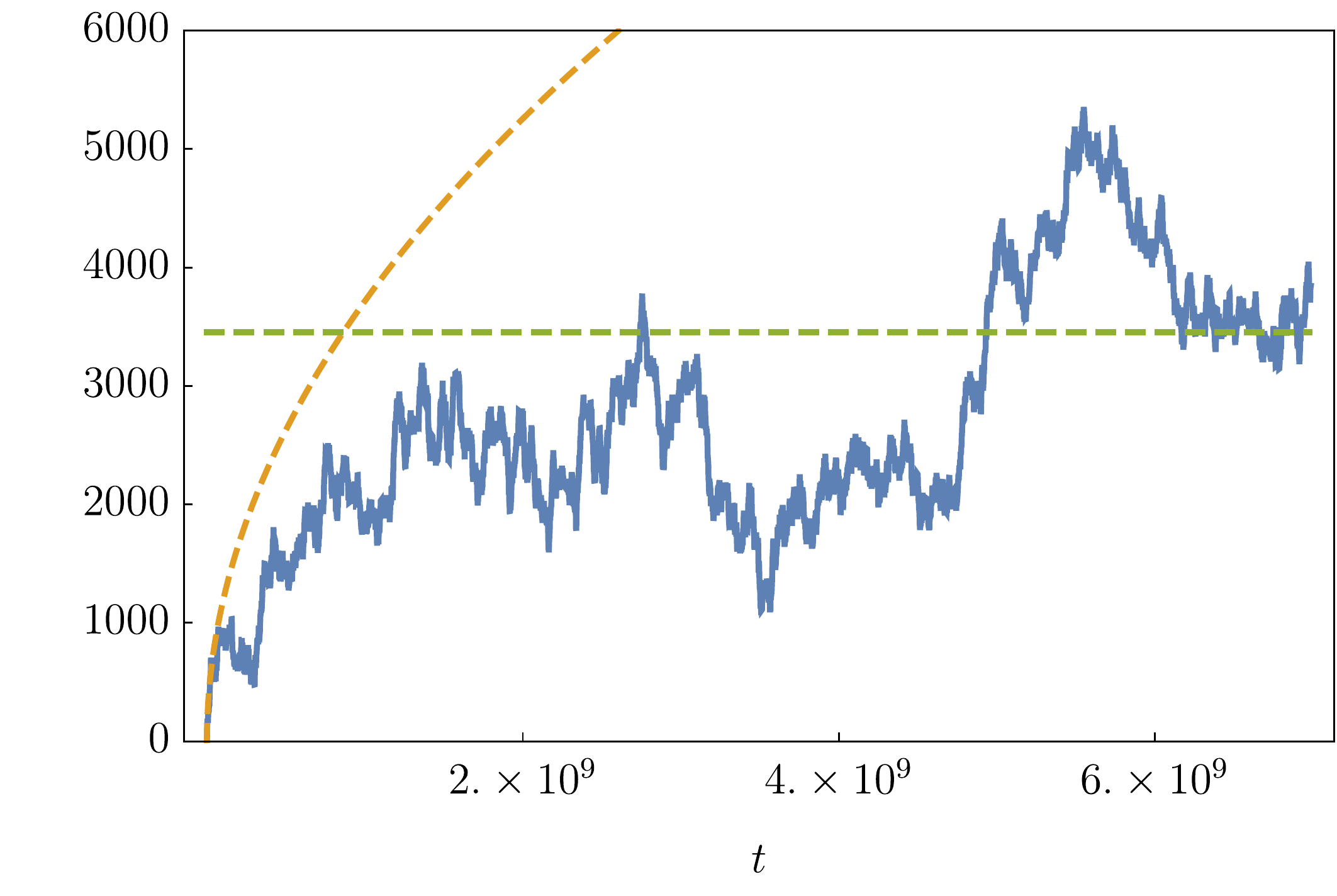}
  \subcaption{$\alpha=1,\beta=1$}\label{fig:a1b1}
  \endminipage\hfill
\minipage{0.32\textwidth}
  \includegraphics[width=\linewidth]{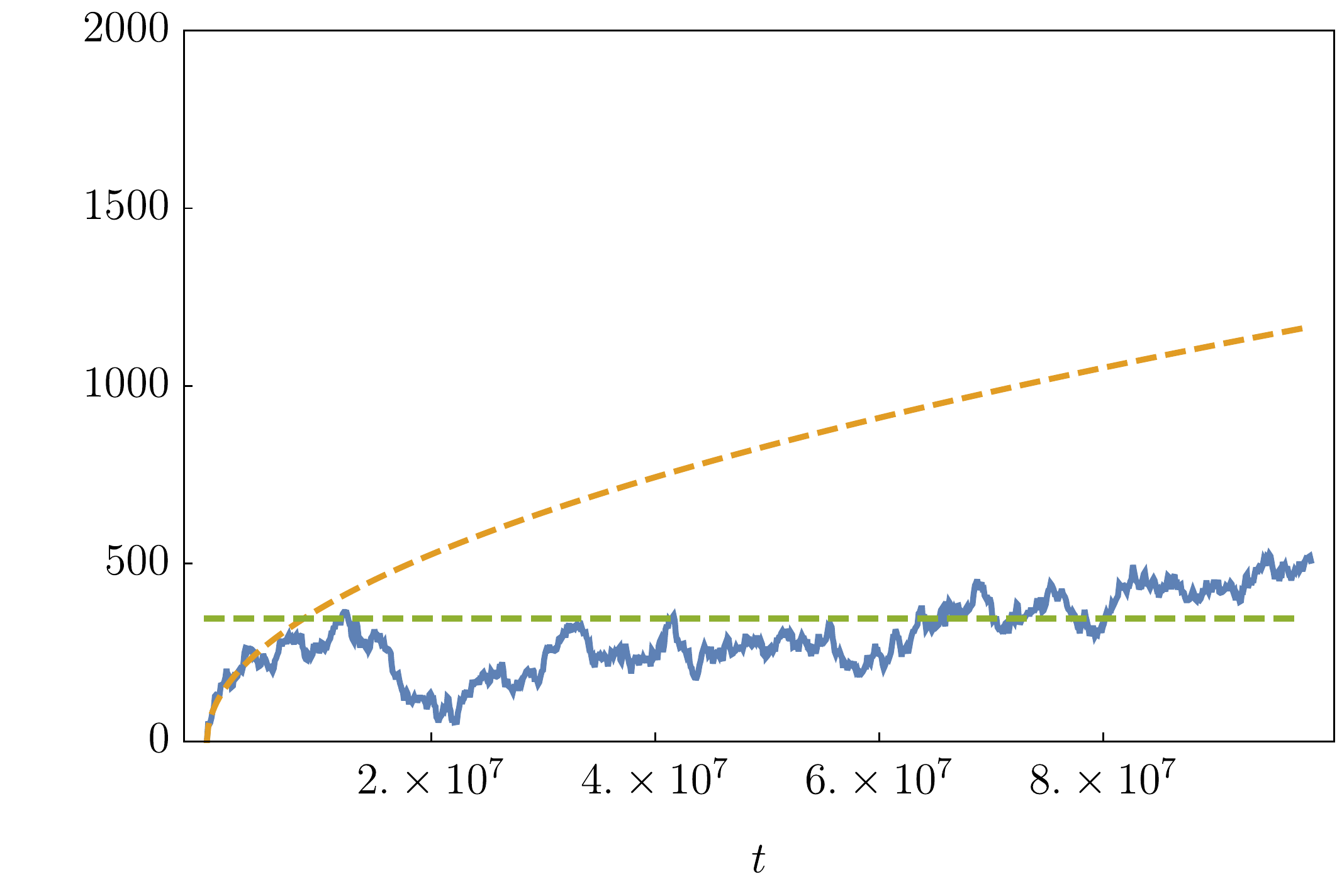}
  \subcaption{$\alpha=1,\beta=10$}\label{fig:a1b10}
\endminipage\hfill
\minipage{0.32\textwidth}%
 \includegraphics[width=\linewidth]{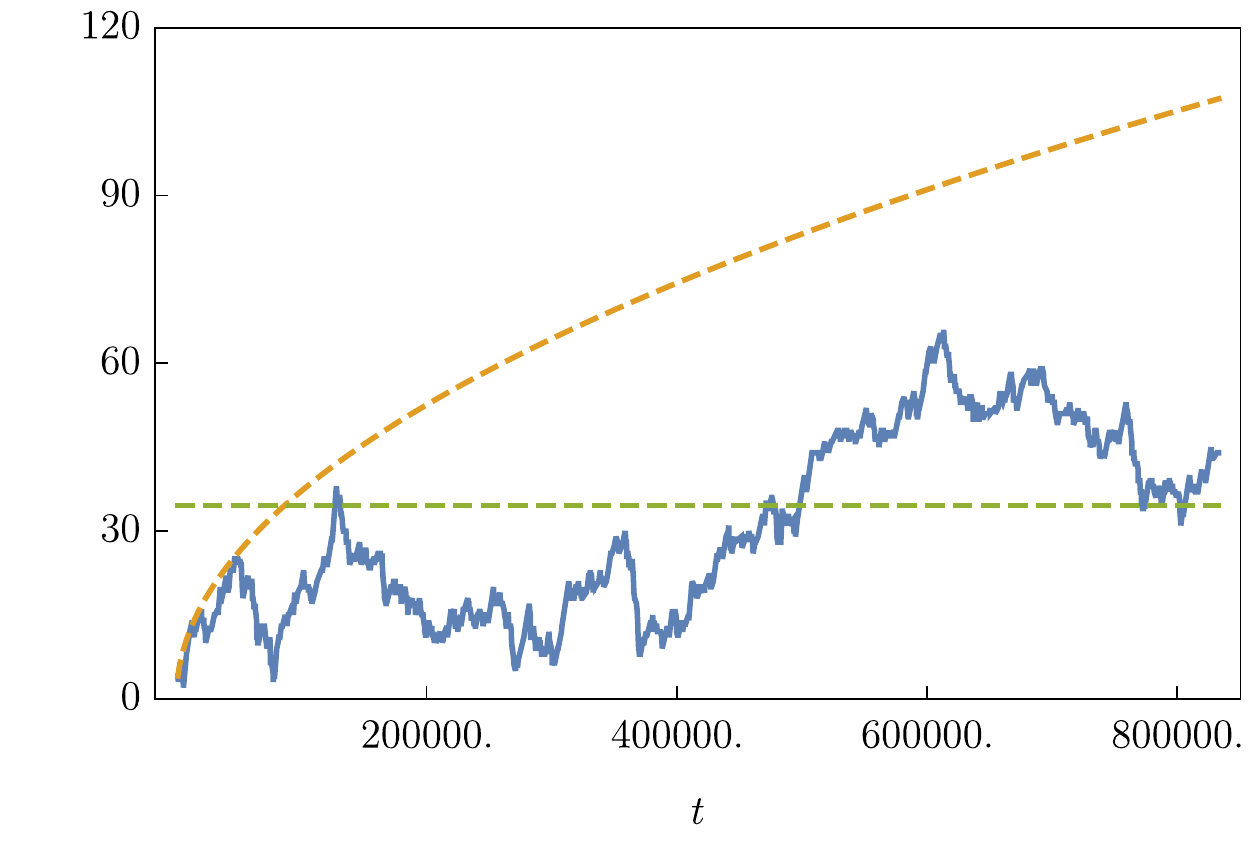}
  \subcaption{$\alpha=1,\beta=100$}\label{fig:a1b100}
\endminipage
  \caption{Maximum queue length, fluid limit approximation (Prop. \ref{cor: other fluid limit}) and steady-state approximation for $N=1000$}
\label{fig:3fluidplots}

\end{figure}

\begin{figure}[H]
\minipage{0.32\textwidth} \includegraphics[width=\linewidth]{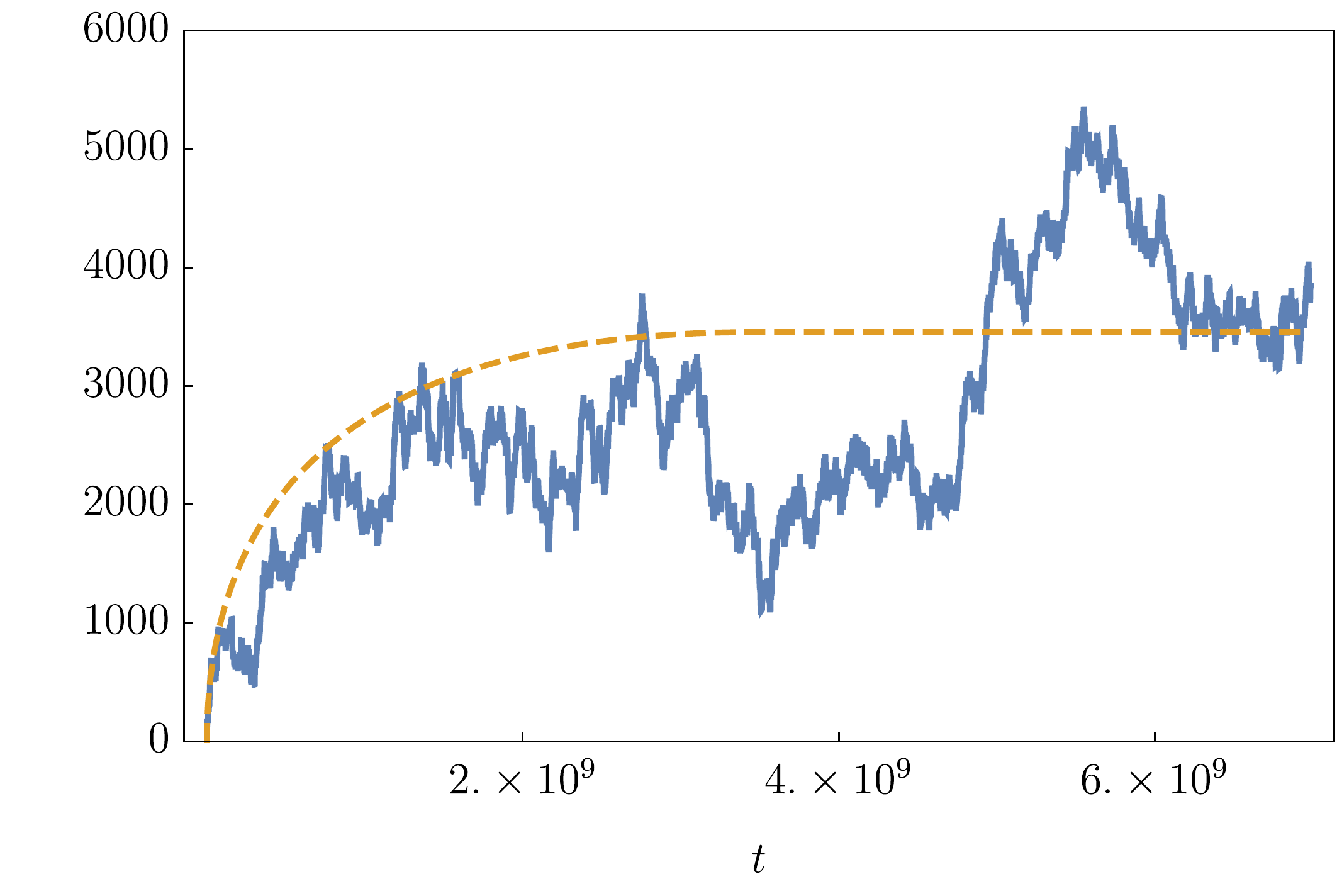}
  \subcaption{$\alpha=1,\beta=1$}\label{fig:a1b1New}
  \endminipage\hfill
\minipage{0.32\textwidth}
  \includegraphics[width=\linewidth]{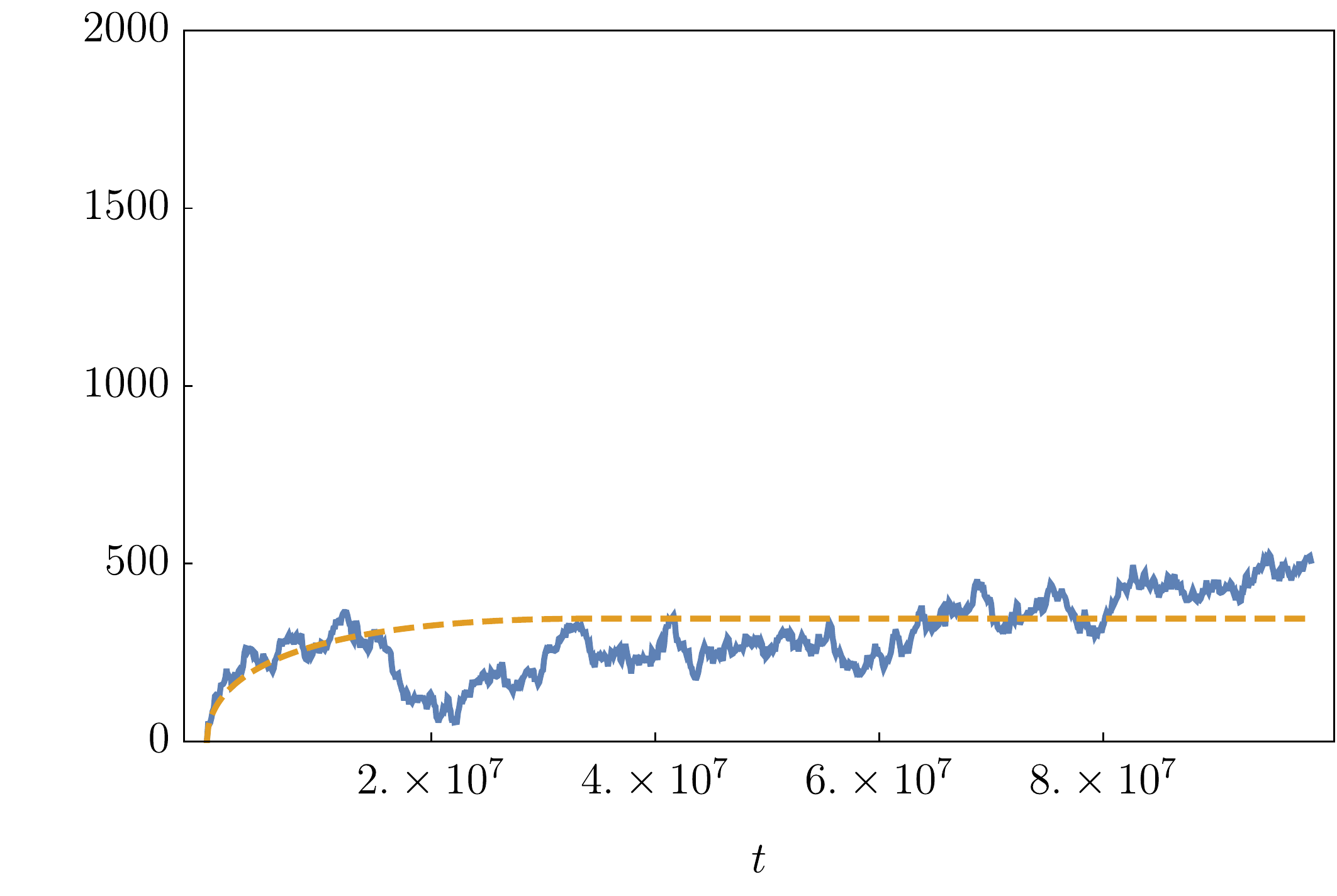}
  \subcaption{$\alpha=1,\beta=10$}\label{fig:a1b10New}
\endminipage\hfill
\minipage{0.32\textwidth}%
 \includegraphics[width=\linewidth]{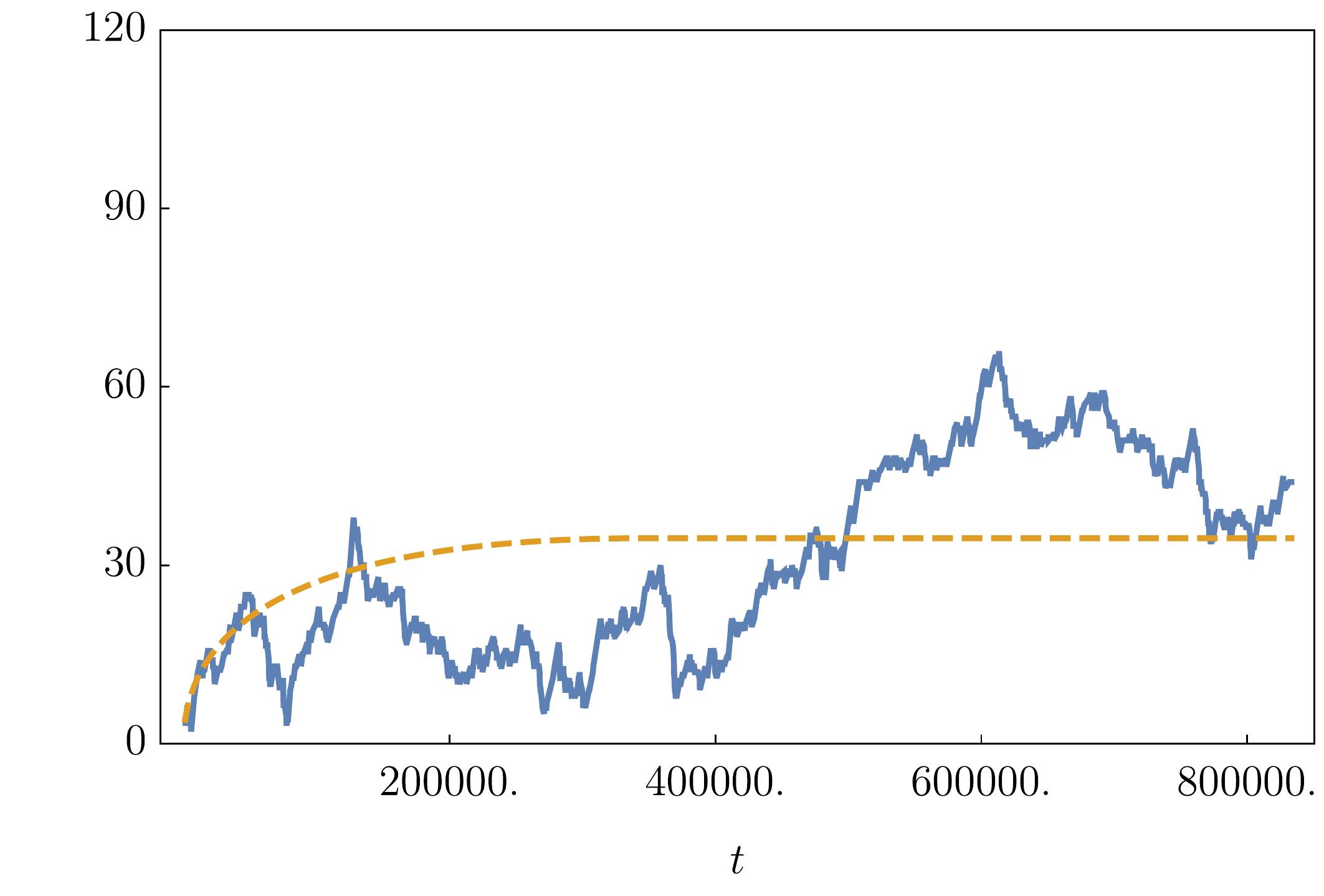}
  \subcaption{$\alpha=1,\beta=100$}\label{fig:a1b100New}
\endminipage
  \caption{Maximum queue length and fluid limit approximation (Thm. \ref{thm: fluid limit}) for $N=1000$}
\label{fig:3newfluidplots}

\end{figure}

\begin{figure}[H]
\minipage{0.5\textwidth} \includegraphics[width=\linewidth]{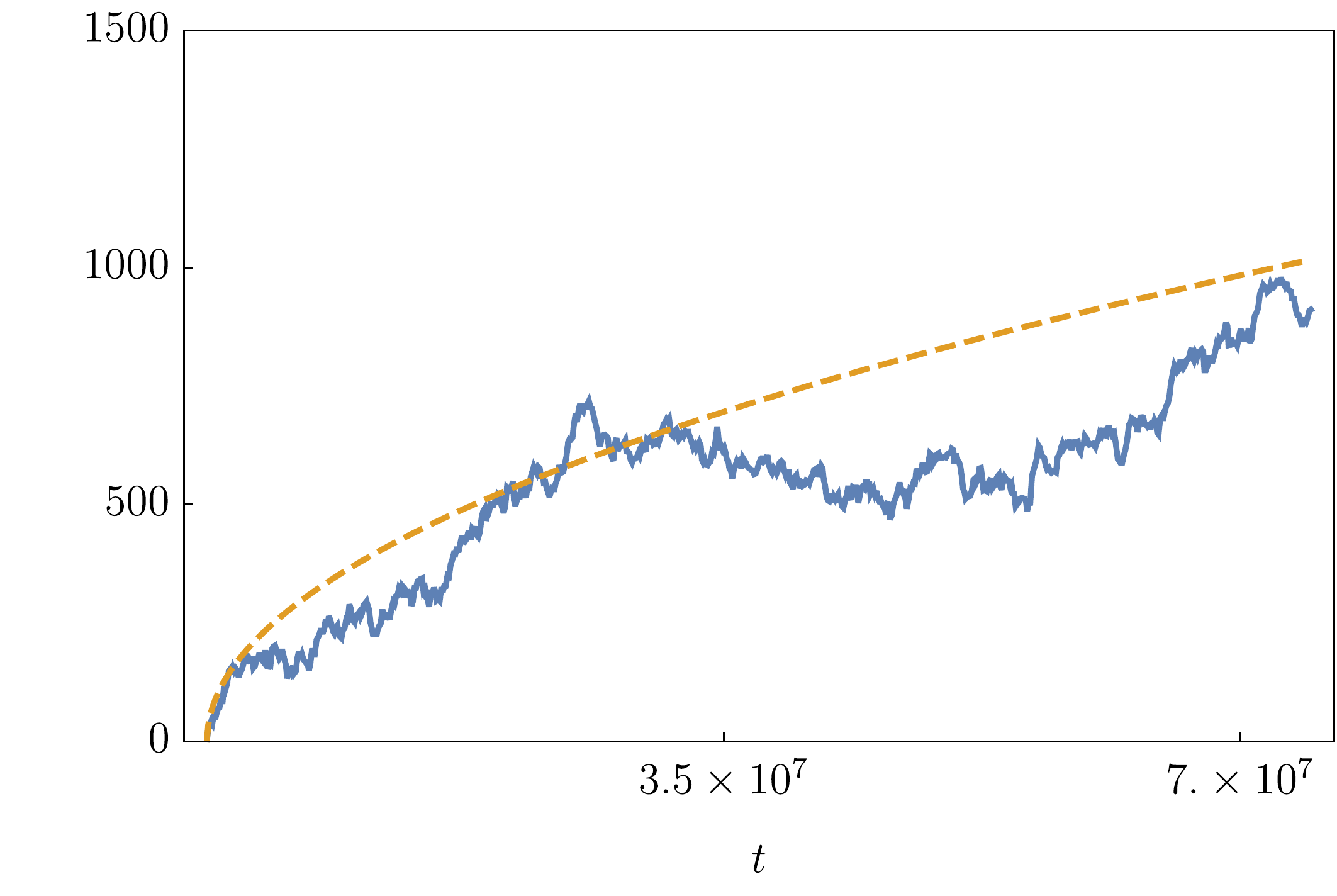}
  \subcaption{$\alpha=1,\beta=1$}\label{fig:a1b1OldSmall}
  \endminipage\hfill
\minipage{0.5\textwidth}
  \includegraphics[width=\linewidth]{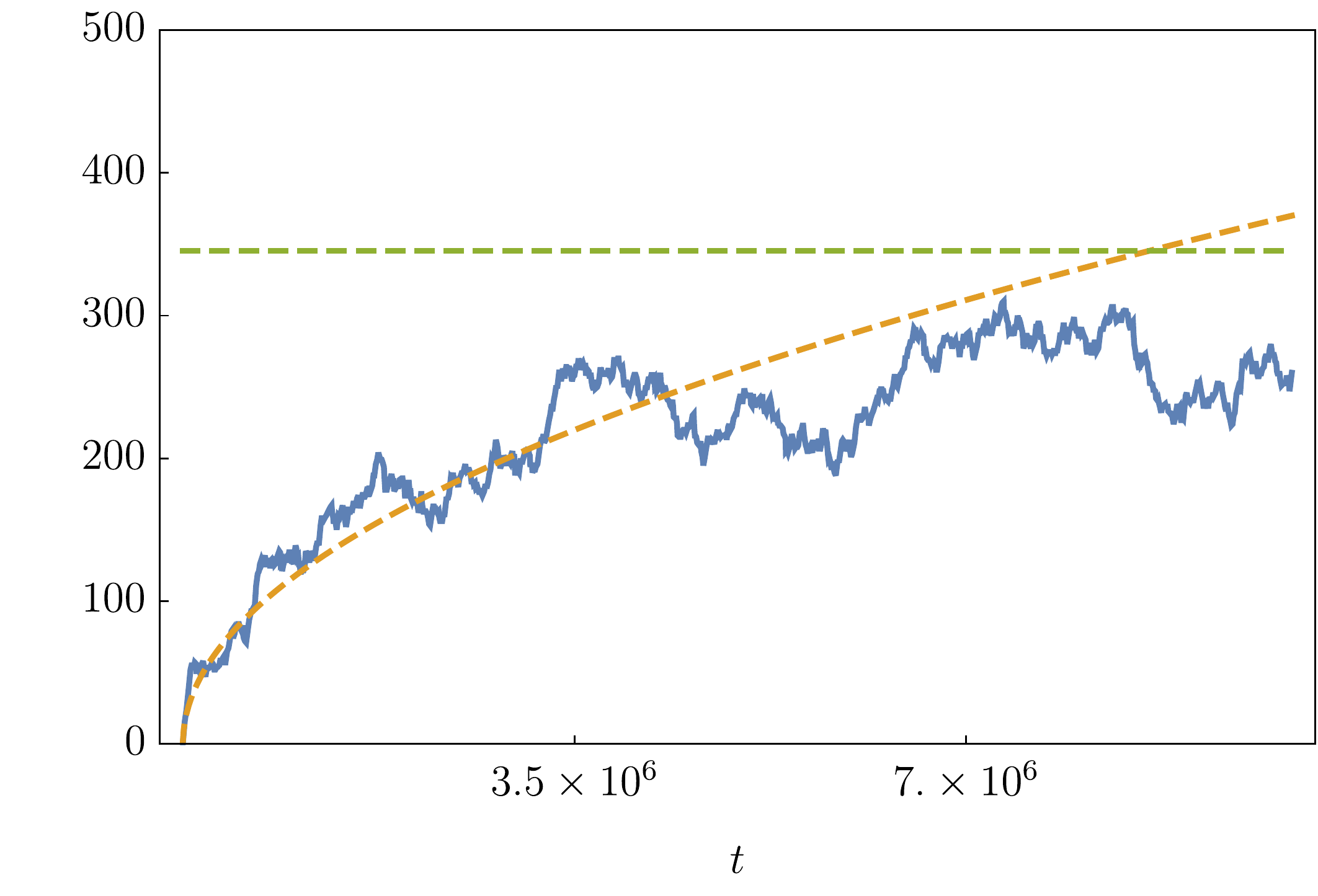}
  \subcaption{$\alpha=1,\beta=10$}\label{fig:a1b10OldSmall}
\endminipage 
\caption{Maximum queue length, fluid limit approximation (Prop. \ref{cor: other fluid limit}) and steady-state approximation for $N=1000$}
\label{fig:2oldfluidplotsSmall}

\end{figure}

As mentioned, one can prove fluid limits under several temporal and spatial scalings. In Figure \ref{fig:3fluidplots}, the maximum queue length is plotted against the rescaled fluid limit given in Proposition \ref{cor: other fluid limit}, which is in orange, and the rescaled steady-state limit, which is in green. In these plots, $N=1000$. The rescaled fluid limit is $\sqrt{\frac{2\alpha t}{N^3}}N\sqrt{\log N}$, and the rescaled steady-state limit satisfies $\frac{\alpha}{(2\beta)}N\log N$.
When we observe Figure \ref{fig:3fluidplots}, we see that for small time instances, the maximum queue length follows the fluid limit described in Proposition \ref{cor: other fluid limit} with a negligible deviation, and we also see that, from the point that the fluid limit and steady state have intersected, the maximum queue length follows the steady state, though with a significant deviation. This latter behavior can be very well explained when we plot the same maximum queue lengths together with the fluid limit in Theorem \ref{thm: fluid limit}, this is shown in Figure \ref{fig:3newfluidplots}.
In Figure \ref{fig:2oldfluidplotsSmall}, we zoom in on the graphs given in Figures \ref{fig:a1b1} and \ref{fig:a1b10}. As these figures show, for small time instances, the maximum queue length follows the fluid limit described in Proposition \ref{cor: other fluid limit} quite well. Again, we can heuristically explain the deviations by approximating the maximum queue length with $\sqrt{1/N^3}N\max_{i\leq N}\vartheta_i$, with $\vartheta_i\sim\mathcal{N}(0,\alpha t)$, i.i.d.\ For $\alpha=1$, and $t=7\cdot 10^7$, simulations show that this approximation has a standard deviation around 95, and for $t=7\cdot 10^6$, we get a standard deviation around 30, this is of the order of magnitude of the errors in Figures \ref{fig:a1b1OldSmall} and \ref{fig:a1b10OldSmall}, respectively.

\section[Conclusion]{Conclusion.}\label{sec: future work}
In this paper, we analyzed a fork-join queue with $N$ servers in heavy traffic. We considered the case of nearly deterministic arrivals and service times, and we derived a fluid limit of the maximum queue length, in Theorem \ref{thm: fluid limit}, as $N$ grows large. 

Furthermore, we assumed delays to be memoryless. However, we are confident that these results can be extended to nearly deterministic settings where the delays have general distributions. Another, but less straightforward extension of this result, would be to assume arrival and service processes that are not Markovian.

Moreover, as the figures in Section \ref{subsec: examples and numerics} show, it should be possible to derive a more refined limit. Therefore, it is interesting to look at the second-order convergence of the maximum queue length. We are currently exploring this for the system in steady state. In other words, we try to gain more insight into the process by finding a convergence result of $\frac{\MaxQueueLength{\infty}}{N}-\frac{\alpha}{(2\beta)}\log N$. For the process limit, proving a second-order convergence result is much harder and more technical because the scaled maximum of $N$ independent Brownian motions converges to a Brown-Resnick process \cite{brown1977extreme}.  
\section[Proofs]{Proofs.}\label{sec: Proofs}
In this section, we prove Theorem \ref{thm: fluid limit}. Since each server has the same arrival process, the queue lengths are dependent. The general idea of proving Theorem \ref{thm: fluid limit} is to approximate the scaled centralized service process in \eqref{eq: normal like random variable} by a normally distributed random variable. We can use extreme value theory to prove convergence of the maximum of these normally distributed random variables in probability. By using the non-uniform version of the Berry-Ess\' een theorem; cf.\ \cite{michel1976constant}, we show that the convergence result of the original process is the same as the convergence result with normally distributed random variables. Furthermore, we prove convergence of the part involving non-zero starting points. This gives us the pointwise convergence of the process, which we prove in Section \ref{subsubsec: pointwise convergence}. In this section, we also prove convergence of the finite-dimensional distributions. Finally, we prove in Section \ref{subsubsec: tightness} that the process is tight. These three results together prove the theorem.
\subsection[Definitions]{Definitions.}
For the sake of notation, we use the expressions given in Definition \ref{def: tilde arrival service} to prove the tightness.
\begin{definition}\label{def: tilde arrival service}

We define the random walk $\TildeQueue{i}{n}$ as
\begin{align}
\TildeQueue{i}{n}=\frac{\TildeArrival{n}+\TildeService{i}{n}}{\log N},
\end{align}
where
\begin{align}
\TildeArrival{n}=\frac{\Arrival{n}}{N}-\left(1-\frac{\alpha}{N}\right)\frac{\lfloor n\rfloor}{N},
\end{align}
and
\begin{align}
\TildeService{i}{n}=-\frac{\Service{i}{n}}{N}+\left(1-\frac{\alpha}{N}\right)\frac{\lfloor n\rfloor}{N}.
\end{align}
Furthermore,
\begin{align}\label{eq: normal like random variable}
\MaximumProcess{i}{t}=\frac{\TildeService{i}{{t N^3\log N}}}{\sqrt{\alpha t(1-\frac{\alpha}{N})\log N}}\frac{\sqrt{tN^3\log N}}{\sqrt{\lfloor{tN^3\log N\rfloor}}},
\end{align}
with $\Arrival{n}$ and $\Service{i}{n}$ given in Definitions \ref{def: arrival process} and Definition \ref{def: service process i server} respectively.
\end{definition}
As mentioned in Section \ref{subsec: shape fluid limit}, when $\MaxQueueLength{0}=0$, the quantity in \eqref{eq: maxqueue length expression} simplifies to
\begin{align*}
&\FluidProcessTwo{t}\nonumber\\
=&\max_{i\leq N}\sup_{0\leq s\leq t}\frac{\Big(\Arrival{tN^3\log N}-\Arrival{sN^3\log N}\Big)-\left(\Service{i}{tN^3\log N}-\Service{i}{sN^3\log N}\right)}{N\log N}.
\end{align*} Consequently, we can rewrite 
\begin{align}\label{eq: Q expressed in tilde Q}
&\FluidProcessTwo{t}\nonumber\\
=&\max_{i\leq N}\sup_{0\leq r\leq t}\frac{\TildeArrival{{tN^3\log N}}-\TildeArrival{{rN^3\log N}}+\TildeService{i}{{tN^3\log N}}-\TildeService{i}{{rN^3\log N}}}{\log N}\nonumber\\
=&\max_{i\leq N}\sup_{0\leq r\leq t} \left(\TildeQueue{i}{tN^3\log N}-\TildeQueue{i}{rN^3\log N}\right).
\end{align}
\subsection[Useful lemmas]{Useful lemmas.}\label{subsubsec: useful lemmas}
In order to prove Theorem \ref{thm: fluid limit}, a few preliminary results are needed. As stated in Definition \ref{def: tilde arrival service}, we can write $\TildeQueue{i}{n}$ as 
\begin{align*}
\frac{\TildeArrival{n}+\TildeService{i}{n}}{\log N}.
\end{align*} 
Observe that $\TildeArrival{n}$ does not depend on $i$, while $\TildeService{i}{n}$ does. Hence, it is intuitively clear that $\TildeArrival{n}$ pays no contribution to the maximum queue length. Therefore, in order to prove the pointwise convergence of the maximum queue length, we need to analyze $\frac{\TildeService{i}{n}}{\log N}$.
Specifically, we use the fact that 
\begin{align*}
\MaximumProcess{i}{t}\LimitD Z\text{ as $N\to\infty$},
\end{align*}
with $Z$ a standard normal random variable, which can be shown by the central limit theorem. We can use this result to approximate the maximum queue length because we know that the scaled maximum of $N$ independent and normally distributed random variables converges to a Gumbel distributed random variable. 
To prove the tightness of the maximum queue length, we have to prove that
\begin{align}\label{eq: billingsley condition}
\lim_{\delta\downarrow 0}\limsup_{N\to\infty}\frac{1}{\delta}\probability*{\sup_{t\leq s\leq t+\delta}\left|\FluidProcessTwo{s}-\FluidProcessTwo{t}\right|>\epsilon}=0.
\end{align}
In Lemma \ref{lem: upper bound |Qs-Qt|}, a useful upper bound for the absolute value in \eqref{eq: billingsley condition} is obtained, which we use to prove the tightness of the process. 
\begin{lemma}\label{lem: upper bound |Qs-Qt|}

For $t>0$, $\delta>0$ and $\MaxQueueLength{0}=0$, we have that
\begin{align}
&\sup_{t\leq s\leq t+\delta}\left|\FluidProcessTwo{s}-\FluidProcessTwo{t}\right|\nonumber\\
\leq &\sup_{t\leq s\leq t+\delta}\max_{i\leq N}\left(\TildeQueue{i}{sN^3\log N}-\TildeQueue{i}{tN^3\log N}\right)\nonumber\\
+&2\sup_{t\leq s\leq t+\delta}\max_{i\leq N}\left(\TildeQueue{i}{tN^3\log N}-\TildeQueue{i}{sN^3\log N}\right).
\end{align}
\end{lemma}
In our proofs, we use the fact that $\MaximumProcess{i}{t}$ converges in distribution to a normally distributed random variable. To be able to use this convergence result, we prove an upper bound of the convergence rate in Lemma \ref{lem: distance cdf normal}.
\begin{lemma}\label{lem: distance cdf normal}
For $t>0$, we have that an upper bound of the rate of convergence of $\frac{\pm\TildeService{i}{{t N^3\log N}}\sqrt{tN^3\log N}}{\sqrt{\alpha t(1-\frac{\alpha}{N})\log N\lfloor{tN^3\log N\rfloor}}}$ to a standard normal random variable is given by
\begin{align}\label{eq: inequality michel}
\left|\probability*{\MaximumProcess{i}{t}<y}-\Phi(y)\right|\leq \DifferenceNormal,
\end{align}
with $c_{t}>0$.
\end{lemma}
Lemma \ref{lem: distance cdf normal} follows from the main result in \cite{michel1976constant}, where the author proves the non-uniform Berry-Ess\' een inequality. To prove tightness, we need the following lemma:
\begin{lemma}\label{lem: convergence moments 1 2 4}
For $t>0$,
\begin{align}\label{eq: convergence expectation first moment to constant}
\limsup_{N\to\infty}\expect*{\max\left(\max_{i\leq N} \frac{\pm\TildeService{i}{t N^3\log N}}{\log N},0\right)^{\frac{5}{2}}}\leq \left(2\alpha t\right)^{\frac{5}{4}}.
\end{align}
\end{lemma}
In order to prove pointwise convergence of the starting position, we show in Lemma \ref{lem: pointwise convergence startposition} that 
\begin{align*}
\max_{i\leq N}\left(\frac{\TildeService{i}{{t N^3\log N}}}{\log N}+\frac{\QueueZero{i}}{N\log N}\right)\approx \max_{i\leq N}\left(\frac{\sqrt{\alpha t}X_i}{\sqrt{\log N}}+\frac{\QueueZero{i}}{N\log N}\right),
\end{align*}
with $X_i\sim\mathcal{N}(0,1)$, as $N$ is large.\\
In Lemma \ref{lem: approximation starting point convergence}, we prove the convergence of
$\max_{i\leq N}\left(\frac{\sqrt{\alpha t}X_i}{\sqrt{\log N}}+\frac{\QueueZero{i}}{(N\log N)}\right)$.
\begin{lemma}[Pointwise convergence approximation starting position]\label{lem: approximation starting point convergence}
\begin{align*}
\max_{i\leq N}\left(\frac{\sqrt{\alpha t}X_i}{\sqrt{\log N}}+\frac{Q^{(N)}_i(0)}{N\log N}\right)
\LimitP g(t,q(0))\text{ as $N\to\infty$},
\end{align*}
with $X_i\sim\mathcal{N}(0,1)$ i.i.d.\ and the function $g$ as given in Theorem \ref{thm: fluid limit}.
\end{lemma}
The proofs of Lemmas \ref{lem: upper bound |Qs-Qt|}, \ref{lem: distance cdf normal}, \ref{lem: convergence moments 1 2 4}, and \ref{lem: approximation starting point convergence} can be found in Appendix \ref{app: proof lemma convergence fourth moment}. Lemma \ref{lem: approximation starting point convergence} follows from Lemma \ref{lem: extreme value convergence}, where a more general result is proven on $\max_{i\leq N}\sum_{j=1}^k\frac{Y^{(j)}_i}{a_N^{(j)}}$.
\subsection[Pointwise convergence]{Pointwise convergence.}\label{subsubsec: pointwise convergence} 
In this section, we prove pointwise convergence of the scaled maximum queue length appearing in Theorem \ref{thm: fluid limit}.
\begin{theorem}[Pointwise convergence]\label{lem: pointwise convergence total process}
For $t>0$,
\begin{align}
\FluidProcessTwo{t}\LimitP q(t)\text{ as $N\to\infty$},
\end{align}
with $q(t)$ given in Equation \eqref{eq: fluid limit function}.
\end{theorem}
As Equation \eqref{eq: maxqueue length expression} shows, we can write the scaled maximum queue length as a maximum of two random variables, namely, one pertaining to a system starting empty and one pertaining to a system starting non-empty. We prove the pointwise convergence of the first part of this maximum in Lemma \ref{lem: pointwise convergence}. In Lemma \ref{lem: pointwise convergence startposition} we prove the pointwise convergence of the second part. In order to do so, we need some extra results, which are stated in Lemmas \ref{lem: approximation starting point convergence}, \ref{lem: sup convergence tn3}, \ref{lem: upper bound pointwise convergence t large}, and \ref{lem: inf convergence tn3}.
\begin{lemma}\label{lem: pointwise convergence}
For $t>0$ and $\MaxQueueLength{0}=0$,
\begin{align*}
\FluidProcessTwo{t}\LimitP \bigg(\sqrt{2\alpha t}-\beta t\bigg)\mathbbm{1}\bigg(t<\frac{\alpha}{2\beta^2}\bigg)+\frac{\alpha}{2\beta}\mathbbm{1}\bigg(t\geq\frac{\alpha}{2\beta^2}\bigg)\text{ as $N\to\infty$}.
\end{align*}
\end{lemma}
To prove convergence of sequences of real-valued random variables to a constant it suffices to show convergence in distribution. Therefore, we use Lemmas \ref{lem: sup convergence tn3}, \ref{lem: upper bound pointwise convergence t large} and \ref{lem: inf convergence tn3}  below to prove that the upper and lower bound of the cumulative distribution function converge to the same function. 
\begin{lemma}\label{lem: sup convergence tn3}

For $\delta>0$, $t<\frac{\alpha}{(2\beta^2)}$ and $\MaxQueueLength{0}=0$,
\begin{align}\label{eq: limsup converg prob tn3}
\limsup_{N\to\infty}\mathbb{P}\left(\FluidProcessTwo{t}> \sqrt{2\alpha t}-\beta t+\delta\right)=0.
\end{align}
\end{lemma}
\proof{Proof}
Let $\delta>0$ be given. Let us assume that $t<\frac{\alpha}{
(2\beta^2)}$.  We then have that
\begin{align*}
&\probability*{\FluidProcessTwo{t}>\sqrt{2\alpha t}-\beta t+\delta}\nonumber\\
=&\probability*{\max_{i\leq N}\sup_{0\leq s\leq t}\left(\frac{\TildeArrival{sN^3\log N}+\TildeService{i}{sN^3\log N}}{\log N}\right)-\sqrt{2\alpha t}+\beta t>\delta}.
\end{align*}
For $t<\frac{\alpha}{(2\beta^2)}$, $\sqrt{2\alpha t}-\beta t$ is an increasing function. Therefore,
\begin{align*}
&\probability*{\max_{i\leq N}\sup_{0\leq s\leq t}\left(\frac{\TildeArrival{sN^3\log N}+\TildeService{i}{sN^3\log N}}{\log N}\right)-\sqrt{2\alpha t}+\beta t>\delta}\nonumber\\
 \leq &\probability*{\max_{i\leq N}\sup_{0\leq s\leq t}\left(\frac{\TildeArrival{sN^3\log N}+\TildeService{i}{sN^3\log N}}{\log N}-\sqrt{2\alpha s}+\beta s\right)>\delta}\nonumber\\
 =& \probability*{\sup_{0\leq s\leq t}\left(\max_{i\leq N}\frac{\TildeArrival{sN^3\log N}+\TildeService{i}{sN^3\log N}}{\log N}-\sqrt{2\alpha s}+\beta s\right)>\delta}.
\end{align*}
Observe that
\begin{align*}
&\probability*{\sup_{0\leq s\leq t}\left(\max_{i\leq N}\frac{\TildeArrival{sN^3\log N}+\TildeService{i}{sN^3\log N}}{\log N}-\sqrt{2\alpha s}+\beta s\right)>\delta}\nonumber\\
\leq&\probability*{\sup_{0\leq s\leq t}\left|\max_{i\leq N}\frac{\TildeArrival{sN^3\log N}+\TildeService{i}{sN^3\log N}}{\log N}-\sqrt{2\alpha s}+\beta s\right|>\delta}\nonumber\\
\leq&\probability*{\sup_{0\leq s\leq t}\left|\frac{\TildeArrival{sN^3\log N}}{\log N}+\beta s\right|>\frac{\delta}{2}}+\probability*{\sup_{0\leq s\leq t}\left|\frac{\max_{i\leq N}\TildeService{i}{sN^3\log N}}{\log N}-\sqrt{2\alpha s}\right|>\frac{\delta}{2}}.
\end{align*}
Moreover, $\frac{\TildeArrival{n}}{\log N}+\frac{\beta n}{(N^3\log N)}$ is a martingale with mean 0. 
Therefore, by Doob's maximal submartingale inequality
\begin{align}\label{eq: convergence arrival process}
&\probability*{\sup_{0\leq s\leq t}\left|\frac{\TildeArrival{sN^3\log N}}{\log N}+\beta s\right|>\frac{\delta}{2}}\nonumber\\
\leq &\probability*{\sup_{0\leq s\leq t}\left|\frac{\TildeArrival{sN^3\log N}}{\log N}+\beta \frac{\lfloor{sN^3\log N\rfloor}}{N^3\log N}\right|+\sup_{0\leq s\leq t}\left|\beta \frac{\lfloor{sN^3\log N\rfloor}}{N^3\log N}-\beta s\right|>\frac{\delta}{2}}\nonumber\\
\leq& \probability*{\sup_{0\leq s\leq t}\left|\frac{\TildeArrival{sN^3\log N}}{\log N}+\beta \frac{\lfloor{sN^3\log N\rfloor}}{N^3\log N}\right|>\frac{\delta}{4}}+\probability*{\sup_{0\leq s\leq t}\left|\beta \frac{\lfloor{sN^3\log N\rfloor}}{N^3\log N}-\beta s\right|>\frac{\delta}{4}}\nonumber\\
\leq&\frac{16}{\delta^2}\text{Var}\left(\frac{\TildeArrival{tN^3\log N}}{\log N}\right)+o_N(1)\nonumber\\
=&\frac{16}{\delta^2}\left(1-\frac{\alpha}{N}-\frac{\beta}{N^2}\right)\left(\frac{\alpha}{N}+\frac{\beta}{N^2}\right)\frac{\lfloor{tN^3\log N\rfloor}}{N^2(\log N)^2}+o_N(1)\LimitN 0.
\end{align}
Furthermore, in order to have
\begin{align}\label{eq: convergence service process}
\probability*{\sup_{0\leq s\leq t}\left|\frac{\max_{i\leq N}\TildeService{i}{sN^3\log N}}{\log N}-\sqrt{2\alpha s}\right|>\frac{\delta}{2}}\LimitN 0,
\end{align}
we need to have that $\left(\frac{\max_{i\leq N}\TildeService{i}{sN^3\log N}}{\log N},s\in[0,t]\right)$ converges to $\left(\sqrt{2\alpha s},s\in[0,t]\right)$ u.o.c. Thus
\begin{align}\label{eq: pw convergence S}
\lim_{N\to\infty}\probability*{\left|\frac{\max_{i\leq N}\TildeService{i}{sN^3\log N}}{\log N}-\sqrt{2\alpha s}\right|>\epsilon}=0,
\end{align}
and for all $r\in[0,t]$,
\begin{align}\label{eq: tightness convergence S}
\lim_{\eta\downarrow 0}\limsup_{N\to\infty}\frac{1}{\eta}\probability*{\sup_{r\leq s\leq r+\eta}\left|\frac{\max_{i\leq N}\TildeService{i}{sN^3\log N}}{\log N}-\frac{\max_{i\leq N}\TildeService{i}{rN^3\log N}}{\log N}\right|>\epsilon}=0.
\end{align} 
To prove the limit in \eqref{eq: pw convergence S}, we use the result of Lemma \ref{lem: distance cdf normal} and observe that for all $\delta>0$,
\begin{align*}
&\probability*{\frac{\max_{i\leq N}\TildeService{i}{sN^3\log N}}{\log N}>\sqrt{2\alpha s}+\delta}\nonumber\\
=&1-\probability*{\frac{\TildeService{i}{sN^3\log N}}{\log N}<\sqrt{2\alpha s}+\delta}^N\nonumber\\
=&1-\probability*{\MaximumProcess{i}{s}<\frac{\sqrt{2\alpha s}+\delta}{\sqrt{\alpha s(1-\frac{\alpha}{N})}}\sqrt{\log N}\frac{\sqrt{sN^3\log N}}{\sqrt{\lfloor{sN^3\log N\rfloor}}}}^N\nonumber\\
\leq& 1-\left(\Phi\left(\frac{\sqrt{2\alpha s}+\delta}{\sqrt{\alpha s(1-\frac{\alpha}{N})}}\sqrt{\log N}\frac{\sqrt{sN^3\log N}}{\sqrt{\lfloor{sN^3\log N\rfloor}}}\right)-\frac{c_s}{N\sqrt{\log N}}\right)^N\nonumber\\
\leq&1-\Phi\left(\frac{\sqrt{2\alpha s}+\delta}{\sqrt{\alpha s(1-\frac{\alpha}{N})}}\sqrt{\log N}\frac{\sqrt{sN^3\log N}}{\sqrt{\lfloor{sN^3\log N\rfloor}}}\right)^N+\left(1+\frac{c_s}{N\sqrt{\log N}}\right)^N-1\nonumber\\
\LimitN&0.
\end{align*}
The proof that
\begin{align*}
\probability*{\frac{\max_{i\leq N}\TildeService{i}{sN^3\log N}}{\log N}<\sqrt{2\alpha s}-\delta}\LimitN 0,
\end{align*}
goes analogously. 
To prove the quantity in \eqref{eq: tightness convergence S}, we observe that due to the facts that $\TildeService{i}{n}$ is a random walk that satisfies the duality principle, $\max_{i\leq N}x_i-\max_{i\leq N}y_i\leq \max_{i\leq N}(x_i-y_i)$, and $\probability*{\lvert X\rvert>\epsilon}\leq \probability*{X>\epsilon}+\probability*{-X>\epsilon}$, we have the upper bound
\begin{align*}
&\frac{1}{\eta}\probability*{\sup_{r\leq s\leq r+\eta}\left|\frac{\max_{i\leq N}\TildeService{i}{sN^3\log N}}{\log N}-\frac{\max_{i\leq N}\TildeService{i}{rN^3\log N}}{\log N}\right|>\epsilon}\nonumber\\
\leq&\frac{1}{\eta}\probability*{\sup_{0\leq s\leq \eta}\max_{i\leq N}\frac{\TildeService{i}{sN^3\log N}}{\log N}>\epsilon}+\frac{1}{\eta}\probability*{\sup_{0\leq s\leq \eta}\max_{i\leq N}\frac{-\TildeService{i}{sN^3\log N}}{\log N}>\epsilon}+o_N(1).
\end{align*}
The $o_N(1)$ term appears since $\lfloor (r+\eta)N^3\log N\rfloor-\lfloor r N^3\log N\rfloor\in\{\lfloor \eta N^3\log N\rfloor,\lfloor \eta N^3\log N\rfloor+1\}$. Now, we have that $\pm\TildeService{i}{n}$ is a martingale with mean 0. The maximum of independent martingales is a submartingale; therefore, 
$\left(\max\left(0,\max_{i\leq N}\frac{\pm\TildeService{i}{{\eta N^3\log N}}}{\log N}\right)\right)^{\frac{5}{2}}$ is a non-negative submartingale. Hence, by use Doob's maximal submartingale inequality we can conclude that
\begin{align*}
&\frac{1}{\eta}\probability*{\sup_{0\leq s\leq \eta}\max_{i\leq N}\frac{\TildeService{i}{sN^3\log N}}{\log N}>\epsilon}+\frac{1}{\eta}\probability*{\sup_{0\leq s\leq \eta}\max_{i\leq N}\frac{-\TildeService{i}{sN^3\log N}}{\log N}>\epsilon}\nonumber\\
\leq &\frac{1}{\eta\epsilon^{\frac{5}{2}}}\expect*{\max\left(\max_{i\leq N}\frac{\TildeService{i}{\eta N^3\log N}}{\log N},0\right)^{\frac{5}{2}}}+\frac{1}{\eta\epsilon^{\frac{5}{2}}}\expect*{\max\left(\max_{i\leq N}\frac{-\TildeService{i}{\eta N^3\log N}}{\log N},0\right)^{\frac{5}{2}}}.
\end{align*}
By taking the $\limsup_{N\to\infty}$ in this expression and applying Lemma \ref{lem: convergence moments 1 2 4}, we see that this is upper bounded by $ 2\frac{\eta^{\frac{1}{4}}(2\alpha)^{\frac{5}{4}}}{\epsilon^{\frac{5}{2}}}$.
This can be made as small as possible when $\eta$ is chosen small enough. We also know that $\frac{\max_{i\leq N}\TildeService{i}{0}}{\log N}=0$, and that the finite-dimensional distributions of $\left(\frac{\max_{i\leq N}\TildeService{i}{sN^3\log N}}{\log N},s\in[0,t]\right)$ converge to the finite-dimensional distributions of $\left(\sqrt{2\alpha s},s\in[0,t]\right)$, which follows from Theorem \ref{lem: finite-dimensional distributions}. The lemma follows.
\hfill\rlap{\hspace*{-2em}\Halmos}
\endproof
Having examined $t\in[0,\frac{\alpha}{(2\beta^2)})$, we now turn to $t\in[\frac{\alpha}{(2\beta^2)},\infty]$.
\begin{lemma}\label{lem: upper bound pointwise convergence t large}
For $\delta>0$, $\frac{\alpha}{(2\beta^2)}\leq t\leq \infty$ and $\MaxQueueLength{0}=0$,
\begin{align*}
\limsup_{N\to \infty}\probability*{\FluidProcessTwo{t}>\frac{\alpha}{2\beta}+\delta}=0.
\end{align*}
\end{lemma}
\proof{Proof}
We write
\begin{align*}
\ArrivalUpper{n}=\sum_{j=1}^n \ArrivalBernoulliUpper{j}
\end{align*}
with
\begin{align*}
\ArrivalBernoulliUpper{j}=\left\{\begin{matrix}
\text{ }\text{ }\text{ }\text{ }\text{ }\text{ }\text{ }\text{ }\frac{\alpha}{N}+\frac{\beta}{N^{2}}-\frac{m}{N^2}&\text{ w.p. }&  1-\frac{\alpha}{N}-\frac{\beta}{N^2},&\\ 
-1+\frac{\alpha}{N}+\frac{\beta}{N^{2}}-\frac{m}{N^2}&\text{ w.p. }& \text{ }\text{ }\text{ }\text{ }\text{ }\frac{\alpha}{N}+\frac{\beta}{N^2},&
\end{matrix}\right.
\end{align*}
with $0<m<\beta$.
Furthermore, we write
\begin{align*}
\ServiceUpper{i}{n}=\sum_{j=1}^n \ServiceBernoulliUpper{i}{j},
\end{align*}
with
\begin{align*}
\ServiceBernoulliUpper{i}{j}=\left\{\begin{matrix}\text{ }
-\frac{\alpha}{N}-\frac{\beta}{N^{2}}+\frac{m}{N^2}&\text{ w.p. }&  1-\frac{\alpha}{N},\\ 
1-\frac{\alpha}{N}-\frac{\beta}{N^{2}}+\frac{m}{N^2}&\text{ w.p. }& \text{ }\text{ }\text{ }\text{ }\text{ }\frac{\alpha}{N}.
\end{matrix}\right.
\end{align*}
Thus,
\begin{align*}
\Arrival{n}-\Service{i}{n}=\ArrivalUpper{n}+\ServiceUpper{i}{n},
\end{align*}
and
\begin{align*}
\sup_{0\leq k\leq n}\left(\Arrival{k}-\Service{i}{k}\right)\leq \sup_{0\leq k\leq n}\ArrivalUpper{k}+\sup_{0\leq k\leq n}\ServiceUpper{i}{k}.
\end{align*}
We obtain by using Doob's maximal submartingale inequality that
\begin{align*}
\probability*{\sup_{0\leq k\leq n} \ArrivalUpper{k}\geq x}\leq \mathbb{E}\left[e^{\ThetaArrivalUpper \ArrivalBernoulliUpper{j}}\right]e^{-\ThetaArrivalUpper x}=e^{-\ThetaArrivalUpper x},
\end{align*}
with $\ThetaArrivalUpper $ the solution to the equation
\begin{align*}
\mathbb{E}\left[e^{\ThetaArrivalUpper \ArrivalBernoulliUpper{j}}\right]&=\left(\frac{\alpha}{N}+\frac{\beta}{N^2}\right)\text{exp}\left\{\ThetaArrivalUpper \left(-1+\frac{\alpha}{N}+\frac{\beta}{N^2}-\frac{m}{N^2}\right)\right\}\nonumber\\
&+\left(1-\frac{\alpha}{N}-\frac{\beta}{N^2}\right)\text{exp}\left\{\ThetaArrivalUpper \left(\frac{\alpha}{N}+\frac{\beta}{N^2}-\frac{m}{N^2}\right)\right\}=1.
\end{align*}
When we consider the second-order Taylor approximation of this expression with $1/N$ around 0, we obtain
\begin{align*}
\ThetaArrivalUpper =\frac{2m N^2 }{-\alpha ^2 N^2+\alpha  N^3-2 \alpha  \beta  N-\beta ^2+m^2+\beta  N^2}+O\left(\frac{1}{N^2}\right).
\end{align*}
Consequently, we have for $N$ large $\ThetaArrivalUpper \approx \frac{2m}{(\alpha N)}$. By the monotone convergence theorem, we know that
\begin{align*}
\probability*{\sup_{k\geq 0}\ArrivalUpper{k}\geq x}\leq e^{-\ThetaArrivalUpper x}\approx e^{-\frac{2m}{(\alpha N)}x}.
\end{align*}
In conclusion,
\begin{align*}
\frac{\sup_{k\geq 0}\ArrivalUpper{k}}{N\log N}\LimitP 0\text{ as $N\to\infty$}.
\end{align*}
Similarly, by using Doob's maximal submartingale inequality, we obtain that
\begin{align*}
\probability*{\sup_{n\geq 0} \ServiceUpper{i}{n}\geq x}\leq e^{-\ThetaServiceUpper{i}x},
\end{align*}
with $\ThetaServiceUpper{i}$ the solution to the equation
\begin{align*}
\mathbb{E}\left[e^{\ThetaServiceUpper{i}\ServiceBernoulliUpper{i}{j}}\right]&=\frac{\alpha}{N}\text{exp}\left\{\ThetaServiceUpper{i}\left(1-\frac{\alpha}{N}-\frac{\beta}{N^2}+\frac{m}{N^2}\right)\right\}\nonumber\\
&+\left(1-\frac{\alpha}{N}\right)\text{exp}\left\{\ThetaServiceUpper{i}\left(-\frac{\alpha}{N}-\frac{\beta}{N^2}+\frac{m}{N^2}\right)\right\}=1.
\end{align*}
The second-order Taylor approximation of $\mathbb{E}\left[e^{\ThetaServiceUpper{i}\ServiceBernoulliUpper{i}{j}}\right]$ with $1/N$ around 0 gives
\begin{align*}
\ThetaServiceUpper{i}=\frac{2 N^2 \left(\beta -m \right)}{-\alpha ^2 N^2+\alpha  N^3+\left(\beta -m \right)^2}+O\left(\frac{1}{N^2}\right).
\end{align*}
Thus, for $N$ large, $\ThetaServiceUpper{i}\approx\frac{2(\beta-m)}{(\alpha N)}$. Concluding, $\sup_{n\geq 0} \ServiceUpper{i}{n}$ is stochastically dominated by an exponentially distributed random variable $E_i^{(u,N)}$ with mean $\frac{\alpha N}{(2(\beta-m))}$. Because $\sup_{n\geq 0} \ServiceUpper{i}{n}\perp \sup_{n\geq 0} \ServiceUpper{j}{n}$ for $i\neq j$, we can conclude that also $E_i^{(u,N)}\perp E_j^{(u,N)}$ for $i\neq j$. Therefore,
\begin{align*}
\probability*{\frac{\max_{i\leq N}E_i^{(u,N)}}{N}\leq \frac{\alpha}{2(\beta-m)}\left(x+\log N\right)}\LimitN e^{-e^{-x}},
\end{align*}
and
\begin{align*}
\frac{\max_{i\leq N}E_i^{(u,N)}}{N\log N}\LimitP\frac{\alpha}{2(\beta-m)}\text{ as $N\to\infty$}.
\end{align*}
Because,
\begin{align*}
\FluidProcessTwo{t}\leq_{st.} \FluidProcessTwo{\infty}\leq \frac{\sup_{k\geq 0}\ArrivalUpper{k}}{N\log N}+\frac{\max_{i\leq N}\sup_{k\geq 0}\Service{i}{k}}{N\log N},
\end{align*}
the lemma follows.
 \hfill\rlap{\hspace*{-2em}\Halmos}
\endproof
\begin{lemma}\label{lem: inf convergence tn3}

For $\delta>0$ and $\MaxQueueLength{0}=0$,
\begin{align}\label{eq: liminf converg prob tn3}
\liminf_{N\to\infty}\mathbb{P}\left(\FluidProcessTwo{t}\geq \bigg(\sqrt{2\alpha t}-\beta t\bigg)\mathbbm{1}\bigg(t<\frac{\alpha}{2\beta^2}\bigg)+\frac{\alpha}{2\beta}\mathbbm{1}\bigg(t\geq\frac{\alpha}{2\beta^2}\bigg)-\delta\right)=1.
\end{align}
\end{lemma}
\proof{Proof}

Let us first assume that $t\leq\frac{\alpha}{(2\beta^2)}$. We have the lower bound
\begin{align*}
\FluidProcessTwo{t}\geq_{st.} \max_{i\leq N}\frac{\Arrival{tN^3\log N}-\Service{i}{tN^3\log N}}{N\log N}.
\end{align*}
By Equations \eqref{eq: convergence arrival process} and \eqref{eq: convergence service process}, we know that
\begin{align*}
\max_{i\leq N}\frac{\Arrival{tN^3\log N}-\Service{i}{tN^3\log N}}{N\log N}\LimitP \sqrt{2\alpha t}-\beta t\text{ as $N\to\infty$}.
\end{align*}
Let us now assume that $t>\frac{\alpha}{(2\beta^2)}$. We have that
\begin{align*}
\FluidProcessTwo{t}\geq_{st.} \max_{i\leq N}\frac{\Arrival{\oldfrac{\alpha}{2\beta^2}N^3\log N}-\Service{i}{\oldfrac{\alpha}{2\beta^2}N^3\log N}}{N\log N}\LimitP \frac{\alpha}{2\beta},
\end{align*} 
as $N\to\infty$, by again using Lemma \ref{lem: sup convergence tn3}. This proves the lemma. \hfill\rlap{\hspace*{-2em}\Halmos}
\endproof
\text{ }
\proof{Proof of Lemma \ref{lem: pointwise convergence}}
By combining the results of Lemmas \ref{lem: sup convergence tn3}, \ref{lem: upper bound pointwise convergence t large} and \ref{lem: inf convergence tn3}, Lemma \ref{lem: pointwise convergence} follows.

\hfill\rlap{\hspace*{-2em}\Halmos}
\endproof 
In Lemma \ref{lem: pointwise convergence startposition}, we connect the convergence of
\begin{align*}
\max_{i\leq N}\StartPositionGeneral{i}{t}
\end{align*}
to the convergence of
\begin{align*}
\max_{i\leq N}\left(\frac{\sqrt{\alpha t}X_i}{\sqrt{\log N}}+\frac{\QueueZero{i}}{N\log N}\right).
\end{align*}
\begin{lemma}[Convergence starting position]\label{lem: pointwise convergence startposition}
Assume that for $X_i$ i.i.d.\ standard normally distributed, 
\begin{align}\label{eq: convergence sum X plus startposition}
\max_{i\leq N}\left(\frac{\sqrt{\alpha t}X_i}{\sqrt{\log N}}+\frac{\QueueZero{i}}{N\log N}\right)\LimitP g(t,q(0))\text{ as $N\to\infty$},
\end{align}
for a certain function $g$. Then
\begin{align*}
\max_{i\leq N}\StartPositionGeneral{i}{t}
\LimitP g(t,q(0))-\beta t\text{ as $N\to\infty$}.
\end{align*}
\end{lemma}
\proof{Proof}
We have
\begin{align}
&\max_{i\leq N}\StartPositionGeneral{i}{t}\\
=&\frac{\Arrival{tN^3\log N}-\left(1-\frac{\alpha}{N}\right)tN^3\log N}{N\log N}+\max_{i\leq N}\frac{\left(1-\frac{\alpha}{N}\right)tN^3\log N-\Service{i}{tN^3\log N}+\QueueZero{i}}{N\log N}.\label{subeq: split up startposition}
\end{align}
We already proved in Equation \eqref{eq: convergence arrival process} that the first term in \eqref{subeq: split up startposition} converges to $-\beta t$. Furthermore, we can rewrite the second term as 
\begin{align*}
\max_{i\leq N}\left(\frac{\TildeService{i}{tN^3\log N}}{\log N}+\frac{\QueueZero{i}}{N\log N}
+O_N\left(\frac{1}{N\log N}\right)\right).
\end{align*}
We can easily deduce from Lemma \ref{lem: distance cdf normal} that
\begin{align*}
\left|\probability*{\frac{\TildeService{i}{tN^3\log N}}{\log N}<y}-\probability*{\frac{\sqrt{\alpha t(1-\frac{\alpha}{N})}}{\sqrt{\log N}}\frac{\sqrt{\lfloor{tN^3\log N\rfloor}}}{\sqrt{tN^3\log N}}X_i<y}\right|
\leq\frac{c_t}{N\sqrt{\log N}},
\end{align*}
with $X_i\sim\mathcal{N}(0,1)$, and $c_t$ given in Lemma \ref{lem: distance cdf normal}. Then, it is easy to see that
\begin{align}\label{eq: bound probabilities starting position}
&\left|\probability*{\frac{\TildeService{i}{tN^3\log N}}{\log N}+\frac{\QueueZero{i}}{N\log N}<y}-\probability*{\frac{\sqrt{\alpha t(1-\frac{\alpha}{N})}}{\sqrt{\log N}}\frac{\sqrt{\lfloor{tN^3\log N\rfloor}}}{\sqrt{tN^3\log N}}X_i+\frac{\QueueZero{i}}{N\log N}<y}\right|\nonumber\\
\leq&\frac{c_t}{N\sqrt{\log N}}.
\end{align}
Now, because of the facts that we assume the convergence result in \eqref{eq: convergence sum X plus startposition}, and
\begin{align*}
\frac{\sqrt{\alpha t(1-\frac{\alpha}{N})}}{\sqrt{\log N}}\frac{\sqrt{\lfloor{tN^3\log N\rfloor}}}{\sqrt{tN^3\log N}}X_i=\frac{\sqrt{\alpha t}X_i}{\sqrt{\log N}}+o_N\left(\frac{1}{\sqrt{\log N}}\right)X_i,
\end{align*}
it is easy to see that
\begin{align*}
\max_{i\leq N}\left(\frac{\sqrt{\alpha t(1-\frac{\alpha}{N})}}{\sqrt{\log N}}\frac{\sqrt{\lfloor{tN^3\log N\rfloor}}}{\sqrt{tN^3\log N}}X_i+\frac{\QueueZero{i}}{N\log N}\right)\LimitP g(t,q(0))\text{ as $N\to\infty$}.
\end{align*}
Let $\epsilon>0$, then because of the bound given in \eqref{eq: bound probabilities starting position}, and the convergence result in \eqref{eq: convergence sum X plus startposition},
\begin{align*}
&\probability*{\max_{i\leq N}\left(\frac{\TildeService{i}{tN^3\log N}}{\log N}+\frac{\QueueZero{i}}{N\log N}\right)<g(t,q(0))-\epsilon}\nonumber\\
=&\probability*{\frac{\TildeService{i}{tN^3\log N}}{\log N}+\frac{\QueueZero{i}}{N\log N}<g(t,q(0))-\epsilon}^N\nonumber\\
\leq& \probability*{\frac{\sqrt{\alpha t(1-\frac{\alpha}{N})}}{\sqrt{\log N}}\frac{\sqrt{\lfloor{tN^3\log N\rfloor}}}{\sqrt{tN^3\log N}}X_i+\frac{\QueueZero{i}}{N\log N}<g(t,q(0))-\epsilon}^N+\left(\frac{c_t}{N\sqrt{\log N}}+1\right)^N-1\nonumber\\
\LimitN& 0.
\end{align*}
The proof that 
\begin{align*}
&\probability*{\max_{i\leq N}\left(\frac{\TildeService{i}{tN^3\log N}}{\log N}+\frac{\QueueZero{i}}{N\log N}\right)>g(t,q(0))+\epsilon}\LimitN 0,
\end{align*}
goes analogously. Hence, the lemma follows.
\hfill\rlap{\hspace*{-2em}\Halmos}
\endproof
\text{ }
\proof{Proof of Theorem \ref{lem: pointwise convergence total process}}
In Lemmas \ref{lem: pointwise convergence} and \ref{lem: pointwise convergence startposition} we have proven that both parts in the maximum in \eqref{eq: maxqueue length expression} converge to a limit. The lemma follows.
\hfill\rlap{\hspace*{-2em}\Halmos}
\endproof
We can easily extend this result to finite-dimensional distributions.
\begin{theorem}[The finite-dimensional distributions converge]\label{lem: finite-dimensional distributions}
If
\begin{align*}
X^{(N)}(t)\LimitP f(t)
\end{align*}
for all $t>0$, then for $(t_1,t_2,\ldots,t_k)$ 
\begin{align*}
\left(X^{(N)}(t_1),X^{(N)}(t_2),\ldots,X^{(N)}(t_k)\right)\LimitP\left(f(t_1),f(t_2),\ldots,f(t_k)\right)\text{ as $N\to\infty$}.
\end{align*}
\end{theorem}
\proof{Proof}
\begin{align*}
&\probability*{\left\lVert\left(X^{(N)}(t_1),X^{(N)}(t_2),\ldots,X^{(N)}(t_k)\right)-\left(f(t_1),f(t_2),\ldots,f(t_k)\right)\right\rVert>\epsilon}\nonumber\\
\leq &\probability*{\left| X^{(N)}(t_1)-f(t_1)\right|+\cdots+\left| X^{(N)}(t_k)-f(t_k)\right|>\epsilon}\nonumber\\
\leq &\probability*{\left|X^{(N)}(t_1)-f(t_1)\right|>\frac{\epsilon}{k}}+\cdots+\probability*{\left|X^{(N)}(t_k)-f(t_k)\right|>\frac{\epsilon}{k}}\LimitN 0,
\end{align*}
with $\left\lVert \cdot\right\rVert$ the Euclidean distance in $\mathbb{R}^k$.
\hfill\rlap{\hspace*{-2em}\Halmos}
\endproof

\subsection[Tightness]{Tightness.}\label{subsubsec: tightness} 
It is known that when a sequence of random processes is tight and its finite-dimensional distributions converge, then this sequence converges u.o.c.; cf.\ \cite[Thm.~7.1, p.~80]{billingsley1968convergence}. From \cite[Thm.~7.3, p.~82]{billingsley1968convergence}, we know that a process $(X^{(N)}(t),t\in [0,T])$ is tight when for all positive $\eta$ there exists an $a$ and an integer $N_0$ such that for all $N\geq N_0$
\begin{align}\label{eq: tightness condition 1}
\probability*{\left|X^{(N)}(0)\right|>a}\leq \eta,
\end{align}
and for all $\epsilon>0$ and $\eta>0$, there exists a $0<\delta<1$ and an integer $N_0$ such that for all $N\geq N_0$
\begin{align}\label{eq: tightness condition 2}
\frac{1}{\delta}\probability*{\sup_{t\leq s\leq t+\delta}\left|X^{(N)}(s)-X^{(N)}(t)\right|>\epsilon}\leq \eta.
\end{align}
The conditions given in Equations \eqref{eq: tightness condition 1} and \eqref{eq: tightness condition 2} hold for stochastic processes in the space of continuous functions. The process $\big(\FluidProcessTwo{t},t\in[0,T]\big)$ does not lie in this space, because $\MaxQueueLength{n}=\MaxQueueLength{\lfloor n\rfloor}$. However, since $q(t)$ is a continuous function, the conditions in \eqref{eq: tightness condition 1} and \eqref{eq: tightness condition 2} do also apply on $\big(\FluidProcessTwo{t},t\in[0,T]\big)$; cf.\ \cite[Cor.~13.4, p.~142]{billingsley1968convergence}.

In order to prove tightness for the process given in Theorem \ref{thm: fluid limit}, we need to prove tightness of the maximum of two processes, as Equation \eqref{eq: maxqueue length expression} shows. In Lemma \ref{lem: tightness max two processes}, we show that it suffices to prove tightness of the two processes separately. Then, in Lemmas \ref{lem: tight process} and \ref{lem: tightness startposition}, we prove the tightness of the two parts.
\begin{lemma}\label{lem: tightness max two processes}
Assume that $(X^{(N)}(s),s\in[0,t])$ and $(Y^{(N)}(s),s\in[0,t])$ converge to functions $(k(s),s\in[0,t])$ and $(l(s),s\in[0,t])$ u.o.c., respectively, then $(\max(X^{(N)}(s),Y^{(N)}(s)),s\in[0,t])$ converges to $(\max(k(s),l(s)),s\in[0,t])$ u.o.c.
\end{lemma}
\proof{Proof} 
The lemma holds because of the fact that
\begin{align*}
&\probability*{\sup_{0\leq s\leq t}\left|\max(X^{(N)}(s),Y^{(N)}(s))-\max(k(s),l(s))\right|>\epsilon}\nonumber\\
\leq&\probability*{\sup_{0\leq s\leq t}(\max(X^{(N)}(s),Y^{(N)}(s))-\max(k(s),l(s)))>\epsilon}\nonumber\\
+&\probability*{\sup_{0\leq s\leq t}(\max(k(s),l(s))-\max(X^{(N)}(s),Y^{(N)}(s)))>\epsilon}\nonumber\\
\leq&\probability*{\sup_{0\leq s\leq t}\max(X^{(N)}(s)-k(s),Y^{(N)}(s)-l(s))>\epsilon}\nonumber\\
+&\probability*{\sup_{0\leq s\leq t}\max(k(s)-X^{(N)}(s),l(s)-Y^{(N)}(s))>\epsilon}\nonumber\\
\leq &2\probability*{\sup_{0\leq s\leq t}\left|X^{(N)}(s)-k(s)\right|>\epsilon}+2\probability*{\sup_{0\leq s\leq t}\left|Y^{(N)}(s)-l(s)\right|>\epsilon}\LimitN 0.
\end{align*}
\hfill\rlap{\hspace*{-2em}\Halmos}
\endproof
\begin{lemma}[Tightness of the first part]\label{lem: tight process}
For $\epsilon>0$, $\eta>0$, $T>0$ and $\MaxQueueLength{0}=0$, $\exists 0<\delta<1$ and an integer $N_0$ such that $\forall N\geq N_0$ and $t\in[0,T]$

\begin{align}
\frac{1}{\delta}\probability*{\sup_{t\leq s\leq t+\delta}\left|\FluidProcessTwo{s}-\FluidProcessTwo{t}\right|\geq\epsilon}\leq \eta.
\end{align}
\end{lemma}
\proof{Proof}
We take $t>0$. From Lemma \ref{lem: upper bound |Qs-Qt|}, and the fact that $\TildeQueue{i}{}$ is a random walk that satisfies the duality principle, we know that for $N$ large enough, 
\begin{align}
&\frac{1}{\delta}\probability*{\sup_{t\leq s\leq t+\delta}\left|\FluidProcessTwo{s}-\FluidProcessTwo{t}\right|\geq\epsilon}\\
\leq&\frac{1}{\delta}\probability*{\sup_{0\leq s\leq \delta}\max_{i\leq N}\TildeQueue{i}{sN^3\log N}+2\sup_{0\leq s\leq \delta}\max_{i\leq N}-\TildeQueue{i}{sN^3\log N}\geq\epsilon}+o_N(1)\\
\leq&\frac{1}{\delta}\probability*{\sup_{0\leq s\leq \delta}\max_{i\leq N}\TildeQueue{i}{sN^3\log N}\geq\frac{\epsilon}{2}}+\frac{1}{\delta}\probability*{2\sup_{0\leq s\leq \delta}\max_{i\leq N}-\TildeQueue{i}{sN^3\log N}\geq\frac{\epsilon}{2}}+o_N(1)\label{subeq: inequality tight}.
\end{align}
Now we focus on the first term in \eqref{subeq: inequality tight}. The analysis of the second term goes analogously.
\begin{align}
&\frac{1}{\delta}\probability*{\sup_{0\leq s\leq \delta}\max_{i\leq N}\TildeQueue{i}{sN^3\log N}\geq\frac{\epsilon}{2}}\\
=&\frac{1}{\delta}\probability*{\sup_{0\leq s\leq \delta}\max_{i\leq N}\frac{\TildeArrival{{sN^3\log N}}+\TildeService{i}{{sN^3\log N}}}{\log N}\geq\frac{\epsilon}{2}}\label{subeq: tightness arrival +service}\\
\leq&\frac{1}{\delta}\probability*{\sup_{0\leq s\leq \delta}\frac{\TildeArrival{{sN^3\log N}}}{\log N}\geq\frac{\epsilon}{4}}+\frac{1}{\delta}\probability*{\sup_{0\leq s\leq \delta}\max_{i\leq N}\frac{\TildeService{i}{{sN^3\log N}}}{\log N}\geq\frac{\epsilon}{4}}\label{subeq: splitting arrival service tightness}.
\end{align}
In the proof of Lemma \ref{lem: sup convergence tn3}, we already showed that the second term in \eqref{subeq: splitting arrival service tightness} is small. With a similar proof as in Lemma \ref{lem: sup convergence tn3}, one can also prove that the first term is small. Concluding, $\big(\FluidProcessTwo{t},t\in[0,T]\big)$ is tight, when $\MaxQueueLength{0}=0$.\hfill\rlap{\hspace*{-2em}\Halmos}
\endproof

\begin{lemma}[Tightness of the second part]\label{lem: tightness startposition}
For $\epsilon>0$, $\eta>0$ and $T>0$, $\exists 0<\delta<1$ and an integer $N_0$ such that $\forall N\geq N_0$ and $t\in[0,T]$
\begin{align}\label{eq: billingsely time >0}
\frac{1}{\delta}\mathbb{P}\bigg(\sup_{t\leq s\leq t+\delta}\bigg|&\max_{i\leq N}\StartPositionGeneral{i}{s}\nonumber\\
-&\max_{i\leq N}\StartPositionGeneral{i}{t}\bigg|>\epsilon\bigg)<\eta.
\end{align}
Furthermore, for all $\eta$ there exists an $a>0$ such that 
\begin{align}\label{eq: billingsley time 0}
\probability*{\frac{\MaxQueueLength{0}}{N\log N}>a}<\eta.
\end{align}
\end{lemma}
\proof{Proof}
First of all, we observe that for a random variable $X$, $\probability*{\lvert X\rvert>\epsilon}\leq \probability*{X>\epsilon}+\probability*{-X>\epsilon}$. Thus, we can remove the absolute values in \eqref{eq: billingsely time >0} and examine both cases. Since both cases satisfy analogous proofs, we only write down the proof for the first case.
\begin{align*}
\frac{1}{\delta}\mathbb{P}\bigg(\sup_{t\leq s\leq t+\delta}\bigg(&\max_{i\leq N}\StartPositionGeneral{i}{s}\nonumber\\
-&\max_{i\leq N}\StartPositionGeneral{i}{t}\bigg)>\epsilon\bigg)\nonumber\\
\leq&\nonumber\\
 \frac{1}{\delta}\mathbb{P}\Bigg(\sup_{t\leq s\leq t+\delta}\Bigg(&\max_{i\leq N}\Bigg(\frac{\Arrival{sN^3\log N}-\Service{i}{sN^3\log N}}{N\log N}\nonumber\\
 &\text{ }\text{ }\text{ }\text{ }\text{ }\text{ }-\frac{\Arrival{tN^3\log N}-\Service{i}{tN^3\log N}}{N\log N}\Bigg)>\epsilon\Bigg)\nonumber\\
=&\frac{1}{\delta}\probability*{\sup_{0\leq s\leq \delta}\left(\max_{i\leq N}\frac{\Arrival{sN^3\log N}-\Service{i}{sN^3\log N}}{N\log N}\right)>\epsilon}+o_N(1).
\end{align*}
This is the same expression as Equation \eqref{subeq: tightness arrival +service}. In Lemma \ref{lem: tight process}, it is proven that this expression will be small.
At $t=0$, we should choose $a>0$ such that \eqref{eq: billingsley time 0} holds for $N\geq N_0$. This is the case because we know that $\frac{\MaxQueueLength{0}}{(N\log N)}\LimitP q(0)\text{ as $N\to\infty$}.$ The lemma follows.
\hfill\rlap{\hspace*{-2em}\Halmos}
\endproof
\begin{corollary}[Tightness of the process]\label{cor: tightness process}
The process $\big(\FluidProcessTwo{t},t\in[0,T]\big)$ is tight.
\end{corollary}
\proof{Proof}
The process $\big(\FluidProcessTwo{t},t\in[0,T]\big)$ can be written as a maximum of two processes. In Lemmas \ref{lem: tight process} and \ref{lem: tightness startposition}, it is proven that these processes are tight. Then from Lemma \ref{lem: tightness max two processes}, it follows that $\big(\FluidProcessTwo{t},t\in[0,T]\big)$ is tight.
\hfill\rlap{\hspace*{-2em}\Halmos}
\endproof
\text{ }
\proof{Proof of Theorem \ref{thm: fluid limit}}
In Theorem \ref{lem: pointwise convergence total process}, we proved that for fixed $t$, the stochastic process converges in probability to a constant, in Theorem \ref{lem: finite-dimensional distributions}, we proved that the finite-dimensional distributions converge and in Corollary \ref{cor: tightness process}, we showed that the process is tight. Thus the convergence holds u.o.c.
\hfill\rlap{\hspace*{-2em}\Halmos}
\endproof
\text{ }

We now prove that the scaled process in steady state converges to the constant $\frac{\alpha}{(2\beta)}$.
\proof{Proof of Proposition \ref{cor: steady state convergence}.}
Since we look at the system in steady state, we can assume w.l.o.g. that $\MaxQueueLength{0}=0$. Then, we have 
\begin{align*}
\FluidProcessTwo{\infty}\geq_{st.} \FluidProcessTwo{\frac{\alpha}{(2\beta^2)}},
\end{align*}
because $\MaxQueueLength{n}\overset{d}{=} \max_{i\leq N}\sup_{0\leq k\leq n}(\Arrival{k}-\Service{i}{k})$.
We know by Lemma \ref{lem: pointwise convergence} that 
\begin{align*}
\FluidProcessTwo{\frac{\alpha}{(2\beta^2)}}\LimitP \frac{\alpha}{2\beta}\text{ as $N\to\infty$}.
\end{align*}
Furthermore, we know by Lemma \ref{lem: upper bound pointwise convergence t large} that for all $\delta>0$,
\begin{align*}
\limsup_{N\to \infty}\probability*{\FluidProcessTwo{\infty}>\frac{\alpha}{2\beta}+\delta}=0.
\end{align*}
The proposition follows. 
\hfill\rlap{\hspace*{-2em}\Halmos}
\endproof

%
%
%
\pdfbookmark[0]{Appendices}{app}
\begin{APPENDICES}
\section[Taylor expansion of $\ThetaArrivalUpper$]{Taylor expansion of $\ThetaArrivalUpper$.}\label{app: taylor approx theta}
The parameter $\ThetaArrivalUpper $ is the strictly positive solution to the equation
\begin{align*}
\mathbb{E}\left[e^{\ThetaArrivalUpper \ArrivalBernoulliUpper{j}}\right]&=\left(\frac{\alpha}{N}+\frac{\beta}{N^2}\right)\text{exp}\left\{\ThetaArrivalUpper \left(-1+\frac{\alpha}{N}+\frac{\beta}{N^2}-\epsilon(N)\right)\right\}\nonumber\\
&+\left(1-\frac{\alpha}{N}-\frac{\beta}{N^2}\right)\text{exp}\left\{\ThetaArrivalUpper \left(\frac{\alpha}{N}+\frac{\beta}{N^2}-\epsilon(N)\right)\right\}=1,
\end{align*}
with $\epsilon(N)=\frac{m}{N^2}$. We found an approximation of $\ThetaArrivalUpper$, of $\frac{2m}{(\alpha N)}$. To investigate the behavior of $\ThetaArrivalUpper $ more carefully, we look at the function $\theta(x)$ such that
\begin{align*}
f(x,\theta(x))&=\left(\alpha x+\beta x^2\right)\text{exp}\left\{\theta(x)\left(-1+\alpha x+\beta x^2-mx^2\right)\right\}\nonumber\\
&+\left(1-\alpha x-\beta x^2\right)\text{exp}\left\{\theta(x)\left(\alpha x+\beta x^2-mx^2\right)\right\}=1.
\end{align*}
When we set $x_N=\frac{1}{N}$, we get $f(x_N,\theta(x_N))=\mathbb{E}\left[e^{\ThetaArrivalUpper \ArrivalBernoulliUpper{j}}\right]$. We are interested in the case that $N$ is large, therefore we have to investigate $f$ for $x$ around 0. Since $f(x,\theta(x))=1$, we know that $f^{(n)}(0,\theta(0))=0$ for all $n\geq 1$. When we solve these equations for $\theta$ iteratively, we can find $\theta^{(i)}(0)$ for all $i\geq 0$ and we get a Taylor expansion of $\theta(x)$ around 0.
Since $f(x,\theta(x))=1$, we know that
\begin{align*}
\left.\frac{d}{dx}f(x,\theta(x))\right|_{x=0}=-\alpha +\alpha  e^{-\theta(0)}+\alpha  \theta(0)=0.
\end{align*}
Hence, $\theta(0)=0$. When we look at the second and the third derivative of $f(x,\theta(x))$ around 0, while using that $\theta(0)=0$, we see 
\begin{align*}
\left.\frac{d^2}{dx^2}f(x,\theta(x))\right|_{x=0}=0,
\end{align*}
and
\begin{align*}
\left.\frac{d^3}{dx^3}f(x,\theta(x))\right|_{x=0}=3 \theta '(0) \left(\alpha  \theta '(0)-2 m\right).
\end{align*}
Because we know that $f(x,\theta(x))=1$, we solve $3 \theta '(0) \left(\alpha  \theta '(0)-2 m\right)=0.$
This gives $\theta'(0)=0$ or $\theta'(0)=\frac{2m}{\alpha}$. $\theta'(0)=0$ indicates the situation that $\theta\equiv 0$. If we now use the information that $\theta'(0)=\frac{2m}{\alpha}$ and look at the fourth derivative of $f$ we see that
\begin{align*}
\left.\frac{d^4}{dx^4}f(x,\theta(x))\right|_{x=0}=4 m \left(3 \theta ''(0)-\frac{4 m \left(3 \alpha ^2-3 \beta +2 m\right)}{\alpha ^2}\right)=0.
\end{align*}
This gives that $\theta''(0)=\frac{4 m \left(3 \alpha ^2-3 \beta +2 m\right)}{3 \alpha ^2}.$
In general, we can compute each derivative of $\theta(0)$ iteratively. This gives
\begin{align*}
\theta(x)=\frac{2m}{\alpha}x+\frac{4 m \left(3 \alpha ^2-3 \beta +2 m\right)}{3 \alpha ^2}\frac{x^2}{2}+O(x^3).
\end{align*}
Since the function $f(x,\theta)-1$ is analytic we know by the implicit function theorem that the solution $\theta(x)$ is also analytic. So for $x=\frac{1}{N}$ and $N$ is large enough we know that $\ThetaArrivalUpper =\frac{2m}{(\alpha N)}+O\left(\frac{1}{N^2}\right).$
\section[Extreme values of sums of random variables]{Extreme values of sums of random variables.}\label{app: extreme values}
In this section, we prove a convergence result of the maximum of $N$ sums of $n$ random variables. In order to do so, we use and extend results from \cite{davis1988almost} and \cite{fisher1969limiting}. 
\begin{lemma}\label{lem: extreme value convergence}
Consider sequences of continuous random variables $(Y^{(1)}_i,i\geq 1)$, $(Y^{(2)}_i,i\geq 1)$, $\ldots,$ $(Y^{(k)}_i,i\geq 1)$, where all random variables in the sequence $(Y^{(j)}_i,i\geq 1)$ are identically and independently distributed and have infinite right endpoints. Furthermore, $Y^{(j)}_i$ and $Y^{(l)}_m$ are independent for all $j,l\in\{1,2,\ldots,k\}$ and $i,m\geq 1$, and $Y_i^{(j)}$ satisfies Assumption \ref{ass: 5} with function $h^{(j)}(u^{(j)})$. Finally, we have sequences $(a_N^{(j)},N\geq 1)$ such that $\probability*{Y^{(j)}_i\geq a_N^{(j)}}=\frac{1}{N}$. We assume that the random variables $Y^{(j)}_i$ are relatively stable, thus $\max_{i\leq N}\frac{Y^{(j)}_i}{a_N^{(j)}}\LimitP 1$, as $N\to\infty$. Then
\begin{align*}
\max_{i\leq N}\left(\sum_{j=1}^k\frac{Y_i^{(j)}}{a_N^{(j)}}\right)\LimitP \sup_{(u^{(j)},j\leq k)}\left\{\sum_{j=1}^ku^{(j)}\Bigg|\sum_{j=1}^kh^{(j)}(u^{(j)})\leq 1,u^{(j)}\leq 1 \forall j\leq k\right\}\text{ as $N\to\infty$}.
\end{align*}
\end{lemma}
\proof{Proof}
First of all, let us choose $u^{(1)},\ldots,u^{(k)}$ such that $u^{(j)}\leq 1$ for all $j$. It is a well-known result that
\begin{align*}
\probability*{\cup_{i=1}^N\cap_{j=1}^k\bigg\{Y^{(j)}_i\geq u^{(j)} a^{(j)}_N\bigg\}}\LimitN 1\iff N\probability*{\cap_{j=1}^k\bigg\{Y^{(j)}_i\geq u^{(j)} a^{(j)}_N\bigg\}}\LimitN \infty.
\end{align*}
From this, it follows that
\begin{align*}
\log N +\sum_{j=1}^k\log\left(\probability*{Y^{(j)}_i\geq u^{(j)} a^{(j)}_N}\right)\LimitN \infty.
\end{align*}
This is the case when
\begin{align*}
\limsup_{N\to\infty}\left(\sum_{j=1}^k\frac{-\log\left(\probability*{Y^{(j)}_i\geq u^{(j)} a^{(j)}_N}\right)}{\log N}\right)< 1.
\end{align*}
Similarly,
\begin{align*}
\liminf_{N\to\infty}&\left(\sum_{j=1}^k\frac{-\log\left(\probability*{Y^{(j)}_i\geq u^{(j)} a^{(j)}_N}\right)}{\log N}\right)>1
\Rightarrow \probability*{\cup_{i=1}^N\cap_{j=1}^k\bigg\{Y^{(j)}_i\geq u^{(j)} a^{(j)}_N\bigg\}}\LimitN 0.
\end{align*}
Because of the fact that we have $\probability*{Y^{(j)}_i\geq a^{(j)}_N}=\frac{1}{N}$, we can conclude that
\begin{align*}
&\lim_{N\to\infty}\left(\sum_{j=1}^k\frac{-\log\left(\probability*{Y^{(j)}_i\geq u^{(j)} a^{(j)}_N}\right)}{\log N}\right)
=\lim_{N\to\infty}\left(\sum_{j=1}^k\frac{-\log\left(\probability*{Y^{(j)}_i\geq u^{(j)} a^{(j)}_N}\right)}{-\log\left(\probability*{Y^{(j)}_i\geq a^{(j)}_N}\right)}\right)
=\sum_{j=1}^kh^{(j)}(u^{(j)}).
\end{align*}
Let us now call 
\begin{align*}
c^{\star}=\sup_{(u^{(j)},j\leq k)}\left\{\sum_{j=1}^ku^{(j)}\Bigg|\sum_{j=1}^kh^{(j)}(u^{(j)})\leq 1,u^{(j)}\leq 1 \forall j\leq k\right\},
\end{align*}
and let $\epsilon>0$ be small. Then, we distinguish two scenarios. First of all, we consider the case that $|\{1,\ldots,k|h^{(j)}(u^{(j)})=\mathbbm{1}\left(u^{(j)}>0\right)\}|\leq k-2$. Then, there exists a sequence $(u^{(1)}_\epsilon,\ldots,u^{(k)}_\epsilon)$ such that $\sum_{j=1}^ku^{(j)}_\epsilon=c^{\star}-\epsilon$, and $\sum_{j=1}^kh^{(j)}(u^{(j)}_\epsilon)<1$. Therefore,
\begin{align*}
\probability*{\max_{i\leq N}\left(\sum_{j=1}^k\frac{Y_i^{(j)}}{a_N^{(j)}}\right)>c^{\star}-\epsilon}>\probability*{\cup_{i=1}^N\cap_{j=1}^k\bigg\{Y^{(j)}_i\geq u_{\epsilon}^{(j)} a^{(j)}_N\bigg\}}\LimitN 1.
\end{align*}
If $|\{1,\ldots,k|h^{(j)}(u^{(j)})=\mathbbm{1}\left(u^{(j)}>0\right)\}|\geq k-1$, we know that $c^{\star}=1$, and we know that
\begin{align*}
\max_{i\leq N}\left(\sum_{j=1}^k\frac{Y_i^{(j)}}{a_N^{(j)}}\right)\geq_{st.} \max_{i\leq N}\left(\frac{Y_i^{(1)}}{a_N^{(1)}}\right)+\sum_{j=2}^k\frac{Y_i^{(j)}}{a_N^{(j)}}\LimitP 1\text{ as $N\to\infty$}.
\end{align*}
Thus, at this moment we can conclude that the limit cannot be smaller than $c^{\star}$. To prove that 
\begin{align}\label{eq: upper bound extreme values sum of pairs}
\probability*{\max_{i\leq N}\left(\sum_{j=1}^k\frac{Y_i^{(j)}}{a_N^{(j)}}\right)>c^{\star}+\epsilon}\LimitN 0,
\end{align}
we first observe that the boundary is given by $\{(u^{(j)},j\leq k)|\sum_{j=1}^ku^{(j)}=c^{\star}+\epsilon\}$. We already know that $N\probability*{Y_i^{(j)}>u^{(j)}a_N}\LimitN 0$, for $u^{(j)}>1$. Hence,  
\begin{align*}
\limsup_{N\to\infty}\probability*{\max_{i\leq N}\left(\sum_{j=1}^k\frac{Y_i^{(j)}}{a_N^{(j)}}\right)>c^{\star}+\epsilon}>0,
\end{align*}
means that there are limiting points in the set $\{(u^{(j)},j\leq k)|\sum_{j=1}^ku^{(j)}=c^{\star}+\epsilon,u^{(j)}\leq 1 \forall j\}$. However, we know that for all $(u^{(j)},j\leq k)$ with $c^{\star}<\sum_{j=1}^ku^{(j)}<c^{\star}+\epsilon$ that $N\probability{\cap_{j=1}^k\{Y^{(j)}_i\geq u_{\epsilon}^{(j)} a^{(j)}_N\}}\LimitN 0$. Thus, we know that there are no limiting points in the positive quadrants with starting points $(u^{(1)},\ldots,u^{(k)})$ with $c^{\star}<\sum_{j=1}^ku^{(j)}<c^{\star}+\epsilon$. The union of a finite number of quadrants covers the set $\{(u^{(j)},j\leq k)|\sum_{j=1}^ku^{(j)}=c^{\star}+\epsilon,u^{(j)}\leq 1 \forall j\}$. For example, in the case that $k=2$,
\begin{align*}
\{(u,v)|u+v\geq c^{\star}+\epsilon,u\in[0,1],v\in[0,1]\}
\subset \cup_{m=1}^{\lceil{\frac{2}{\epsilon}-\frac{1}{2}\rceil}+1}\left\{(u,v)\bigg| u\geq c^{\star}-1+\frac{m\epsilon}{2},v\geq 1+\frac{\epsilon}{4}-\frac{m\epsilon}{2}\right\}.
\end{align*}
For $k>2$, an analogous proof can be given.
Hence, the limit in \eqref{eq: upper bound extreme values sum of pairs} and the lemma follows.
\hfill\rlap{\hspace*{-2em}\Halmos}
\endproof
\section[Proof of Lemmas \ref{lem: upper bound |Qs-Qt|}, \ref{lem: distance cdf normal}, \ref{lem: convergence moments 1 2 4} and \ref{lem: approximation starting point convergence}]{Proof of Lemmas \ref{lem: upper bound |Qs-Qt|}, \ref{lem: distance cdf normal}, \ref{lem: convergence moments 1 2 4}, and \ref{lem: approximation starting point convergence}.}\label{app: proof lemma convergence fourth moment}

\proof{Proof of Lemma \ref{lem: upper bound |Qs-Qt|}.}

We take $s>t>0$. We write $t^{(N)}=tN^3\log N$, $s^{(N)}=sN^3\log N$, etc. 
We first prove that for $s^{(N)}>t^{(N)}$, the following upper bound holds:
\begin{align}\label{eq: upper bound Qs-Qt}
\FluidProcessThree{s}-\FluidProcessThree{t}\leq &\max_{i\leq N}\left|\TildeQueue{i}{s^{(N)}}-\TildeQueue{i}{t^{(N)}}\right|\nonumber\\
+&\max_{i\leq N}\sup_{t^{(N)}\leq r\leq s^{(N)}}\left(\TildeQueue{i}{t^{(N)}}-\TildeQueue{i}{r}\right).
\end{align}
Due to the defined auxiliary processes in Definition \ref{def: tilde arrival service}, we can write the maximum queue length in terms of $\TildeQueue{i}{}$ as in Equation \eqref{eq: Q expressed in tilde Q}. Similarly, we can rewrite $\FluidProcessThree{s}-\FluidProcessThree{t}$ as
\begin{align*}
&\max_{i\leq N}\sup_{0\leq r\leq s^{(N)}}\left(\TildeQueue{i}{s^{(N)}}
-\TildeQueue{i}{r}\right)
-\max_{i\leq N}\sup_{0\leq u\leq t^{(N)}}\left(\TildeQueue{i}{t^{(N)}}-\TildeQueue{i}{u}\right)\nonumber\\
= & \max_{i\leq N}\left[\TildeQueue{i}{s^{(N)}}-\TildeQueue{i}{t^{(N)}}+\sup_{0\leq r\leq s^{(N)}}\left(\TildeQueue{i}{t^{(N)}}-\TildeQueue{i}{r}\right)\right]\nonumber\\
-&\max_{i\leq N}\sup_{0\leq u\leq t^{(N)}}\left(\TildeQueue{i}{t^{(N)}}-\TildeQueue{i}{u}\right).
\end{align*}
Now, the following upper bounds for $\FluidProcessThree{s}-\FluidProcessThree{t}$ hold:
\begin{align*}
&\FluidProcessThree{s}-\FluidProcessThree{t}\nonumber\\
\leq & \max_{i\leq N}\left(\TildeQueue{i}{s^{(N)}}-\TildeQueue{i}{t^{(N)}}\right)+\max_{i\leq N}\sup_{0\leq r\leq s^{(N)}}\left(\TildeQueue{i}{t^{(N)}}-\TildeQueue{i}{r}\right)\nonumber\\
-&\max_{i\leq N}\sup_{0\leq u\leq t^{(N)}}\left(\TildeQueue{i}{t^{(N)}}-\TildeQueue{i}{u}\right)\nonumber\\
\leq &\max_{i\leq N}\left(\TildeQueue{i}{s^{(N)}}-\TildeQueue{i}{t^{(N)}}\right)\nonumber\\
+&\max_{i\leq N}\left[\sup_{0\leq r\leq s^{(N)}}\left(\TildeQueue{i}{t^{(N)}}-\TildeQueue{i}{r}\right)-\sup_{0\leq u\leq t^{(N)}}\left(\TildeQueue{i}{t^{(N)}}-\TildeQueue{i}{u}\right)\right].
\end{align*}
Observe that both $\sup_{0\leq r\leq s^{(N)}}\left(\TildeQueue{i}{t^{(N)}}-\TildeQueue{i}{r}\right)$ and 
$\sup_{0\leq u\leq t^{(N)}}\left(\TildeQueue{i}{t^{(N)}}-\TildeQueue{i}{u}\right)$
are non-negative random variables. Furthermore, 
\begin{align*}
&\sup_{0\leq r\leq s^{(N)}}\left(\TildeQueue{i}{t^{(N)}}-\TildeQueue{i}{r}\right)-\sup_{0\leq u\leq t^{(N)}}\left(\TildeQueue{i}{t^{(N)}}-\TildeQueue{i}{u}\right)\nonumber\\
\leq &\sup_{t^{(N)}\leq r\leq s^{(N)}}\left(\TildeQueue{i}{t^{(N)}}-\TildeQueue{i}{r}\right).
\end{align*}
Now, we can conclude that
\begin{align*}
&\FluidProcessThree{s}-\FluidProcessThree{t}\nonumber\\
\leq &\max_{i\leq N}\left(\TildeQueue{i}{s^{(N)}}-\TildeQueue{i}{t^{(N)}}\right)\nonumber\\
+&\max_{i\leq N}\left[\sup_{0\leq r\leq s^{(N)}}\left(\TildeQueue{i}{t^{(N)}}-\TildeQueue{i}{r}\right)-\sup_{0\leq u\leq t^{(N)}}\left(\TildeQueue{i}{t^{(N)}}-\TildeQueue{i}{u}\right)\right]\nonumber\\
\leq & \max_{i\leq N}\left|\TildeQueue{i}{s^{(N)}}-\TildeQueue{i}{t^{(N)}}\right|+\max_{i\leq N}\sup_{t^{(N)}\leq r\leq s^{(N)}}\left(\TildeQueue{i}{t^{(N)}}-\TildeQueue{i}{r}\right),
\end{align*}
and hence the inequality in Equation \eqref{eq: upper bound Qs-Qt} is satisfied. We can similarly deduce the lower bound
\begin{align}\label{eq: lower bound Qs-Qt}
\FluidProcessThree{s}-\FluidProcessThree{t}\geq -\max_{i\leq N}\left|\TildeQueue{i}{t^{(N)}}-\TildeQueue{i}{s^{(N)}}\right|.
\end{align}
To show this, we write 
\begin{align*}
&\FluidProcessThree{s}-\FluidProcessThree{t}\nonumber\\
=&\max_{i\leq N}\sup_{0\leq r\leq s^{(N)}}\left(\TildeQueue{i}{s^{(N)}}-\TildeQueue{i}{r}\right)
-\max_{i\leq N}\sup_{0\leq u\leq t^{(N)}}\left(\TildeQueue{i}{t^{(N)}}-\TildeQueue{i}{u}\right)\nonumber\\
= &\max_{i\leq N}\sup_{0\leq r\leq s^{(N)}}\left(\TildeQueue{i}{s^{(N)}}-\TildeQueue{i}{r}\right)\nonumber\\
-&\max_{i\leq N}\left[\TildeQueue{i}{t^{(N)}}-\TildeQueue{i}{s^{(N)}}+\sup_{0\leq u\leq t^{(N)}}\left(\TildeQueue{i}{s^{(N)}}-\TildeQueue{i}{u}\right)\right]\nonumber\\
\geq &\max_{i\leq N}\sup_{0\leq r\leq s^{(N)}}\left(\TildeQueue{i}{s^{(N)}}-\TildeQueue{i}{r}\right)\nonumber\\
-&\max_{i\leq N}\left(\TildeQueue{i}{t^{(N)}}-\TildeQueue{i}{s^{(N)}}\right)
-\max_{i\leq N}\sup_{0\leq u\leq t^{(N)}}\left(\TildeQueue{i}{s^{(N)}}-\TildeQueue{i}{u}\right).
\end{align*}
Observe that
\begin{align*}
\sup_{0\leq r\leq s^{(N)}}\left(\TildeQueue{i}{s^{(N)}}-\TildeQueue{i}{r}\right)\geq \sup_{0\leq u\leq t^{(N)}}\left(\TildeQueue{i}{s^{(N)}}-\TildeQueue{i}{u}\right),
\end{align*}
because $s^{(N)}>t^{(N)}$, so on the left side of the inequality, the supremum is taken over a larger interval than on the right side of the inequality. From this we can conclude that
\begin{align*}
&\FluidProcessThree{s}-\FluidProcessThree{t}\nonumber\\
\geq &\max_{i\leq N}\sup_{0\leq r\leq s^{(N)}}\left(\TildeQueue{i}{s^{(N)}}-\TildeQueue{i}{r}\right)\nonumber\\
-&\max_{i\leq N}\left(\TildeQueue{i}{t^{(N)}}-\TildeQueue{i}{s^{(N)}}\right)
-\max_{i\leq N}\sup_{0\leq u\leq t^{(N)}}\left(\TildeQueue{i}{s^{(N)}}-\TildeQueue{i}{u}\right)\nonumber\\
\geq & -\max_{i\leq N}\left(\TildeQueue{i}{t^{(N)}}-\TildeQueue{i}{s^{(N)}}\right)
\geq -\max_{i\leq N}\left|\TildeQueue{i}{t^{(N)}}-\TildeQueue{i}{s^{(N)}}\right|,
\end{align*}
and indeed \eqref{eq: lower bound Qs-Qt} holds. Combining \eqref{eq: upper bound Qs-Qt} and \eqref{eq: lower bound Qs-Qt} gives
\begin{align*}
&\left|\FluidProcessThree{s}-\FluidProcessThree{t}\right|\nonumber\\
\leq& \max_{i\leq N}\left|\TildeQueue{i}{s^{(N)}}-\TildeQueue{i}{t^{(N)}}\right|+\max_{i\leq N}\sup_{t^{(N)}\leq r\leq s^{(N)}}\left(\TildeQueue{i}{t^{(N)}}-\TildeQueue{i}{r}\right).
\end{align*}
Thus,
\begin{align}\label{eq: inequality 1}
\sup_{t^{(N)}\leq s \leq t^{(N)}+\delta^{(N)}}\left|\frac{\MaxQueueLength{s}}{N\log N}-\FluidProcessThree{t}\right|
\leq& \sup_{t^{(N)}\leq s\leq t^{(N)}+\delta^{(N)}}\max_{i\leq N}\left|\TildeQueue{i}{s}-\TildeQueue{i}{t^{(N)}}\right|\nonumber\\
+&\sup_{t^{(N)}\leq s\leq t^{(N)}+\delta^{(N)}}\max_{i\leq N}\left(\TildeQueue{i}{t^{(N)}}-\TildeQueue{i}{s}\right).
\end{align}
Since both $\sup_{t^{(N)}\leq s \leq t^{(N)}+\delta^{(N)}}\left(\TildeQueue{i}{t^{(N)}}-\TildeQueue{i}{s}\right)$ and \\ $\sup_{t^{(N)}\leq s \leq t^{(N)}+\delta^{(N)}}\left(\TildeQueue{i}{s}-\TildeQueue{i}{t^{(N)}}\right)$ are non-negative random variables, we have that
\begin{align}\label{eq: inequality 2}
\sup_{t^{(N)}\leq s\leq t^{(N)}+\delta^{(N)}}\max_{i\leq N}\left|\TildeQueue{i}{s}-\TildeQueue{i}{t^{(N)}}\right|\leq &\sup_{t^{(N)}\leq s\leq t^{(N)}+\delta^{(N)}}\max_{i\leq N}\left(\TildeQueue{i}{s}-\TildeQueue{i}{t^{(N)}}\right)\nonumber\\
+&\sup_{t^{(N)}\leq s\leq t^{(N)}+\delta^{(N)}}\max_{i\leq N}\left(\TildeQueue{i}{t^{(N)}}-\TildeQueue{i}{s}\right).
\end{align}
Combining the inequalities in \eqref{eq: inequality 1} and \eqref{eq: inequality 2} gives us the desired result.
\hfill\rlap{\hspace*{-2em}\Halmos}
\endproof
\text{ }
\proof{Proof of Lemma \ref{lem: distance cdf normal}.}
$\TildeService{i}{n}$ is a sum of independent and identically distributed random variables with $\expect*{\pm\TildeService{i}{1}}=0$,
and $\text{Var}\left(\pm\TildeService{i}{1}\right)=\left(1-\frac{\alpha}{N}\right)\frac{\alpha}{N^3}.$
So, $\pm\MaximumProcess{i}{t}=\frac{\pm\TildeService{i}{tN^3\log N}\sqrt{tN^3\log N}}{\sqrt{\alpha t(1-\frac{\alpha}{N})\log N\lfloor tN^3\log N\rfloor}}$ has mean 0 and variance 1, and satisfies the central limit theorem.
From \cite{michel1976constant} it follows that for all $y$,
\begin{align*}
\left|\probability*{\pm\MaximumProcess{i}{t}<y}-\Phi(y)\right|
\leq & C\frac{1}{\sqrt{\left\lfloor{tN^3\log N}\right\rfloor}}\expect*{\left|\frac{\pm\TildeService{i}{1}}{\sqrt{\alpha t(1-\frac{\alpha}{N})\log N}}\sqrt{tN^3\log N}\right|^3}\frac{1}{1+|y|^3}.
\end{align*}
Observe that for $N$ large enough and $0<\epsilon<t$, $\left\lfloor{tN^3\log N}\right\rfloor>(t-\epsilon)N^3\log N$. We also have that
\begin{align*}
&\expect*{\left|\frac{\pm\TildeService{i}{1}}{\sqrt{\alpha t(1-\frac{\alpha}{N})\log N}}\sqrt{tN^3\log N}\right|^3}\nonumber\\
=&\frac{N^4\sqrt{N}}{\alpha (1-\frac{\alpha}{N})\sqrt{\alpha (1-\frac{\alpha}{N})}N^3}\left(\left(1-\frac{\alpha}{N}\right)^3\frac{\alpha}{N}+\frac{\alpha^3}{N^3}\left(1-\frac{\alpha}{N}\right)\right)
\leq 2\sqrt{N}\frac{(1+\alpha^2)}{\sqrt{\alpha}},
\end{align*}
which holds for $N>\max(1,2\alpha)$. Thus, the statement of Lemma \ref{lem: distance cdf normal} follows for $N$ large enough, with $c_t=2C\frac{(1+\alpha^2)}{\sqrt{\alpha (t-\epsilon)}}$.
\hfill\rlap{\hspace*{-2em}\Halmos}
\endproof
\text{ }
\proof{Proof of Lemma \ref{lem: convergence moments 1 2 4}.}
We have
\begin{align}
&\expect*{\max\left(0,\frac{\max_{i\leq N}\pm\TildeService{i}{t N^3\log N}}{\log N}\right)^{5/2}}
=\int_{0}^{\infty}\probability*{\frac{\max_{i\leq N}\pm\TildeService{i}{t N^3\log N}}{\log N}>x^{2/5}}dx\nonumber\\
=&\int_{0}^{\infty}\probability*{\max_{i\leq N}\pm\MaximumProcess{i}{t}>x^{2/5}\frac{\log N}{\sqrt{\alpha t\left(1-\frac{\alpha}{N}\right)\log N}}\frac{\sqrt{tN^3\log N}}{\sqrt{\lfloor{tN^3\log N\rfloor}}}}dx\nonumber\\
=&\int_{0}^{\infty}1-\probability*{\pm\MaximumProcess{i}{t}<\frac{x^{2/5}\sqrt{N^3}\log N}{\sqrt{\alpha \left(1-\frac{\alpha}{N}\right)\lfloor{tN^3\log N\rfloor}}}}^Ndx\nonumber\\
\leq&\int_{0}^{\infty}1
-\Vast(\Phi\left(\frac{x^{2/5}\sqrt{N^3}\log N}{\sqrt{\alpha \left(1-\frac{\alpha}{N}\right)\lfloor{tN^3\log N\rfloor}}}\right)
-\frac{c_t}{N\sqrt{\log N}}\frac{1}{1+\left(\oldfrac{x^{2/5}\sqrt{N^3}\log N}{\sqrt{\alpha \left(1-\frac{\alpha}{N}\right)\lfloor{tN^3\log N\rfloor}}}\right)^3}\Vast)^Ndx\nonumber\\
\leq&\int_{0}^{\infty}-\Phi\left(\frac{x^{2/5}\sqrt{N^3}\log N}{\sqrt{\alpha \left(1-\frac{\alpha}{N}\right)\lfloor{tN^3\log N\rfloor}}}\right)^N
+\left(1+\frac{c_t}{N\sqrt{\log N}}\frac{1}{1+\left(\oldfrac{x^{2/5}\sqrt{N^3}\log N}{\sqrt{\alpha \left(1-\frac{\alpha}{N}\right)\lfloor{tN^3\log N\rfloor}}}\right)^3}\right)^Ndx\nonumber\\
=& \expect*{\max\left(0,\frac{\sqrt{\alpha t(1-\frac{\alpha}{N})}\max_{i\leq N}X_i}{\sqrt{\log N}}\frac{\sqrt{\lfloor{tN^3\log N\rfloor}}}{\sqrt{tN^3\log N}}\right)^{5/2}}\label{subeq: pickands expectation}\\
+&\int_{0}^{\infty}-1+\left(1+\frac{c_t}{N\sqrt{\log N}}\frac{1}{1+\left(\oldfrac{x^{2/5}\sqrt{N^3}\log N}{\sqrt{\alpha \left(1-\frac{\alpha}{N}\right)\lfloor{tN^3\log N\rfloor}}}\right)^3}\right)^Ndx,\label{subeq: pickands expectation error}
\end{align}
with $X_i$ standard normally distributed.
By Pickands \cite[Thm.~3.2, p.\ 888]{pickands1968moment}, we know that the expectation in \eqref{subeq: pickands expectation} converges to $(2\alpha t)^{5/4}$. Furthermore, the term in \eqref{subeq: pickands expectation error} is upper bounded by
\begin{align}\label{eq: integral bound michel}
\int_{0}^{\infty}-1+\text{exp}\left(\frac{c_t}{\sqrt{\log N}\left(1+\left(\oldfrac{x^{2/5}\sqrt{N^3}\log N}{\sqrt{\alpha \left(1-\frac{\alpha}{N}\right)\lfloor{tN^3\log N\rfloor}}}\right)^3\right)}\right)dx.
\end{align}
We substitute $y=1\bigg/\left(1+\left(\oldfrac{x^{2/5}\sqrt{N^3}\log N}{\sqrt{\alpha\left(1-\frac{\alpha}{N}\right) \lfloor{tN^3\log N\rfloor}}}\right)^3\right)$, then the term in \eqref{eq: integral bound michel} can be rewritten as 
\begin{align*}
\left(\frac{\lfloor{tN^3\log N\rfloor}}{N^3\log N(\log N)}\right)^{5/4}\int_{0}^1\frac{5 (\sqrt{\alpha\left(1-\frac{\alpha}{N}\right) })^{5/2}}{6 (1-y)^{1/6} y^{11/6}}\left(-1+e^{\oldfrac{c_t}{\sqrt{\log N}} y}\right)dy\LimitN 0.
\end{align*}
The lemma follows.
\hfill\rlap{\hspace*{-2em}\Halmos}
\endproof
\text{ }
\proof{Proof of Lemma \ref{lem: approximation starting point convergence}}
In order to prove that
\begin{align*}
\max_{i\leq N}\left(\frac{\sqrt{\alpha t}X_i}{\sqrt{\log N}}+\frac{\QueueZero{i}}{N\log N}\right)\LimitP g(t,q(0)),
\end{align*}
we first observe that, from the definition of $\QueueZero{i}$ in Theorem \ref{thm: fluid limit}, it is easy to see that
\begin{align*}
&\left|\max_{i\leq N}\left(\frac{\sqrt{\alpha t}X_i}{\sqrt{\log N}}+\frac{\QueueZero{i}}{N\log N}\right)-\max_{i\leq N}\left(\frac{\sqrt{\alpha t}X_i}{\sqrt{\log N}}+\frac{r_N\StartRVIndep{i}}{N\log N}\right)\right|\leq\max_{i\leq N}\frac{\StartRVTriangDep{i}}{N\log N}+\frac{1}{N\log N}\LimitP 0
\end{align*}
as $N\to\infty$. Thus, from this it follows that 
\begin{align*}
\max_{i\leq N}\left(\frac{\sqrt{\alpha t}X_i}{\sqrt{\log N}}+\frac{\QueueZero{i}}{N\log N}\right)\LimitP g(t,q(0))\iff \max_{i\leq N}\left(\frac{\sqrt{\alpha t}X_i+\oldfrac{r_N\StartRVIndep{i}}{N\sqrt{\log N}}}{\sqrt{\log N}}\right)\LimitP g(t,q(0))\text{ as $N\to\infty$}.
\end{align*}
Let us first consider that $\StartRVIndep{i}$ satisfies Assumption \ref{ass: 4}, thus $\StartRVIndep{i}$ has a finite right endpoint. Theorem \ref{thm: fluid limit} says that when $\StartRVIndep{i}$ has a finite right endpoint, that $g(t,q(0))=\sqrt{2\alpha t}+q(0)$. To prove this, first observe that $g(t,q(0))\leq \sqrt{2\alpha t}+q(0)$ because $\max_{i\leq N}\frac{\sqrt{\alpha t}X_i}{\sqrt{\log N}}\LimitP \sqrt{2\alpha t}$ and $\frac{\MaxQueueLength{0}}{(N\log N)}\LimitP q(0)$. Hence, the only thing we need to establish is that for all $\gamma<\sqrt{2\alpha t}+q(0)$,
\begin{align*}
N\probability*{\sqrt{\alpha t}X_i+\frac{r_N\StartRVIndep{i}}{N\sqrt{\log N}}\geq \gamma\sqrt{\log N}}\LimitN \infty.
\end{align*}
When $\gamma<\sqrt{2\alpha t}$, this is obvious, because $\StartRVIndep{i}>0$, and $\max_{i\leq N}\frac{\sqrt{\alpha t}X_i}{\sqrt{\log N}}\LimitP \sqrt{2\alpha t}$. So, let us assume that $\sqrt{2\alpha t}\leq\gamma<\sqrt{2\alpha t}+q(0)$. Because $\StartRVIndep{i}$ has a finite right endpoint, $\frac{r_N}{(N\sqrt{\log N})}=\sqrt{\log N}$. By convolution, we have that
\begin{align*}
&N\probability*{\sqrt{\alpha t}X_i+\sqrt{\log N}\StartRVIndep{i}\geq \gamma\sqrt{\log N}}\nonumber\\
=& N\probability*{\sqrt{\alpha t}X_i\geq \gamma\sqrt{\log N}}
+N\int_{-\infty}^{\gamma\sqrt{\log N}}\probability*{\sqrt{\log N}\StartRVIndep{i}>\gamma\sqrt{\log N}-z}\frac{e^{\frac{-z^2}{(2\alpha t)}}}{\sqrt{2\alpha t\pi}} dz\nonumber\\
\geq&N\int_{-\infty}^{\gamma}\probability*{\StartRVIndep{i}>\gamma-v}\frac{N^{\frac{-v^2}{(2\alpha t)}}}{\sqrt{2\alpha t\pi}}\sqrt{\log N} dv
=\int_{\gamma -q(0)}^{\gamma}\probability*{\StartRVIndep{i}>\gamma-v}\frac{N^{1-\frac{v^2}{(2\alpha t)}}}{\sqrt{2\alpha t\pi}}\sqrt{\log N} dv.
\end{align*}
From this it follows, that when $1-\frac{v^2}{(2\alpha t)}>0$, this integral converges to $\infty$. We chose $\sqrt{2\alpha t}\leq\gamma<\sqrt{2\alpha t}+q(0)$, thus the lower bound $\gamma -q(0)$ in the integral is smaller than $\sqrt{2\alpha t}$ and hence this integral converges to $\infty$. Thus $g(t,q(0))=\sqrt{2\alpha t}+q(0)$.

Let us now consider the scenario described in Assumption \ref{ass: 5}.
Then $g(t,q(0))$ satisfies the limit given in \eqref{eq: limit scenario 2}.
We have the straightforward limit result that for standard normally distributed $X_i$, $\lim_{t\to\infty}\frac{-\log(\probability*{X_i\geq ut})}{-\log(\probability*{X_i\geq t})}=u^2$. Furthermore, following the assumptions on $\StartRVIndep{i}$ in Theorem \ref{thm: fluid limit}, we know that $\lim_{t\to\infty}\frac{-\log(\probability*{\StartRVIndep{i}\geq vt})}{-\log(\probability*{\StartRVIndep{i}\geq t})}=h(v)$. Thus from Lemma \ref{lem: extreme value convergence}, we know that for sequences $(a_N,N\geq 1), (b_N,N\geq 1)$ with $\probability*{X_i\geq a_N}=\probability*{\StartRVIndep{i}\geq b_N}=\frac{1}{N}$, that
\begin{align*}
\max_{i\leq N}\bigg(\frac{X_i}{a_N}+\frac{\StartRVIndep{i}}{b_N}\bigg)\LimitP \sup_{(u,v)}\{u+v|u^2+h(v)\leq 1,0\leq u\leq 1, 0\leq v\leq 1\} \text{ as $N\to\infty$}.
\end{align*}
Now, we can use this result to prove that $\max_{i\leq N}\left(\frac{\sqrt{\alpha t}X_i}{\sqrt{\log N}}+\frac{r_N \StartRVIndep{i}}{(N\log N)}\right)$ converges to the limit in \eqref{eq: limit scenario 2}.
We first observe that
\begin{align*}
\max_{i\leq N}\bigg(\frac{\sqrt{\alpha t}X_i}{\sqrt{\log N}}+\frac{r_N \StartRVIndep{i}}{N\log N}\bigg)=\max_{i\leq N}\bigg(\sqrt{2\alpha t}\frac{X_i}{\sqrt{2\log N}}+q(0)\frac{r_N \StartRVIndep{i}}{q(0)N\log N}\bigg).
\end{align*}
We have that $\frac{a_N}{\sqrt{2\log N}}\LimitN 1$, because
$\max_{i\leq N}\frac{X_i}{a_N}\LimitP 1$,
and $\max_{i\leq N}\frac{X_i}{\sqrt{2\log N}}\LimitP 1\text{ as $N\to\infty$}$.
Analogously, $\frac{b_Nq(0)N\log N}{r_N}\LimitN 1.$
Thus,
\begin{align*}
\left|\max_{i\leq N}\bigg(\sqrt{2\alpha t}\frac{X_i}{a_N}+q(0)\frac{\StartRVIndep{i}}{b_N}\bigg)-\max_{i\leq N}\bigg(\frac{\sqrt{\alpha t}X_i}{\sqrt{\log N}}+\frac{r_N \StartRVIndep{i}}{N\log N}\bigg)\right|\LimitP 0 \text{ as $N\to\infty$}.
\end{align*}
With an analogous proof as before, $\max_{i\leq N}\left(\sqrt{2\alpha t}\frac{X_i}{a_N}+q(0)\frac{\StartRVIndep{i}}{b_N}\right)$ converges to the limit in \eqref{eq: limit scenario 2}.
\hfill\rlap{\hspace*{-2em}\Halmos}
\endproof
\section[Notation]{Notation.}\label{app: Notation}
\begin{itemize}
\item $N$: the number of servers.
\item $\Arrival{n}$: the number of arrivals up to time $\lfloor n\rfloor$.
\item $\ArrivalBernoulli{n}$: Bernoulli random variable indicating a potential arrival at time $n\in\mathbb{N}$.
\item $\Service{i}{n}$: the number of finished services of server $i$ up to time $\lfloor n\rfloor$.
\item $\ServiceBernoulli{i}{n}$: Bernoulli distributed random variable indicating a potential completed service at server $i$ at time $n\in\mathbb{N}$.
\item $\alpha,\beta$ are system parameters.
\item $\ArrivalProbability$: the arrival probability, $\ArrivalProbability=1-\frac{\alpha}{N}-\frac{\beta}{N^2}$.
\item $\SuccessProbability$: the service probability, $\SuccessProbability=1-\frac{\alpha}{N}$.
\item $\MaxQueueLength{n}$: the maximum queue length at time $\lfloor n\rfloor$.
\item $\QueueZero{i}$: the number of tasks at time 0 at queue $i$, $\QueueZero{i}=\StartRVTriangIndep{i}+\StartRVTriangDep{i}$.
\item $\StartRVTriangIndep{i}$: the independent part of the number of tasks at time 0 at queue $i$, $\StartRVTriangIndep{i}=\lfloor r_N \StartRVIndep{i}\rfloor$.
\item $\StartRVTriangDep{i}$: the dependent part of the number of tasks at time 0 at queue $i$.
\item $\StartRVIndep{i}$: continuously distributed and positive random variable.
\item $r_N$: positive scaling sequence.
\item $h(v)=\lim_{t\to\infty}\frac{-\log \left(\probability*{\StartRVIndep{i}>vt}\right)}{-\log \left(\probability*{\StartRVIndep{i}>t}\right)}$.
\item $q(t)$: fluid limit of the process.
\item $g(t,q(0))$: limit of $\max_{i\leq N}\frac{\big(\Arrival{tN^3\log N}-\Service{i}{tN^3\log N}+\QueueZero{i}\big)}{\big(N\log N\big)}$.
\item $\TildeQueue{i}{n}=\frac{(\TildeArrival{n}+\TildeService{i}{n})}{\log N}$.
\item $\TildeArrival{n}=\frac{\Arrival{n}}{N}-\left(1-\frac{\alpha}{N}\right)\frac{\lfloor n\rfloor}{N}$.
\item $\TildeService{i}{n}=-\frac{\Service{i}{n}}{N}+\left(1-\frac{\alpha}{N}\right)\frac{\lfloor n\rfloor}{N}$.
\item $\MaximumProcess{i}{t}=\frac{\TildeService{i}{{t N^3\log N}}\sqrt{tN^3\log N}}{\big(\sqrt{\alpha t(1-\frac{\alpha}{N})\log N}\sqrt{\lfloor{tN^3\log N\rfloor}}\big)}$.
\item $\ArrivalUpper{n}=\sum_{j=1}^n \ArrivalBernoulliUpper{j}$.
\item $\ArrivalBernoulliUpper{j}$:
\begin{align*}
\ArrivalBernoulliUpper{j}=\left\{\begin{matrix}
\text{ }\text{ }\text{ }\text{ }\text{ }\text{ }\text{ }\text{ }\frac{\alpha}{N}+\frac{\beta}{N^{2}}-\frac{m}{N^2}&\text{ w.p. }&  1-\frac{\alpha}{N}-\frac{\beta}{N^2},&\\ 
-1+\frac{\alpha}{N}+\frac{\beta}{N^{2}}-\frac{m}{N^2}&\text{ w.p. }& \text{ }\text{ }\text{ }\text{ }\text{ }\frac{\alpha}{N}+\frac{\beta}{N^2},&
\end{matrix}\right.
\end{align*}
with $0<m<\beta$.
\item $\ServiceUpper{i}{n}=\sum_{j=1}^n \ServiceBernoulliUpper{i}{j}$.
\item $\ServiceBernoulliUpper{i}{j}$:
\begin{align*}
\ServiceBernoulliUpper{i}{j}=\left\{\begin{matrix}\text{ }
-\frac{\alpha}{N}-\frac{\beta}{N^{2}}+\frac{m}{N^2}&\text{ w.p. }&  1-\frac{\alpha}{N},\\ 
1-\frac{\alpha}{N}-\frac{\beta}{N^{2}}+\frac{m}{N^2}&\text{ w.p. }& \text{ }\text{ }\text{ }\text{ }\text{ }\frac{\alpha}{N}.
\end{matrix}\right.
\end{align*}
\item $\ThetaArrivalUpper$ solves:
\begin{align*}
\mathbb{E}\left[e^{\ThetaArrivalUpper \ArrivalBernoulliUpper{j}}\right]=1.
\end{align*}
\item $E_i^{(u,N)}\sim\text{Exp}\left(\frac{2(\beta-m)}{(\alpha N)}\right)$.
\end{itemize}
\end{APPENDICES}


\pdfbookmark[0]{Acknowledgments}{acknowledgements}
\section*{Acknowledgments.}

This research is  supported by the Netherlands Organisation for Scientific Research through the programmes Grip on Complexity [Schol: 438.16.121], MEERVOUD [Vlasiou: 632.003.002], and Talent VICI [Zwart: 639.033.413]. We thank the referees for prompting us to investigate the system where the number of jobs at time 0 increases with $N$. We thank Dr. Guus Balkema for referring us to literature on the convergence of samples.


\bibliographystyle{informs2014} 

\begin{thebibliography}{26}
\providecommand{\natexlab}[1]{#1}
\providecommand{\url}[1]{\texttt{#1}}
\providecommand{\urlprefix}{URL }

\bibitem[{Anderson et~al.(1997)Anderson, Coles, H{\"u}sler
  et~al.}]{anderson1997maxima}
Anderson CW, Coles SG, H{\"u}sler J, et~al. (1997) Maxima of Poisson-like
  variables and related triangular arrays. \emph{Ann. Appl. Probab.} 7(4):953--971.

\bibitem[{Atar et~al.(2012)Atar, Mandelbaum, \protect\BIBand{}
  Zviran}]{atar2012control}
Atar R, Mandelbaum A, Zviran A (2012) Control of fork-join networks in heavy
  traffic. \emph{Proc. 50th Allerton Conf. Comm.,
Control Comput.} (IEEE, Piscataway, NJ), 823--830.

\bibitem[{Baccelli(1985)}]{baccelli1985two}
Baccelli F (1985) Two parallel queues created by arrivals with two demands: The
  M/G/2 symmetrical case. Technical report RR-0426, INRIA, Montbonnot-Saint-Martin, France.

\bibitem[{Baccelli \protect\BIBand{} Makowski(1989)}]{baccelli1989queueing}
Baccelli F, Makowski AM (1989) Queueing models for systems with synchronization
  constraints. \emph{Proc. of the IEEE} 77(1):138--161.


\bibitem[{Billingsley(1968)}]{billingsley1968convergence}
Billingsley P (1968) \emph{Convergence of probability measures}, 2nd ed. (John Wiley \&
  Sons, New York).

\bibitem[{Brown \protect\BIBand{} Resnick(1977)}]{brown1977extreme}
Brown BM, Resnick SI (1977) Extreme values of independent stochastic processes. \emph{J. Appl. Probab.} 14(4):732--739.


\bibitem[{Davis et~al.(1988)Davis, Mulrow, \protect\BIBand{}
  Resnick}]{davis1988almost}
Davis RA, Mulrow E, Resnick SI (1988) Almost sure limit sets of random samples
  in $\mathbb{R}^d$. \emph{Adv. Appl. Probab.} 20(3):573--599.


\bibitem[{Fisher L (1969)}]{fisher1969limiting}
Fisher L (1969) Limiting sets and convex hulls of samples from product
  measures. \emph{Ann. Math. Statist.} 40(5):1824--1832.

\bibitem[{Flatto \protect\BIBand{} Hahn(1984)}]{flatto1984two}
Flatto L, Hahn S (1984) Two parallel queues created by arrivals with two
  demands I. \emph{SIAM J. Appl. Math.} 44(5):1041--1053.
  
  
  \bibitem[{Harrison(1985)}]{harrison1985brownian}
Harrison JM (1985) \emph{Brownian Motion and Stochastic Flow Systems} (John Wiley \& Sons, New York).
 

\bibitem[{Klein(1988)}]{de1988fredholm}
de Klein SJ (1988) \emph{Fredholm integral equations in queueing analysis}. Ph.D.
  thesis, Rijksuniversiteit Utrecht.

\bibitem[{Ko \protect\BIBand{} Serfozo(2004)}]{ko2004response}
Ko SS, Serfozo RF (2004) Response times in M/M/s fork-join networks.
  \emph{Adv. Appl. Probab.} 36(3):854--871.

\bibitem[{Lu \protect\BIBand{} Pang(2015)}]{lu2015gaussian}
Lu H, Pang G (2015) Gaussian limits for a fork-join network with
  nonexchangeable synchronization in heavy traffic. \emph{Math. Oper. Res.} 41(2):560--595.

\bibitem[{Lu \protect\BIBand{} Pang(2017{\natexlab{a}})}]{lu2017heavy}
Lu H, Pang G (2017{\natexlab{a}}) Heavy-traffic limits for a fork-join network
  in the Halfin-Whitt regime. \emph{Stochastic Systems} 6(2):519--600.

\bibitem[{Lu \protect\BIBand{} Pang(2017{\natexlab{b}})}]{lu2017heavy2}
Lu H, Pang G (2017{\natexlab{b}}) Heavy-traffic limits for an infinite-server
  fork--join queueing system with dependent and disruptive services.
  \emph{Queueing Systems} 85(1-2):67--115.

\bibitem[{Michel(1976)}]{michel1976constant}
Michel R (1976) On the constant in the nonuniform version of the Berry-Ess\'een
  theorem. \emph{Zeitschrift f{\"u}r Wahrscheinlichkeitstheorie und verwandte
  Gebiete} 55(1):109--117.

\bibitem[{Nelson \protect\BIBand{} Tantawi(1988)}]{nelson1988approximate}
Nelson R, Tantawi AN (1988) Approximate analysis of fork/join synchronization
  in parallel queues. \emph{IEEE Trans. Comput.} 37(6):739--743.

\bibitem[{Nguyen(1993)}]{nguyen1993processing}
Nguyen V (1993) Processing networks with parallel and sequential tasks: Heavy
  traffic analysis and Brownian limits. \emph{Ann. Appl. Probab.} 3(1):28--55.

\bibitem[{Nguyen(1994)}]{nguyen1994trouble}
Nguyen V (1994) The trouble with diversity: Fork-join networks with
  heterogeneous customer population. \emph{Ann. Appl. Probab.}
  4(1):1--25.

\bibitem[{Pickands~III(1968)}]{pickands1968moment}
Pickands~III J (1968) Moment convergence of sample extremes. \emph{Ann. Math. Statist.} 39(3):881--889.

\bibitem[{Sigman \protect\BIBand{} Whitt(2011{\natexlab{a}})}]{sigman2011heavy}
Sigman K, Whitt W (2011{\natexlab{a}}) Heavy-traffic limits for nearly
  deterministic queues. \emph{J. Appl. Probab.} 48(3):657--678.

\bibitem[{Sigman \protect\BIBand{}
  Whitt(2011{\natexlab{b}})}]{sigman2011heavystat}
Sigman K, Whitt W (2011{\natexlab{b}}) Heavy-traffic limits for nearly
  deterministic queues: stationary distributions. \emph{Queueing Systems}
  69(2):145.

\bibitem[{Varma(1990)}]{varma1990heavy}
Varma S (1990) \emph{Heavy and light traffic approximations for queues with
  synchronization constraints}. Ph.D. thesis, University of Maryland.

\bibitem[{Wright(1992)}]{wright1992two}
Wright PE (1992) Two parallel processors with coupled inputs. \emph{Adv. Appl. Probab.} 24(4):986--1007.

\end{thebibliography}



\end{document}